\documentclass[preprint,3p,12pt]{elsarticle}

\usepackage{color}
\usepackage{amsmath}
\usepackage{amsfonts}
\usepackage{amssymb}
\usepackage{latexsym}
\usepackage{latexsym}
\usepackage{array}
\usepackage{amsthm}
\usepackage{graphics}
\usepackage{comment}
\usepackage{epsf,graphicx,subfigure}
\usepackage{epsfig}
\usepackage{epstopdf}
\usepackage{tikz}
\usetikzlibrary{arrows,positioning,shapes.geometric}

\usepackage{longtable}
\usepackage{rotating}
\usepackage{booktabs}

\newtheorem{Theorem}{Theorem}
\newtheorem{Proposition}{Proposition}
\newtheorem{Lemma}{Lemma}
\newtheorem{Remark}{Remark}
\newtheorem{Example}{Example}
\newtheorem{Assumption}{Assumption}
\makeatletter
\def\ps@pprintTitle{%
  \let\@oddhead\@empty
  \let\@evenhead\@empty
  \let\@oddfoot\@empty
  \let\@evenfoot\@oddfoot
}
\makeatother

\begin{document}

\begin{frontmatter}

\title{When Scattering Transform Meets Non-Gaussian Random Processes, a Double Scaling Limit Result}

\author[GRLiu]{Gi-Ren Liu}
\author[YCSheu]{Yuan-Chung Sheu}
\author[HTWu]{Hau-Tieng Wu\corref{mycorrespondingauthor}}
\cortext[mycorrespondingauthor]{Hau-Tieng Wu}
\ead{hauwu@math.duke.edu}

\address[GRLiu]{Department of Mathematics, National Chen-Kung University, Tainan, Taiwan}
\address[YCSheu]{Department of Applied Mathematics, National Yang Ming Chiao Tung University, Hsinchu, Taiwan}
\address[HTWu]{Department of Mathematics and Department of Statistical Science, Duke University, Durham, NC, USA}

\begin{abstract}
%Scattering transform is a novel signal processing tool aiming to capture the spirit of widely applied convolutional neural network.
We explore the finite dimensional distributions of the second-order scattering transform of a class of non-Gaussian processes
when all the scaling parameters go to infinity simultaneously.
For frequently used wavelets, we find a coupling rule for the scale parameters of the wavelet transform within the first and second layers
such that the limit exists. The coupling rule is explicitly expressed in terms of the Hurst index of the long range dependent inputs.
More importantly, the spectral density function of the limiting process
can be explicitly expressed in terms of the Hurst index of the long-range dependent input process and the Fourier transform of the mother wavelet. To obtain these results, we first show that the scattering transform of a class of non-Gaussian processes
can be approximated by the second-order scattering transform of Gaussian processes when the scale parameters are sufficiently large.
Secondly, for each approximation process, which is still a non-Gaussian random process, we show that
its doubling scaling limit exists by the moment method.
The long-range dependent non-Gaussian model considered in this work is realized by real data, and assumptions are supported by numerical results.

\end{abstract}

\begin{keyword}
long-range dependent; non-Gaussian processes; wavelet transform;
scattering transform; Wiener-It$\hat{\textup{o}}$ decomposition;
Feynman diagram; double scaling limits.
%% MSC codes here, in the form: \MSC code \sep code
%% or \MSC[2008] code \sep code (2000 is the default)
\MSC[2010] Primary 60G60, 60H05, 62M15; Secondary 35K15.

\end{keyword}

\end{frontmatter}

%% \linenumbers

%% main text
\section{Introduction}\label{sec:introduction}

For the purpose of feature extraction, the second-order scattering transforms (ST) has been applied to various signals, for example, fetal heart rate \cite{chudavcek2013scattering}, brain waves \cite{liu2020diffuse,liu2020hospitals}, respiration \cite{wu2014assess}, marine bioacoustics \cite{balestriero2017linear}, and audio \cite{anden2011multiscale,li2019heart}.
Given a signal $X$, the second-order ST of $X$ is defined through the convolution and modulus operators as follows
\begin{equation*}
U[j_{1},j_{2}]X(t)=\big||X\star\psi_{j_{1}}|\star\psi_{j_{2}}(t)\big|,\ j_{1},j_{2}\in \mathbb{Z},\ t\in\mathbb{R},
\end{equation*}
where
$\psi_{j_{1}}(\cdot)=2^{-j_{1}}\psi(2^{-j_{1}}\cdot)$ and $\psi_{j_{2}}(\cdot)=2^{-j_{2}}\psi(2^{-j_{2}}\cdot)$ are wavelets generated from a selected mother wavelet $\psi$.
In the work \cite{bruna2015intermittent}, $U[j_{1},j_{2}]X$ has been
used as multi-scale measurements of intermittency, i.e., the frequency of occurrence of activity at different scales.

There have been several theoretical supports established for ST. If the mother wavelet and the low-pass filter $\phi_{J}$ satisfy the Littlewood-Paley condition, the ST coefficients $\{U[j_{1},j_{2}]X\star\phi_{J}\mid j_{1},j_{2}\leq J\}$ are approximately invariant to time shifts and stable to small time-warping deformation \cite{mallat2012group} when the integer $J$ is large enough. Similar time invariance and stability to deformation results are discussed in \cite{wiatowski2017mathematical} under a different setting.
The limiting distributions of the second-order ST of random processes
have ever been explored in \cite{bruna2015intermittent} for the fractional Brownian motions
and \cite{liu2020central} for stationary Gaussian processes.
In these works, the first scale parameter $j_{1}$ remains fixed and the second scale parameter $j_{2}$ goes to infinity.
%\begin{figure}[h]
%\centering
%\subfigure[][Awake]
%{\includegraphics[scale=0.4]{EEG_raw_awake.eps}}
%\subfigure[][REM]
%{\includegraphics[scale=0.4]{EEG_raw_REM.eps}}
%\subfigure[][N1]
%{\includegraphics[scale=0.4]{EEG_raw_N1.eps}}
%\subfigure[][N2]
%{\includegraphics[scale=0.4]{EEG_raw_N2.eps}}
%\subfigure[][N3]
%{\includegraphics[scale=0.4]{EEG_raw_N3.eps}}
%\caption{30-s EEG signals during different stages of sleep}
%\label{fig:EEG_raw}
%\end{figure}

\begin{figure}
 \centering
  \includegraphics[scale=0.4]{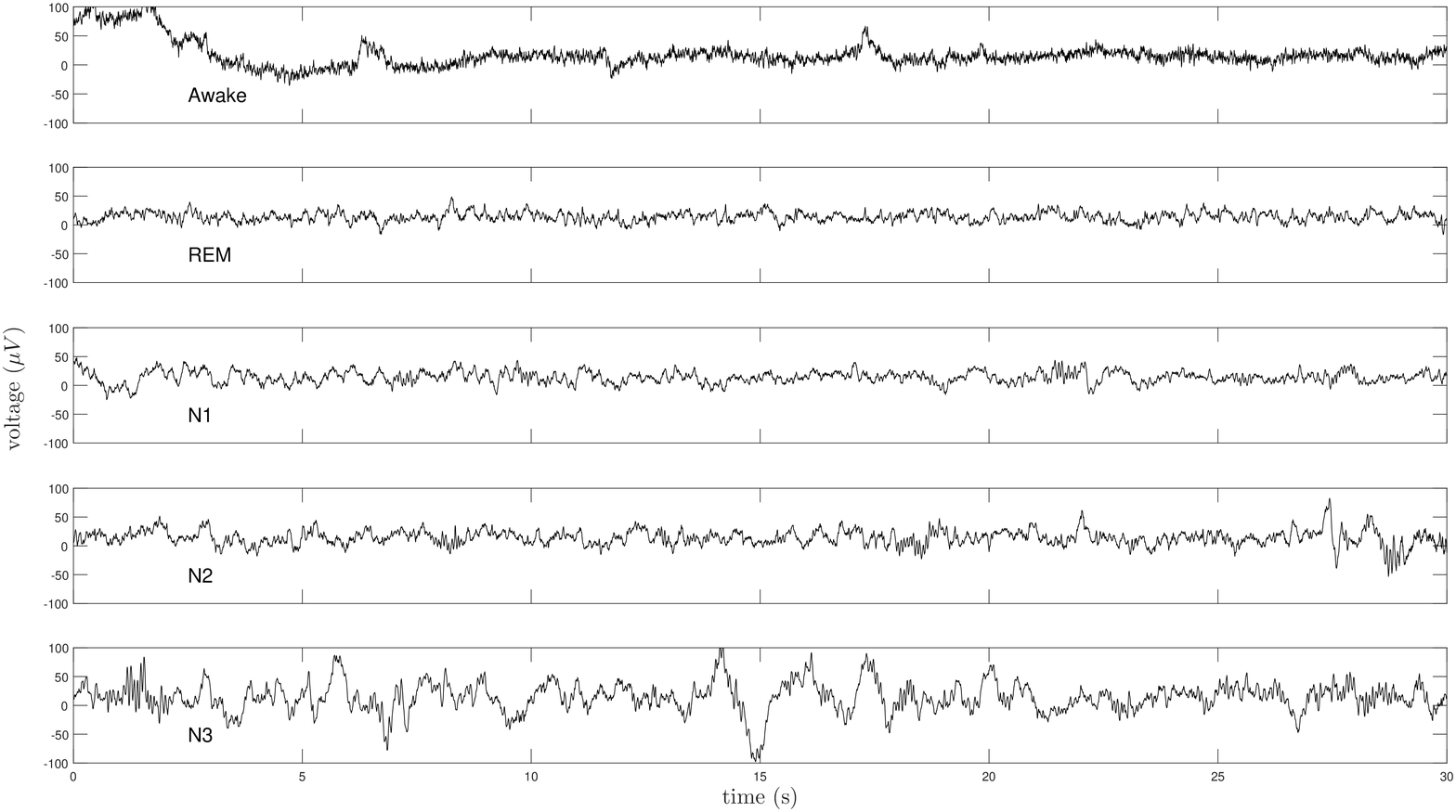}
\caption{30-s EEG signals during different stages of sleep}
\label{fig:EEG_raw}
\end{figure}

In the existing literature, the Gaussian assumption and the stationarity are indispensable.
However, for most data in the real world, its distribution may have the familiar bell-shape,
but it fails one or more statistical normality tests.
Let's take the overnight sleep electroencephalography (EEG) signals in the Sleep-EDF database \cite{kemp2000analysis,goldberger2000physiobank}
as an example. It is well known that overnight EEG signals are non-stationary \cite{kawabata1973nonstationary}.
On the other hand, it is also well known that sleep experts and physicians can ``see into'' the non-stationarity, and summarize some ``patterns'' over each 30-s epochs \cite{rechtschaffen1968manual}. Specifically, there exist some patterns inside the EEG signal so that different sleep stages can be defined, including Awake, rapid eyeball movement (REM) and non-REM, which is further classified into N1, N2 and N3. See the latest rules announced by American Academy of Sleep Medicine (AASM) \cite{iber2007aasm} for details. This ``pattern'' can be viewed as some kind of ``stationarity''. Hence, it reasonable to assume that the EEG signal over each 30-s segment can be well-approximated by a sample path of a stationary process.

\begin{figure}[h]
\centering
\subfigure[][Distribution of EEG data]
{\includegraphics[scale=0.4]{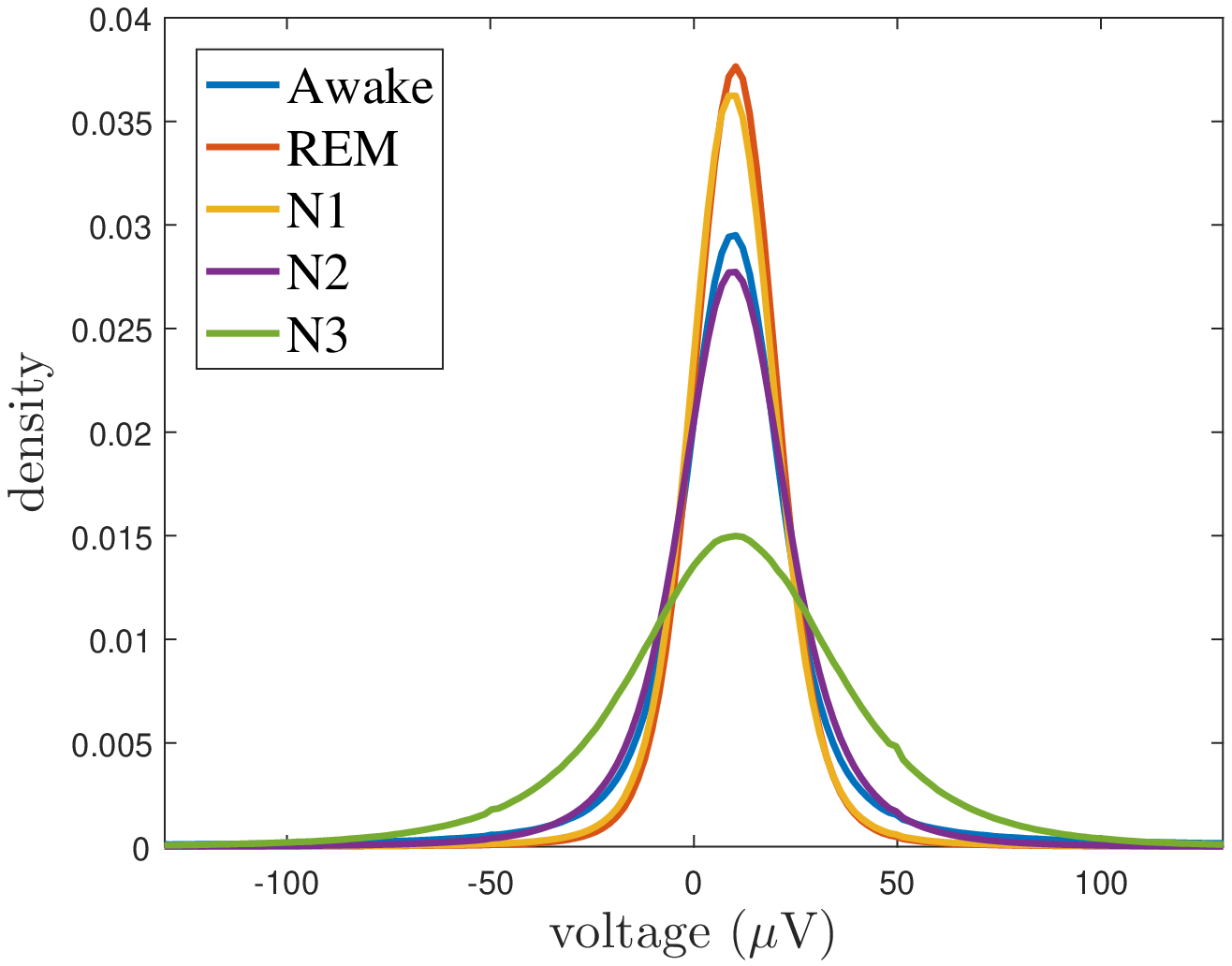}}
\hspace{0.15cm}
\subfigure[][Q-Q Plot]
{\includegraphics[scale=0.4]{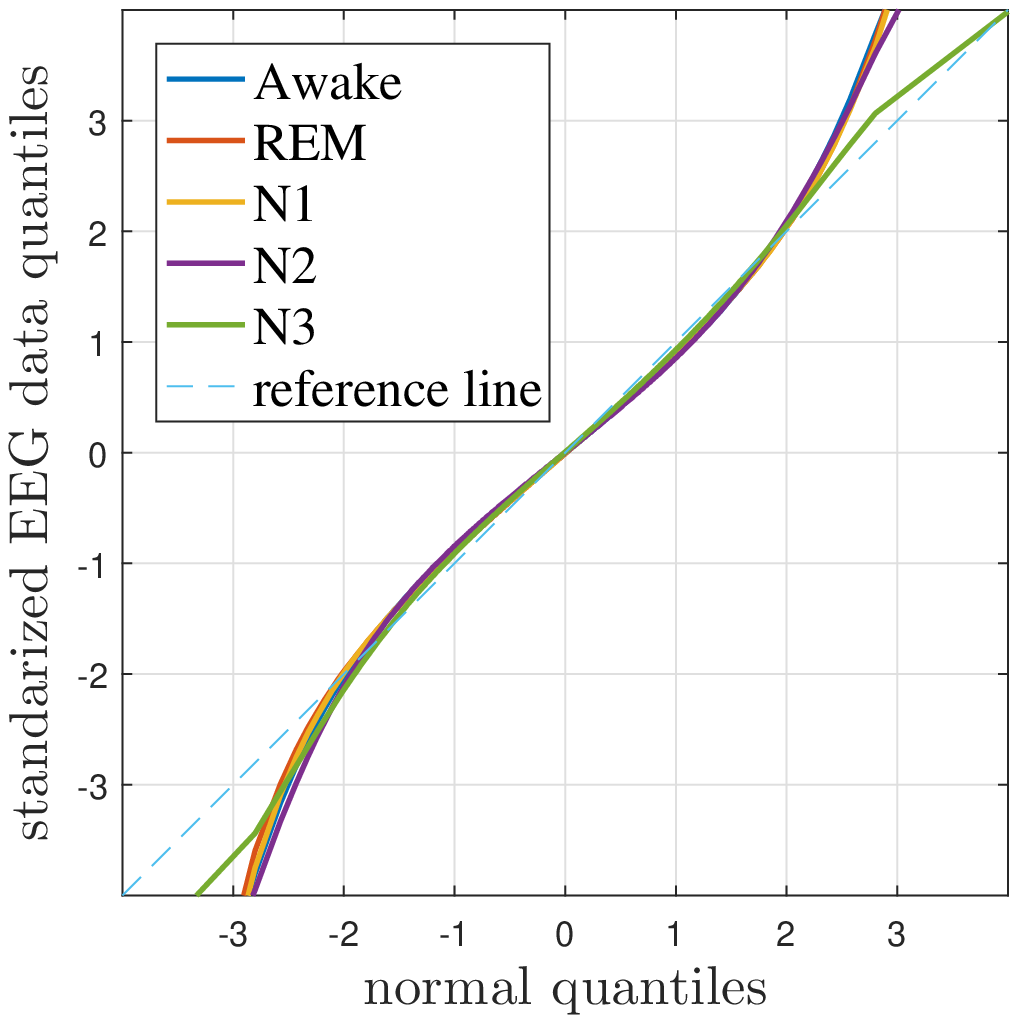}}
\subfigure[][Transformation function $A$]
{\includegraphics[scale=0.4]{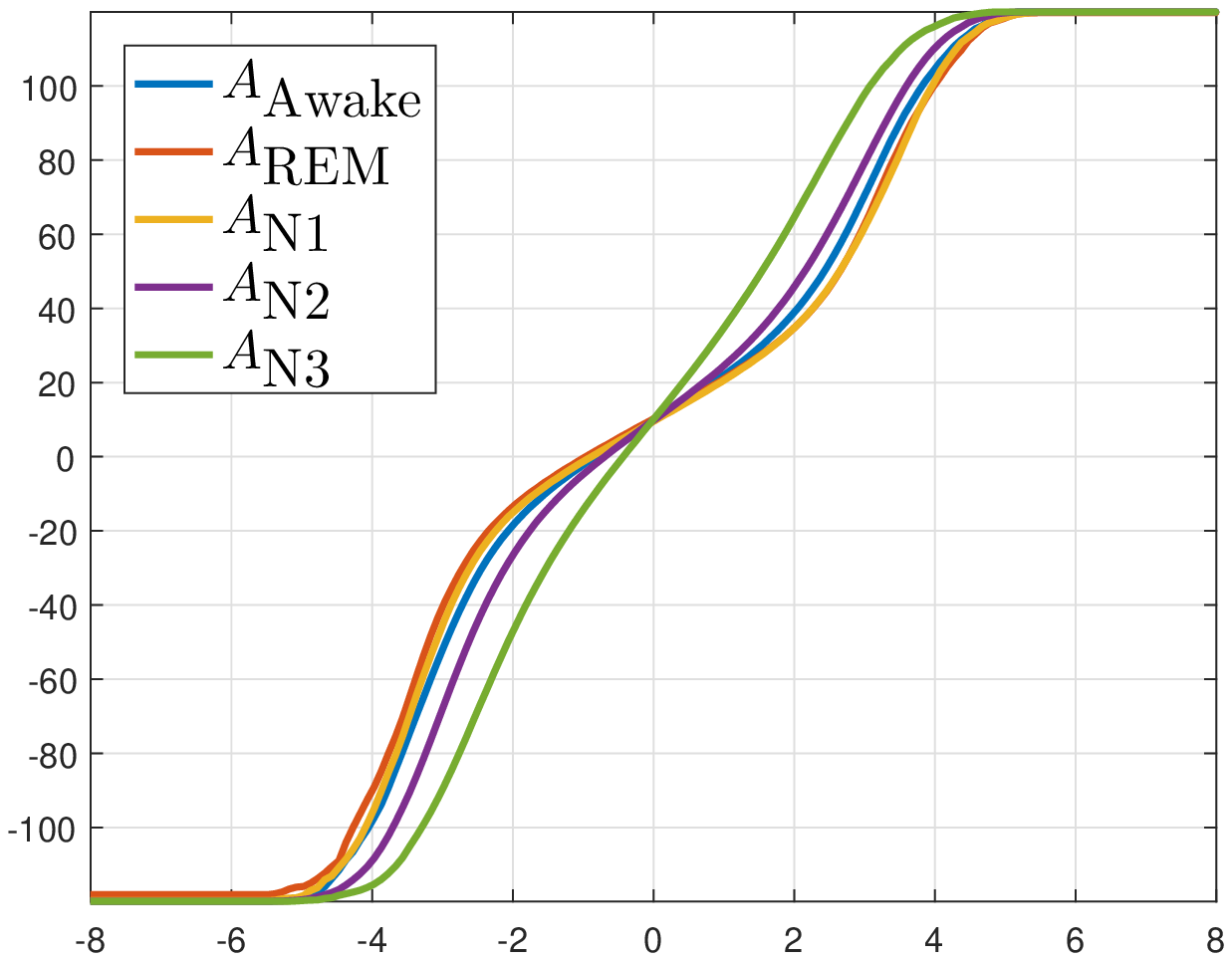}}
\caption{Distributions and Q-Q plots of EEG data during different stages of sleep.}
\label{fig:EEG_stat}
\end{figure}

Figure \ref{fig:EEG_raw} is an illustration of 30-s EEG signals during different stages of sleep.
Let $D_{\textup{Awake}}$, $D_{\textup{REM}}$, $D_{\textup{N1}}$, $D_{\textup{N2}}$ and $D_{\textup{N3}}$ represent the collection of 30-s EEG signal segments recorded during
the stage Awake, REM, N1, N2, and N3, respectively.
Their distributions are plotted in Figure \ref{fig:EEG_stat}(a).
The normal quantile-quantile (Q-Q) plots in Figure \ref{fig:EEG_stat}(b) show that
the tail distributions of these data sets do not look like Gaussian.
This observation motivates us to fit the distribution of the data by introducing a non-random function $A$
to squeeze the Gaussian distribution slightly.
The shape of these functions is shown in Figure \ref{fig:EEG_stat}(c).

%Different from \cite{bruna2015intermittent,liu2020central},
%we consider a non-Gaussian model for the input signal $X$ and let $j_{1}$ and $j_{2}$ go to infinity simultaneously.
%The second scale parameter $j_{2}$ is viewed as a function of $j_{1}$,
%denoted by $j_{2}= j_{2}(j_{1})$,
%and satisfies (\ref{def:coupling}).
%Under some reasonable assumptions and appropriate rescaling, we obtain a non-Gaussian limit from the process $U[j_{1},j_{2}]X$.

\begin{table}\label{List_symbols}
\begin{center}
\begin{tabular}{p{3.9cm}p{13cm}}
\toprule
{\bf Symbol} & {\bf Description}
\\
\midrule
$\star$  & Convolution operator\\
%$j_{k}$ & Scale parameter for the wavelet transform in the $k$-th layer\\
%$\pm\lambda_{1:p}$    & Abbreviation for $\pm(\lambda_{1},\lambda_{2},\ldots,\lambda_{p})\in \mathbb{R}^{p}$\\
%     & Abbreviation for $(-\lambda_{1},-\lambda_{2},\ldots,-\lambda_{p})\in \mathbb{R}^{p}$\\
%$\lambda^{+}_{1:p}$, $d\lambda_{1:p}$,      & Abbreviation for $\lambda_{1}+\lambda_{2}+\cdots+\lambda_{p}$ and $d\lambda_{1}\ d\lambda_{2}\cdots %d\lambda_{p}$\\
$j_{1}$ & Scale parameter for the wavelet transform in the first layer of ST\\
$j_{2}=j_{2}(j_{1})$ & Scale parameter for the wavelet transform in the second layer of ST, which is viewed as a function of $j_{1}$.\\
$\psi$, $\hat{\psi}$, $\psi_{j}(\cdot)$ & Mother wavelet, its Fourier transform, and its scaled version $2^{-j}\psi(2^{-j}\cdot)$\\
$\alpha$  &  Vanishing-moment parameter for $\psi$ defined through $\hat{\psi}(\cdot) = C_{\hat{\psi}}(\cdot)|\cdot|^{\alpha}$\\
$G$       & Stationary Gaussian process\\
$f_{G}$, $f_{G\star\psi_{j_{1}}}$   & Spectral density functions of the random processes $G$ and $G\star\psi_{j_{1}}$\\
$\beta$       & Singularity parameter for $f_{G}$ defined through $f_{G}(\cdot) = C_{G}(\cdot)|\cdot|^{\beta-1}$\\
$A$       & Non-random function\\
$X$       & Input of ST modeled as a subordinated Gaussian process $X=A(G)$\\
$\{\frac{H_{\ell}}{\sqrt{\ell!}}\}_{\ell=0}^{\infty}$ & Orthonormal Hermite polynomials\\
$\{C_{A,\ell}\}_{\ell=0}^{\infty}$ & Hermite coefficients of $A$\\
$\{C_{||,\ell}\}_{\ell=0}^{\infty}$ & Hermite coefficients of the absolute value function $|\cdot|$\\
$S_{j_{1}}$  &  $C_{1}G\star \psi_{j_{1}}$\\
$T_{j_{1}}$  &  $\overset{\infty}{\underset{\ell=2}{\sum}}
C_{\ell}\frac{H_{\ell}\left(G\right)}{\sqrt{\ell!}}\star \psi_{j_{1}}$\\
$\sigma^{2}_{Z}$ & Variance of a stationary random process $Z$, e.g., $\sigma^{2}_{S_{j_{1}}}$, $\sigma^{2}_{G\star \psi_{j_{1}}}$.\\
$Y_{j_{1}}$  &  Normalized wavelet coefficients of $G$, i.e., $G\star\psi_{j_{1}}/\sigma_{G\star\psi_{j_{1}}}.$\\
$U[j_{1}]X$   & Absolute value of the wavelet coefficient of  $X$, i.e., $|X \star \psi_{j_{1}}|$.\\
$U[j_{1},j_{2}]X$ & The second-order ST of $X$ defined by $U[j_{2}]U[j_{1}]X$\\
$a_{n} = O(b_{n})$  & For sequences $\{a_{n}\}$ and $\{b_{n}\}$, there exists
a constant $C>0$ such that $|a_{n}| \leq C |b_{n}|$ when $n$ is sufficiently large.\\
$f(a_{n},x) = O(b_{n}g(x))$  & For functions $f$ and $g$ with the same domain $D$, there exists a
constant $C>0$ such that $|f(a_{n},x)| \leq C |b_{n}g(x)|$ for all $x\in D$ when $n$ is sufficiently large.\\
\bottomrule
\end{tabular}
\end{center}
\caption{List of frequently used symbols and abbreviations}
\end{table}

In this paper we are concerned with the double scaling limit for the feature processes generated from
the second-order ST of strict-sense stationary non-Gaussian processes.
A double scaling limit is a limit in which two or more specific variables are sent to zero or infinity at the same moment.
Double scaling limits are often applied to the random matrix models \cite{bleher2003double,bleher2007large}
and the string theory \cite{giveon1999little}.
The double scaling limit was also considered in
the estimation of high dimensional covariance matrices \cite{ledoit2012nonlinear,couillet2014large}.
Different from \cite{bruna2015intermittent,liu2020central}, the current work considers a class of non-Gaussian processes $X$, each of which is generated by a non-random function $A$ of a stationary Gaussian process
$G$. For commonly used mother wavelets, in Theorem \ref{thm5a}, we proved that
the rescaled process $2^{\frac{j_{1}(\beta-1)}{2}}2^{\frac{j_{2}}{2}}U[j_{1},j_{2}]X(2^{j_{2}}\cdot)$
converges to the absolute value of a Gaussian process in the finite dimensional distribution sense when $j_{1},j_{2}\rightarrow\infty$
with $j_{2}=j_{2}(j_{1})$ and
\begin{equation}\label{def:coupling}
1<\underset{j_{1}\rightarrow\infty}{\lim\inf}\ \frac{j_{2}(j_{1})}{j_{1}}\leq \underset{j_{1}\rightarrow\infty}{\lim\sup}\ \frac{j_{2}(j_{1})}{j_{1}}< \frac{1}{1-\beta},
\end{equation}
where $\beta\in(0,1)$ is the Hurst index used to quantify the long-range dependence of $X$.
To verify Theorem \ref{thm5a}, our approach consists of two parts.
First, in Proposition \ref{lemma:estimateED}, we utilize the Feynman diagram technique to show that the second-order ST of the non-Gaussian process $X$ can be approximated by the second-order ST of a Gaussian process as follows
\begin{equation}
\mathbb{E}\left[\big|U[j_{1},j_{2}]X(2^{j_{2}}\cdot)-U[j_{1},j_{2}]C_{A,1}G(2^{j_{2}}\cdot)\big|^{2}\right] = O\left( 2^{-j_{1}\beta}2^{-j_{2}\beta} +2^{j_{1}(1-2\beta)}2^{-j_{2}}\right)
\end{equation}
when $j_{1},j_{2}\rightarrow \infty$ with $j_{2}=j_{2}(j_{1})$ and $\underset{j_{1}\rightarrow\infty}{\lim\inf}\ \frac{j_{2}(j_{1})}{j_{1}}>1$,
where $C_{A,1}=\int_{\mathbb{R}}A(z)(2\pi)^{-\frac{1}{2}}e^{-\frac{z^{2}}{2}}dz$.
Second, in Proposition \ref{thm:third_order_Gaussian}, we apply the Wiener-It$\hat{\textup{o}}$ expansion to the absolute value function
and use the Feynman diagram technique again to show that
\begin{equation}\label{normalized_UG}
2^{\frac{j_{1}(\beta-1)}{2}}2^{\frac{j_{2}}{2}}\big|G\star\psi_{j_{1}}\big|\star\psi_{j_{2}}(2^{j_{2}}t)\Rightarrow\kappa\int_{\mathbb{R}}
e^{i\lambda t}
\hat{\psi}(\lambda)W(d\lambda)
\end{equation}
in the finite dimensional distribution sense when
$j_{1},j_{2}\rightarrow \infty$ with $j_{2}=j_{2}(j_{1})$ and $\underset{j_{1}\rightarrow\infty}{\lim\inf}\ \frac{j_{2}(j_{1})}{j_{1}}>1$,
where
$\kappa$ is defined in (\ref{def:kappa}) and $W$ is a Gaussian random measure.
From requiring the approximation error to tend to zero when the scale parameters go to infinity and
ensuring the existence of the doubling scaling limit of the approximation process, we obtained the coupling relationship  (\ref{def:coupling}) on the scale parameters under which
the double scaling limit of the second-order scattering transform of non-Gaussian stationary processes
exist.
Moreover, we observed that the spectral density function of the double scaling limit of
$2^{\frac{j_{1}(\beta-1)}{2}}2^{\frac{j_{2}}{2}}U[j_{1},j_{2}]X(2^{j_{2}}\cdot)$
can be expressed in terms of the Hurst index $\beta$ of the non-Gaussian process $X$ and Fourier transform of the mother wavelet $\psi$.

%In these works, the first scale parameter $j_{1}$ remains fixed and the second scale parameter $j_{2}$ goes to infinity.
%Different from \cite{bruna2015intermittent,liu2020central},
%we consider a non-Gaussian model for the input signal $X$ and let $j_{1}$ and $j_{2}$ go to infinity simultaneously.
%The second scale parameter $j_{2}$ is viewed as a function of $j_{1}$,
%denoted by $j_{2}= j_{2}(j_{1})$,
%and satisfies (\ref{def:coupling}).
%Under some reasonable assumptions and appropriate rescaling, we obtain a non-Gaussian limit from the process $U[j_{1},j_{2}]X$.

The rest of the paper is organized as follows.
In Section \ref{sec:preliminary}, we present some preliminaries about the strictly stationary process
and summarize necessary material for ST.
In Section \ref{sec:mainresult}, we state our main
results, including Proposition \ref{lemma:estimateED},
Proposition \ref{thm:third_order_Gaussian} and Theorem \ref{thm5a}.
The proofs of Proposition \ref{lemma:estimateED} and
Proposition \ref{thm:third_order_Gaussian} are given
in Section \ref{sec:proof}. The proofs of technical lemmas are given in the appendix.
Table \ref{List_symbols} contains a list of frequently used symbols and abbreviations.

\vskip 20 pt

\section{Preliminaries}\label{sec:preliminary}

\subsection{Stationary random processes and Wiener-It$\hat{o}$ integrals}

Let $(\Omega, \mathcal{F}, P)$  be an underlying
probability space such that all random elements appeared in this article are defined on it.
Let $\{X(t)\mid t\in \mathbb{R}\}$ be a
strict-sense stationary random process, i.e., for any $d\in \mathbb{N}$, and $\tau,x_{1},...,x_{d},t_{1},...,t_{d}\in \mathbb{R}$,
\begin{equation}\label{def:stationary}
P\left(X(t_{1})\leq x_{1},...,X(t_{d})\leq x_{d}\right) =
P\left(X(t_{1}+\tau)\leq x_{1},...,X(t_{d}+\tau)\leq x_{d}\right).
\end{equation}
Especially, (\ref{def:stationary}) implies that
$P\left(X(t_{1})\leq x\right) =
P\left(X(t_{2})\leq x\right)$ for any $t_{1},t_{2},x\in \mathbb{R}$.
Denote the cumulative distribution function of $X(t)$ by $\Phi_{X}$.
Let $\{G(t)\mid t\in \mathbb{R}\}$ be a stationary Gaussian process with $\mathbb{E}[G(t)]=0$ and $\mathbb{E}[|G(t)|^{2}]=1$.
The cumulative distribution function of $G(t)$ is denoted by $\Phi_{G}$.
%Denote
%\begin{equation*}
%\Phi_{G}(x) = \int_{-\infty}^{x}\frac{1}{\sqrt{2\pi}}e^{-\frac{y^{2}}{2}}dy,\ x\in \mathbb{R}.
%\end{equation*}
If the inverse function $\Phi^{-1}_{X}$ of $\Phi_{X}$ exists,
then $X(t)$ has the same marginal distribution as $\Phi^{-1}_{X}\left(\Phi_{G}(G(t))\right)$
because
\begin{equation}\label{construct_A}
P\left(\Phi^{-1}_{X}\left(\Phi_{G}(G(t))\right)\leq x  \right)
= P\left(\Phi_{G}(G(t))\leq \Phi_{X}(x)  \right)
 = P\left(G(t)\leq \Phi_{G}^{-1}(\Phi_{X}(x))  \right)
 = \Phi_{X}(x)
\end{equation}
for all $x\in \mathbb{R}$ and $t\in \mathbb{R}$.

\begin{Example}\label{example1}
If $X(t)$ has the Gumbel distribution, i.e.,
\begin{equation*}
\Phi_{X}(x) = \textup{Exp}\left(-e^{-\frac{x-c_{1}}{c_{2}}}\right),\ x\in \mathbb{R}
\end{equation*}
where $c_{1},c_{2}\in \mathbb{R}$ are parameters
with $c_{2}>0$, then
\begin{equation*}%\label{phi_inverse_example}
\Phi^{-1}_{X}(z) = c_{1}-c_{2}\ln\left(-\ln(z)\right)\ \textup{for}\ z\in (0,1).
\end{equation*}
\end{Example}
\begin{Example}\label{example2}
If $X(t)$ has the Laplace distribution, i.e.,
\begin{equation*}
\Phi_{X}(x) = \left\{
\begin{array}{ll}
\frac{1}{2}e^{\frac{x-c_{1}}{c_{2}}} & \textup{if}\ x\leq c_{1},\\
1-\frac{1}{2}e^{-\frac{x-c_{1}}{c_{2}}} & \textup{if}\ x> c_{1},
\end{array}\right.
\end{equation*}
where $c_{1},c_{2}\in \mathbb{R}$ are parameters
with $c_{2}>0$, then
\begin{equation*}%\label{phi_inverse_example2}
\Phi^{-1}_{X}(z) = \left\{
\begin{array}{ll}
c_{1}+c_{2}\ln(2z) & \textup{if}\ z\in(0,\frac{1}{2}],\\
c_{1}-c_{2}\ln\left(2(1-z)\right) & \textup{if}\ z\in(\frac{1}{2},1).
\end{array}\right.
\end{equation*}
\end{Example}
Hence, in this work, the input signal $X$ is modeled as
a nonlinear function $A$ of the Gaussian process $G$, that is,
\begin{equation}\label{form:X}
X(t) = A(G(t)),\ t\in \mathbb{R}.
\end{equation}
For the non-random function $A$, we make the following assumption.
\begin{Assumption}\label{assumption:GH}
The function $A$ belongs to the Gaussian Hilbert space
$L^{2}(\mathbb{R},(2\pi)^{-1/2}e^{-\frac{z^{2}}{2}}dz)$,
that is,
\begin{equation}\label{eq:assumption:GH}
\int_{\mathbb{R}}|A(z)|^2 \frac{1}{\sqrt{2\pi}}e^{-\frac{z^{2}}{2}}dz<\infty,
\end{equation}
\end{Assumption}
Note that for the cumulative distributions in Examples \ref{example1} and \ref{example2},
by the well-known inequality
\begin{equation*}%\label{tail_normal}
\frac{x}{x^{2}+1}\frac{1}{\sqrt{2\pi}}e^{-\frac{x^{2}}{2}}\leq 1-\Phi_{G}(x) \leq  \frac{1}{x}\frac{1}{\sqrt{2\pi}}e^{-\frac{x^{2}}{2}},\ \ x>0,
\end{equation*}
$\Phi_{X}^{-1}\Phi_{G}(z)$ has polynomial growth at most when $|z|\rightarrow\infty$.
As shown in Figure \ref{fig:EEG_stat}(c), if the range of data is finite, the function $A$ is bounded,
which automatically satisfies Assumption  \ref{assumption:GH}.

It is known that the Hermite
polynomials
\[H_{\ell}(z)=(-1)^{\ell}e^{\frac{z^{2}}{2}}\frac{d^{\ell}}{dz^{\ell}}e^{-\frac{z^{2}}{2}},\ z\in \mathbb{R},\ \ell=0,1,2,...,
\]
form an orthogonal basis for the Gaussian Hilbert space
with $\int_{\mathbb{R}}(H_{\ell}(z))^{2}\frac{1}{\sqrt{2\pi}}e^{-\frac{z^{2}}{2}}dz = \ell!$.
Under Assumption \ref{assumption:GH},
$A(\cdot)$ has the following expansion
\begin{align}\label{hermiteexpansion}
A(z)=C_{A,0}+\sum_{\ell=1}^{\infty}{C_{A,\ell}}
\frac{H_{\ell}(z)}{\sqrt{\ell!}},
\end{align}
where
\begin{align}\label{hermitecoeff}
C_{A,\ell}=\int_{\mathbb{R}}A(z)\frac{H_{\ell}(z)}{\sqrt{\ell!}}\frac{1}{\sqrt{2\pi}}e^{-\frac{z^{2}}{2}}dz.
\end{align}
%By simple calculation, $H_{0}(z) = 1$, $H_{1}(z) = z$, and $H_{2}(z) = z^{2}-1$.
%Note that the spectral density function $f_{X}$ is nonnegative and integrable with $\int_{\mathbb{R}} f_{X}(\lambda) d\lambda = R_{X}(0)$.
%If the function $C_{X}$ in Assumption 2 satisfies $C_{X}(0)>0$, then $X$ is a long-range dependent process %\cite{doukhan2002theory,pipiras2017long} because
%$f_{X}$ has a singularity at the origin.
%In \cite{gao2004modelling}, the long-range dependent Gaussian processes having spectral densities of the form
%\begin{equation}\label{example_spectral}
%f_{X}(\lambda) = \frac{c_{1}}{|\lambda|^{1-\beta_{1}}(|\lambda|^{2}+c_{2})^{\beta_{2}}},
%\end{equation}
%where $\beta_{1}\in(0,1)$, $\beta_{2}\in(1/2,\infty)$ and $c_{1},c_{2}\in(0,\infty)$, were applied to model the interest rate
%and the Standard and Poor 500 data.
%{\color{red}More discussion about the stationary random processes with singular spectral densities can be referred to %\cite{leonenko1999limit,anh1999possible,doukhan2002theory}.
%The long-range dependence parameter $\beta$ can be obtained by regressing the logged squared wavelet
%coefficients \cite{bardet2000wavelet}.}
%For $\beta_{1} = \beta_{2} = 1$,
%(\ref{example_spectral}) is the spectral density
%of the well-known Ornstein-Uhlenbeck process, which is short-range dependent.
%%
The Hermite rank of $A$ is defined by
\begin{align*}
\textup{rank}\left(A\right) = \textup{inf}\Big\{\ell\geq 1:\ C_{A,\ell}\neq 0\Big\}.
\end{align*}
In view of the strong instability of the Hermite rank \cite{bai2019sensitivity,bai2018instability} caused
by the translation, e.g.,
$\textup{rank}\left(|\cdot|\right) = 2$ but $\textup{rank}\left(|\cdot\pm\varepsilon|\right) = 1$
for any $\varepsilon>0$, we consider the more stable case in the current work.

\begin{Assumption}\label{assumption:rank}
For the function $A$, we suppose that $\textup{rank}\left(A\right) = 1.$
\end{Assumption}

Another key element to describe a random process is its temporal covariance.
Let $R_{G}$ be the covariance function of the underlying Gaussian process $G$.
By the Bochner-Khinchin theorem \cite[Chapter 4]{krylov2002introduction},
there exists a unique nonnegative measure $F_{G}$ on $\mathbb{R}$ such that
$R_{G}$ has the following spectral representation
\begin{equation}\label{eq:BK}
R_{G}(t) = \int_{\mathbb{R}}e^{i\lambda t}F_{G}(d\lambda),\ t\in \mathbb{R}.
\end{equation}
The measure $F_{G}$ is called the spectral measure of $R_{G}$.
\begin{Assumption}\label{assumption:spectral}
The spectral measure $F_{G}$ has the density
$f_{G}$ and
\begin{equation}\label{condition_f}
f_{G}(\lambda) = \frac{C_{G}(\lambda)}{|\lambda|^{1-\beta}},\ \lambda\in \mathbb{R}\setminus\{0\},
\end{equation}
where $\beta\in(0,1)$ is the {\em Hurst index} of long-range dependence and $C_{G}$ is a bounded and continuous function from $\mathbb{R}$ to $[0,\infty)$
such that $f_{G}\in L^{1}$.
\end{Assumption}
Under Assumptions \ref{assumption:rank} and \ref{assumption:spectral},
\cite[claim 2.11]{bai2018instability} shows that the Hurst index of the non-Gaussian process $X$ is just equal to $\beta$.
%In practice, the Hurst index $\beta$ is estimated from the speed of decay of the
%covariance function of $X$.
%The Lebesgue measure of $\mathcal{B} = \{\beta\in(0,1)\mid\ \beta\ell=1\ \textup{for a certain}\ \ell\in \mathbb{N}\}$ is equal to zero.
%Hence, we focus on $\beta$ belonging to the complement of $\mathcal{B}$ in this work.

\begin{figure}[h]
\centering
\subfigure[][correlation functions]
{\includegraphics[scale=0.6]{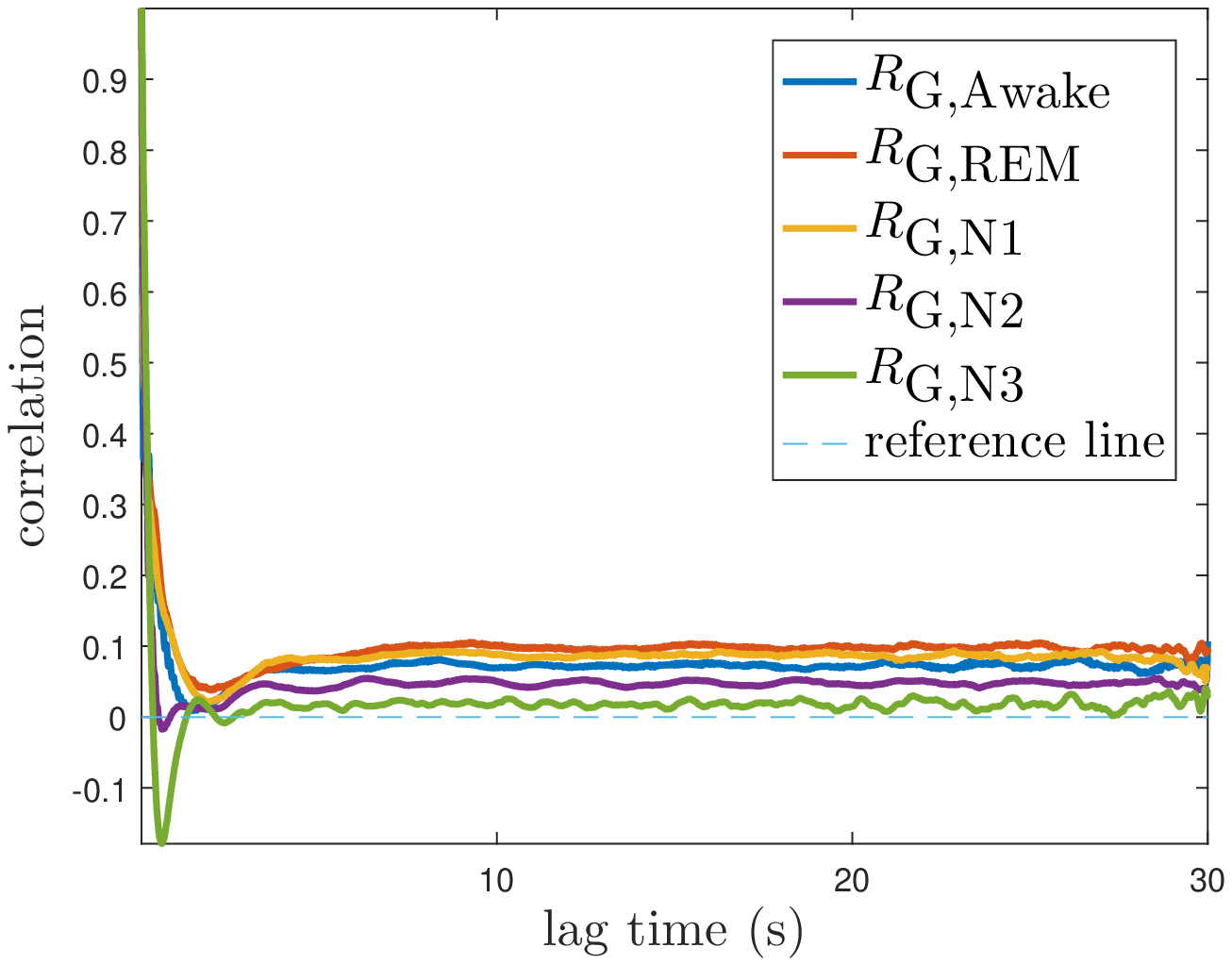}}
\hspace{0.15cm}
\subfigure[][spectral density functions]
{\includegraphics[scale=0.6]{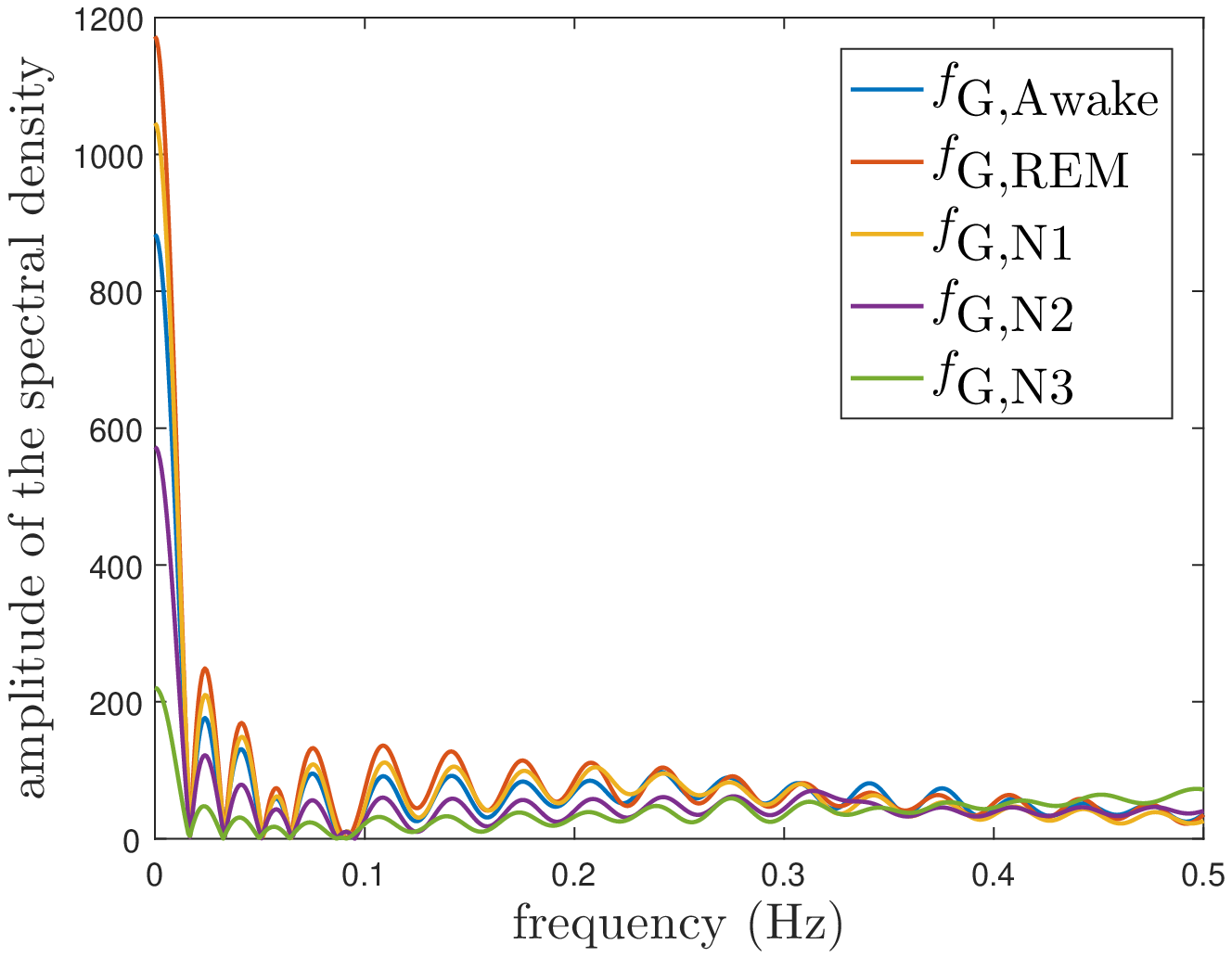}}
\caption{Illustration of the correlation functions and spectral density functions of the underlying Gaussian processes.
}
\label{fig:EEG_cov}
\end{figure}

To provide a more concrete example, we illustrate the shapes of the covariance functions of the underlying Gaussian processes for the EEG example mentioned in Figure \ref{fig:EEG_raw} in Figure \ref{fig:EEG_cov}(a).
First of all, each 30-s signal segment in the data set $D_{\textup{Awake}}$ is transformed by the inverse function of $A_{\textup{Awake}}$ (see Figure \ref{fig:EEG_stat}(d)).
The obtained signal segments are viewed as 30-s sample pathes of a Gaussian process, which is denoted by $G_{\textup{Awake}}$.  We use these sample pathes to estimate the correlation function $R_{G,\textup{Awake}}$ of $G_{\textup{Awake}}$.
The correlation functions, including $R_{G,\textup{REM}}$, $R_{G,\textup{N1}}$, $R_{G,\textup{N2}}$, and $R_{G,\textup{N3}}$, are generated from the data sets $D_{\textup{REM}}$, $D_{\textup{N1}}$, $D_{\textup{N2}}$ and $D_{\textup{N3}}$ in the same way.
The spectral density functions corresponding to these estimated correlation functions are plotted in Figure \ref{fig:EEG_cov}(b),
which shows that the restriction $\beta\in(0,1)$ in (\ref{condition_f}) is a practical assumption.
%It is worth noting that the stationarity does not always hold for the overnight EEG signals,
%but each 30-s EEG signal segment can be viewed as a sample path of a stationary process.

Under Assumption \ref{assumption:spectral}, (\ref{eq:BK}) can be rewritten as
\begin{equation}\label{eq:BK2}
R_{G}(t) = \int_{\mathbb{R}}e^{i\lambda t}f_{G}(\lambda)d\lambda,\ t\in \mathbb{R}.
\end{equation}
By the Karhunen Theorem, the Gaussian process $G$ has the
representation
\begin{align}\label{sample path represent}
G(t)=\int_{\mathbb{R}}
e^{i\lambda t}\sqrt{f_{G}(\lambda)}W(d\lambda),\ t\in \mathbb{R},
\end{align}
where $W(d\lambda)$ is the standard complex-valued Gaussian random measure on $\mathbb{R}$
%; that is, a centered
%orthogonal-scattered Gaussian random  measure on $\mathbb{R}^{n}$
satisfying
\begin{align}\label{ortho}
W(\Delta_{1})=\overline{W(-\Delta_{1})}\ \ \textup{and}\ \ \mathbb{E}\left[W(\Delta_{1})
\overline{W(\Delta_{2})}\right]=\textup{Leb}(\Delta_{1}\cap\Delta_{2})
\end{align}
for all Borel sets
 $\Delta_{1}$  and $\Delta_{2}$
in $\mathbb{R}$, where Leb is the Lebesgue measure on $\mathbb{R}$.

%%%%%%%%%%%%%%%%%%%%%%%%%%%%%%%%%%%%%%%%%%%%%%%55

\subsection{Scattering transform}
The scattering transform was originally proposed by Mallat \cite{mallat2012group} as a tool for extracting information from high-dimensional data. It has several attractive provable properties, including translational invariance, non-expanding variance, and Lipschitz continuous to spatial deformation \cite{mallat2012group}. Its definition is restated here in order to self-contain the paper.

Let $\psi\in L^{1}(\mathbb{R})\cap L^{2}(\mathbb{R})$ be the mother wavelet, which satisfies $\int_{\mathbb{R}} \psi(t) dt = 0$.
A family of real-valued functions $\{\psi_{j}(t)\mid j\in \mathbb{Z}, t \in \mathbb{R}\}$
is called a wavelet family if it is generated from $\psi$ through the dilation procedure
\begin{equation}
\psi_{j}(t) = 2^{-j}\psi(2^{-j}t).
\end{equation}
Denote the Fourier transform of $\psi$ by $\hat{\psi}$, which is defined as
\begin{equation}\label{df:Fourier}
\hat{\psi}(\lambda) = \int_{\mathbb{R}}e^{-i\lambda t}\psi(t) dt.
\end{equation}
Because $\psi\in L^{1}(\mathbb{R})$ and $\int_{\mathbb{R}} \psi(t) dt = 0$,
$\hat{\psi}$ is a continuous function with $\hat{\psi}(0) = 0$. This observation motivates us to make the following assumption.
%for the case $N_{\psi}=1$.

\begin{Assumption}\label{Assumption:1:wavelet}
For the Fourier transform of $\psi$, we assume that there exists a continuous function $C_{\hat{\psi}}$, which is positive at the origin (i.e., $C_{\hat{\psi}}(0)>0$) and has exponential decay at infinity,
such that
\begin{equation}\label{eq:Psi}
\hat{\psi}(\lambda) = C_{\hat{\psi}}(\lambda)|\lambda|^{\alpha}
\end{equation}
for some $\alpha\geq 1$.
\end{Assumption}

%In the following discussion,  $\alpha$ is called the {\it generalized} vanishing moment of $\psi$.

By the Taylor expansion, we know that the parameter $\alpha$ in Assumption \ref{Assumption:1:wavelet} is greater than or equal to the vanishing moment $N$ of $\psi$, where
\begin{equation}\label{df:vanish_moment}
N = \max\left\{n\in \mathbb{N}\,\, \Big|\, \int_{\mathbb{R}}t^{\ell}\psi(t)dt=0\ \textup{for}\ \ell = 0,1,\ldots,n-1,\ \textup{and}\  \int_{\mathbb{R}}t^{n}\psi(t)dt\neq0\right\}.
\end{equation}
Assumption \ref{Assumption:1:wavelet} with $\alpha\geq 1$ holds for commonly used wavelets, including the real part of complex Morlet wavelet ($\alpha=1$), the Mexican hat wavelet ($\alpha=2$), and the $K$-th order Daubechies wavelet ($\alpha=K/2$), where $K = 2,4,\ldots,20.$
%A vanishing moment limits the wavelets ability to represent polynomial behaviour or information in a signal.
%Assumption \ref{Assumption:1:wavelet} with $\alpha\in(0,1)$ is realized by
%the real part of Morse wavelets \cite{olhede2002generalized}, which is a family of analytic wavelets parameterized by two parameters in the %Fourier domain via
%\begin{equation}\label{eq:morse}
%\Psi(\lambda) = \left\{\begin{array}{ll}K_{\alpha,\gamma}\mathcal{H}(\lambda)\lambda^{\alpha}e^{-\lambda^{\gamma}}\ & \textup{for}\ \lambda\geq %0,\\
% 0 \ & \textup{for}\ \lambda< 0,\end{array}
% \right.
%\end{equation}
%where $\alpha>0$, $\gamma>0$, $K_{\alpha,\gamma}$ is a normalizing constant, and $\mathcal{H}$ is the Heaviside unit step function.
For the wavelet function, whose Fourier transform $\hat{\psi}$ vanishes in a neighborhood of the origin, e.g., the Mayer wavelet,
all results in this work still hold.

Given a locally bounded function $X: \mathbb{R}\rightarrow \mathbb{R}$, the continuous wavelet transform  of $X$ is defined as
\begin{equation}\label{df:Wavelet}
     X\star \psi_{j} (t) = 2^{-j}\int_{\mathbb{R}}X(s)\psi(\frac{t-s}{2^{j}})ds\, ,
\end{equation}
where $j\in \mathbb{Z}$ indicates the scale and $t\in \mathbb{R}$ indicates time.
%Compared with the windowed Fourier transform, the continuous wavelet transform is better in capturing short-lived high frequency phenomena, e.g., %the singularities and transient structures, because the time-width of $2^{-j_{1}}\psi(2^{-j_{1}}t)$ is adapted to the scale variable $2^{j_{1}}$ %\cite{daubechies1992ten}.
The first-order ST of $X$ is defined by composing
the wavelet transform of $X$ with the absolute value function
as follows
\begin{equation}\label{def:generalized_1st_ST}
U[j_{1}]X(t) = \Big|X\star \psi_{j_{1}}(t)\Big|,\ j_{1}\in \mathbb{Z}.
\end{equation}
The second-order ST of $X$ is defined by
\begin{equation}\label{df:second_order_scattering}
     U[j_{1},j_{2}]X(t) = \Big|U[j_{1}]X\star \psi_{j_{2}}(t)\Big|
     =
     \Big|2^{-j_{2}}\int_{\mathbb{R}} U[j_{1}]X(s)\psi(\frac{t-s}{2^{j_{2}}})ds\Big|,\ \  j_{1},j_{2}\in \mathbb{Z}.
\end{equation}
In this work, the second scale parameter $j_{2}=j_{2}(j_{1})$ is viewed as a function of $j_{1}$.
%The higher-order scattering transforms
%$$\{U[j_{1},j_{2},...,j_{n}]X(t)\}_{j_{1},j_{2},...,j_{n}\in \mathbb{Z}},\ n\in \mathbb{N},$$ are defined
%in an iterative way.
Note that the continuous wavelet transform (\ref{df:Wavelet}) is better able than the time-frequency representation obtained from the windowed Fourier transform
to capture very short-lived high frequency phenomena, e.g., the singularities and transient structures, because the time-width of $2^{-j_{1}}\psi(2^{-j_{1}}t)$ is adapted to the scale variable $2^{j_{1}}$ \cite{daubechies1992ten,meyer1992wavelets}.
The authors of \cite{mallat2012group}, \cite{anden2014deep},
and \cite{bruna2015intermittent} showed that the second-order ST can be applied to
quantize the frequency of occurrence of various scale phenomena.
In \cite{bruna2015intermittent}, $\mathbb{E}\left[U[j_{1},j_{2}]X(t)\right]$ was estimated with empirical averages
for various stochastic processes, including the fractional Brownian motion, the Poisson process, the $\alpha$-stable Levy processes,
and the multifractal processes.

%In practice, the chosen $j_{2}$ is greater than $j_{1}$ because high frequency characteristics have been filtered
%after convoluting $X$ with $\psi_{j_{1}}$.

\subsection{An approximation of the second-order ST of $X$}

Based on Assumption \ref{assumption:rank},
the second-order ST of $X$ can be rewritten as
\begin{equation}\label{df:second_order_scattering_ST}
     U[j_{1},j_{2}]X =U[j_{1},j_{2}]A(G) =  \Big||S_{j_{1}}+T_{j_{1}}|\star \psi_{j_{2}}\Big|,
\end{equation}
where $S_{j_{1}}$ and $T_{j_{1}}$ are stationary processes defined by
\begin{equation}\label{def:S}
S_{j_{1}}(\tau) = C_{A,1}G\star \psi_{j_{1}}(\tau)
\end{equation}
and
\begin{equation}\label{def:T}
T_{j_{1}}(\tau) =  \overset{\infty}{\underset{\ell=2}{\sum}}
C_{A,\ell}\frac{H_{\ell}\left(G\right)}{\sqrt{\ell!}}\star \psi_{j_{1}}(\tau)
\end{equation}
for $\tau\in \mathbb{R}$.
Because the computation of the modulus operator $|\cdot|$ is followed by a convolution with the wavelet $\psi_{j_{2}}$,
we need a calculable expression for the absolute value of
the sum of the processes $S_{j_{1}}$ and $T_{j_{1}}$.
To reach it, we need the following property \cite{major1981lecture,ito1951multiple,houdre1994chaos}:
\begin{align}\label{expectionhermite}
\mathbb{E}[H_{\ell_{1}}(G(t_{1}))H_{\ell_{2}}(G(t_{2}))]=
\delta^{\ell_{1}}_{\ell_{2}} \ell_{1}!
R_{G}^{\ell_{1}}(t_{1}-t_{2}),
\ \ \ t_{1},\ t_{2}\in \mathbb{R},
\end{align}
where $\delta^{\ell_{1}}_{\ell_{2}}$ is the Kronecker delta.
%is used to get the asymptotic behavior of $\sigma_{S_{j_{1}}}^{2} = \mathbb{E}\left[|S_{j_{1}}(\tau)|^{2}\right]$
%and $\sigma_{T_{j_{1}}}^{2} = \mathbb{E}\left[|T_{j_{1}}(\tau)|^{2}\right]$ when $j_{1}$ is sufficiently large.
%which do not depend on $\tau$ due to the stationarity of
%$S_{j_{1}}$ and $T_{j_{1}}$.
%The proof of (\ref{expectionhermite})
%can be referred to
%\cite[Corollary 5.5 and p. 30]{major1981lecture}, \cite{ito1951multiple} or \cite{houdre1994chaos}.
%\section{Properties of multi-fold convolutions of spectral densities}\label{sec:appendix_conv}
The lemma below is about the behavior of the spectral density of the $\ell$-th power of the covariance function $R_{G}$ near the origin.
\begin{Lemma}{\cite[Lemma 1]{liu2015multi}}\label{lemma:convolution}
%\cite[Appendix H]{liu2020central}
For any $\ell\geq 2$,
there exists a bounded and continuous function $B_{\ell}$ such that
\begin{equation}\label{singularconvolution1}
f_{G}^{\star \ell}(\lambda)=
\left\{\begin{array}{lr}
B_{\ell}(\lambda)|\lambda|^{\ell\beta-1}\ \ & \textup{for}\ \ell\beta<1,
\\
B_{\ell}(\lambda)\textup{ln}(2+\frac{1}{|\lambda|})\ \ & \textup{for}\ \ell\beta=1,
\\
B_{\ell}(\lambda)\ \ & \textup{for}\ \ell\beta> 1,
\end{array}
\right.
\end{equation}
almost everywhere.
Moreover, for any $\ell_{2} >\ell_{1} > 1/\beta$,
$\underset{\lambda\in \mathbb{R}}{\sup} B_{\ell_{2}}(\lambda)
\leq \underset{\lambda\in \mathbb{R}}{\sup} B_{\ell_{1}}(\lambda)$.
\end{Lemma}
The lemma below is about the asymptotic behavior of $\sigma_{S_{j_{1}}}^{2} = \mathbb{E}\left[|S_{j_{1}}(\tau)|^{2}\right]$
and $\sigma_{T_{j_{1}}}^{2} = \mathbb{E}\left[|T_{j_{1}}(\tau)|^{2}\right]$ when $j_{1}$ is sufficiently large.
\begin{Lemma}\label{lemma:var_S_T}
Under Assumptions \ref{assumption:spectral} and \ref{Assumption:1:wavelet},
the covariance function of $S_{j_{1}}$ has the spectral representation
\begin{equation}\label{spectral_S}
R_{S_{j_{1}}}(t) = C_{A,1}^{2}\int_{\mathbb{R}}e^{it\lambda}f_{G}(\lambda)|\hat{\psi}(2^{j_{1}}\lambda)|^{2}\ d\lambda.
\end{equation}
When $j_{1}\rightarrow\infty$,
\begin{align}\label{lemma:varS}
\underset{j_{1}\rightarrow\infty}{\lim}2^{j_{1}\beta}\sigma_{S_{j_{1}}}^{2}=
\underset{j_{1}\rightarrow\infty}{\lim}2^{j_{1}\beta}R_{S_{j_{1}}}(0) =C_{A,1}^{2}\sigma^{2},
\end{align}
where $\sigma^{2}=C_{G}(0)\int_{\mathbb{R}}|\hat{\psi}(\lambda)|^{2}|\lambda|^{\beta-1} d\lambda$.
On the other hand, the variance of $T_{j_{1}}$ has the asymptotic behavior:
\begin{align}\label{lemma:varT}
\sigma_{T_{j_{1}}}^{2} = \left\{\begin{array}{ll}
O\left(2^{-2j_{1}\beta}\right) & \ \textup{if}\ \beta<1/2,
\\
O\left(j_{1}2^{-j_{1}}\right) & \ \textup{if}\ \beta=1/2,
\\
O\left(2^{-j_{1}}\right) & \ \textup{if}\ \beta>1/2,
\end{array}
\right.
\end{align}
where $O(\cdot)$ is the asymptotic notation defined in Table \ref{List_symbols}.
\end{Lemma}
The proof of Lemma \ref{lemma:var_S_T} can be found in \ref{sec:proof:var_S_T}.
Based on Lemma \ref{lemma:var_S_T}, we attempt to show that $|S_{j_{1}}| \geq |T_{j_{1}}|$ for almost all sample pathes of the Gaussian process
$G$ whenever $j_{1}$ is large enough.
\begin{Lemma}\label{lemma:path_dominate}
Suppose that the function $A$ satisfies Assumptions \ref{assumption:GH}-\ref{assumption:rank}, the Gaussian process $G$ satisfies Assumption \ref{assumption:spectral}, and the mother wavelet $\psi$ satisfies Assumption \ref{Assumption:1:wavelet}.
For each $j_{1}\in \mathbb{N}$, let $\tau_{j_{1},1},...,\tau_{j_{1},n_{j_{1}}}$ be points randomly sampled
from $\mathbb{R}$, where $n_{j_{1}}=[2^{\frac{\beta}{6}j_{1}}]$ for $\beta\in(0,\frac{1}{2}]$ and
$n_{j_{1}}=[2^{\frac{1-\beta}{6}j_{1}}]$ for $\beta\in(\frac{1}{2},1)$.
Define the event
\begin{align}\label{def:event}
\mathcal{E}_{j_{1}}:=\left\{|S_{j_{1}}(\tau_{j_{1},k})|<|T_{j_{1}}(\tau_{j_{1},k})|\ \textup{for some}\ k\in\{1,2,...,n_{j_{1}}\}\right\}
\end{align}
for each $j_{1}\in \mathbb{N}$, where $S_{j_{1}}$ and $T_{j_{1}}$ are defined in (\ref{def:S}) and (\ref{def:T}), respectively.
We have
\begin{align}\label{statement_io}
P\left(\mathcal{E}_{j_{1}}\ \textup{occurs infinitely often as}\ j_{1}\rightarrow\infty\right)=0.
\end{align}

\end{Lemma}
The proof of Lemma \ref{lemma:path_dominate} can be found in \ref{sec:proof:lemma:io_path}.
%The time points $\tau_{j_{1},1},...,\tau_{j_{1},n}$ mentioned in Lemma \ref{lemma:path_dominate} can also be uniformly sampled from
%$[-2^{j_{1}\frac{\beta}{24}},2^{j_{1}\frac{\beta}{24}}]$, i.e.,
%$\tau_{j_{1},k}=-2^{j_{1}\frac{\beta}{24}}+2^{-j_{1}\frac{\beta}{12}}k$ for $k=1,...,n$.
Note that
\begin{equation}\label{eq:abs-abs}
|x+y|-|x| = \textup{sign}(x)y\ \ \textup{if}\  \ |x|\geq |y|.
\end{equation}
By Lemma \ref{lemma:path_dominate} and (\ref{eq:abs-abs}), for almost every sample path of the Gaussian process $G$, when $j_{1}$ is sufficiently large,
\begin{align}\notag
&\big|X\star \psi_{j_{1}}(\tau_{j_{1},k})\big|-\big|C_{A,1}G\star \psi_{j_{1}}(\tau_{j_{1},k})\big|
\\\notag=&\big|S_{j_{1}}(\tau_{j_{1},k})+T_{j_{1}}(\tau_{j_{1},k})\big|-\big|S_{j_{1}}(\tau_{j_{1},k})\big|
\\\label{|S+T|-|S|}=& \textup{sign}\left(S_{j_{1}}(\tau_{j_{1},k})\right)T_{j_{1}}(\tau_{j_{1},k}),\ k=1,2,...,n_{j_{1}},
\end{align}
where
$\{\tau_{j_{1},k}\mid k=1,...,n_{j_{1}}\}$
are  selected from $\mathbb{R}$ randomly for each $j_{1}$.
Denote the difference between the wavelet coefficients of $\big|X\star \psi_{j_{1}}\big|$ and
$\big|C_{A,1}G\star\psi_{j_{1}}\big|$ by
\begin{align}\notag
\mathcal{D}_{j_{1},j_{2}}(t)=&\int_{\mathbb{R}}\left[\big|X\star \psi_{j_{1}}\left(\tau\right)\big| -\big|C_{A,1}G\star\psi_{j_{1}}\left(\tau\right)\big|\right]\psi_{j_{2}}(t-\tau)d\tau
\\\label{D_true}=&\int_{\mathbb{R}} \left[\big|S_{j_{1}}(\tau)+T_{j_{1}}(\tau)\big|-\big|S_{j_{1}}(\tau)\big|\right]\psi_{j_{2}}(t-\tau)d\tau
\end{align}
and
\begin{align}\notag
\widetilde{\mathcal{D}}_{j_{1},j_{2}}(t)=&\int_{\mathbb{R}}
\textup{sign}\left(S_{j_{1}}\left(\tau\right)\right)T_{j_{1}}(\tau)\psi_{j_{2}}(t-\tau)d\tau.
\end{align}
Note that $\mathbb{E}[|\mathcal{D}_{j_{1},j_{2}}(t)|^{2}]$  and
$\mathbb{E}[|\widetilde{\mathcal{D}}_{j_{1},j_{2}}(t)|^{2}]$ do not depend on $t$ due to the strict stationarity of $X$.
From the observation (\ref{|S+T|-|S|}), we make the following assumption.
\begin{Assumption}\label{Assumption:passlimit}
When $j_{1}$ is sufficiently large, there exists a constant $M>0$ such that
\begin{align}
\mathbb{E}[|\mathcal{D}_{j_{1},j_{2}}(t)|^{2}] \leq M
\mathbb{E}[|\widetilde{\mathcal{D}}_{j_{1},j_{2}}(t)|^{2}]
%\mathbb{E}\left[\Big|\int_{\mathbb{R}}\textup{sign}\left(S_{j_{1}}(\tau)\right)T_{j_{1}}(\tau)\psi_{j_{2}}(t-\tau)d\tau\Big|^{2}\right].
\end{align}
for all $j_{2}\in \mathbb{Z}$ and $t\in \mathbb{R}$.
\end{Assumption}
\begin{figure}[h]
\centering
\subfigure[][$j_{2}(j_{1})=j_{1}$]
{\includegraphics[scale=0.6]{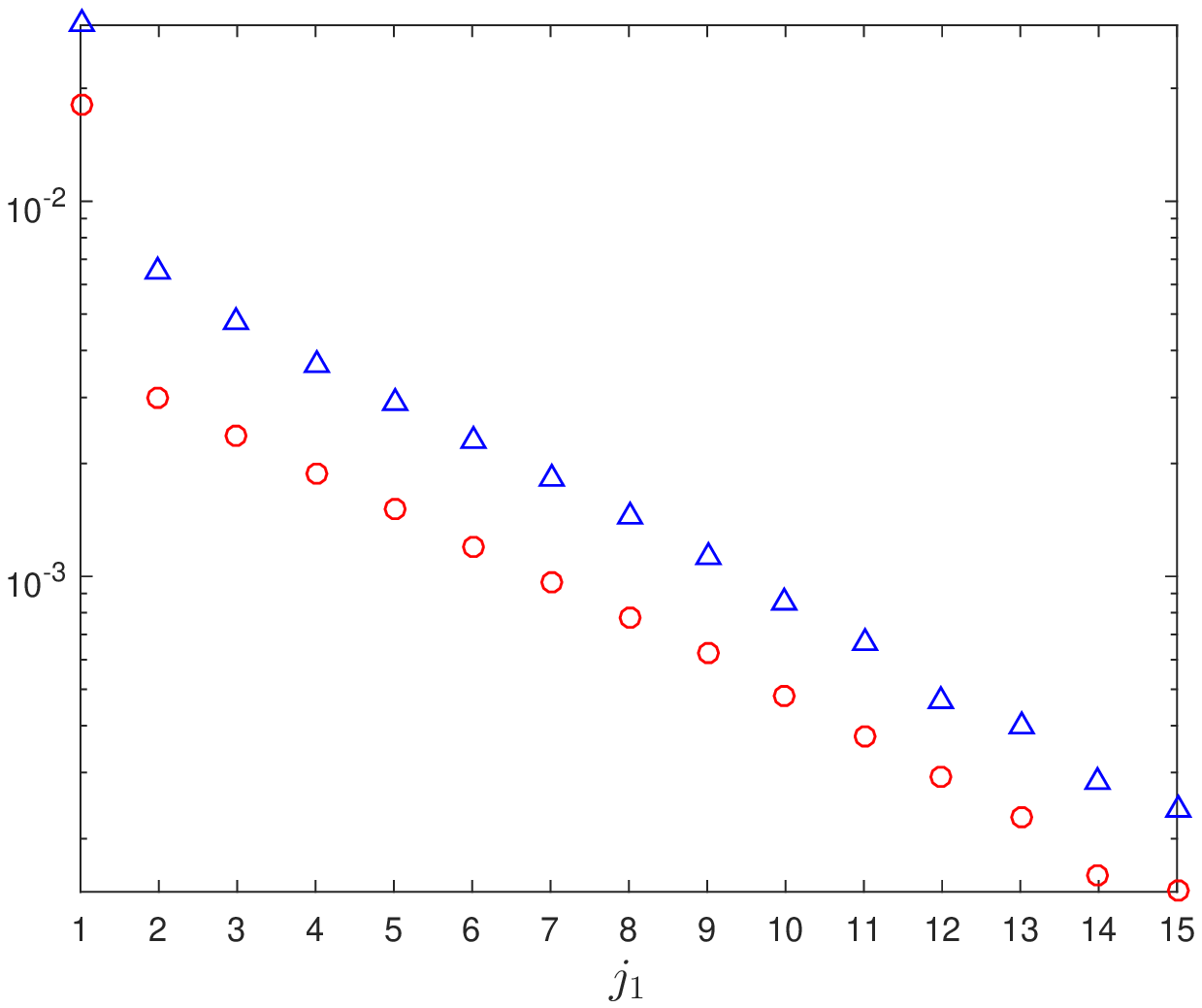}}
\hspace{0.1cm}
\subfigure[][$j_{2}(j_{1})=1.1j_{1}$]
{\includegraphics[scale=0.6]{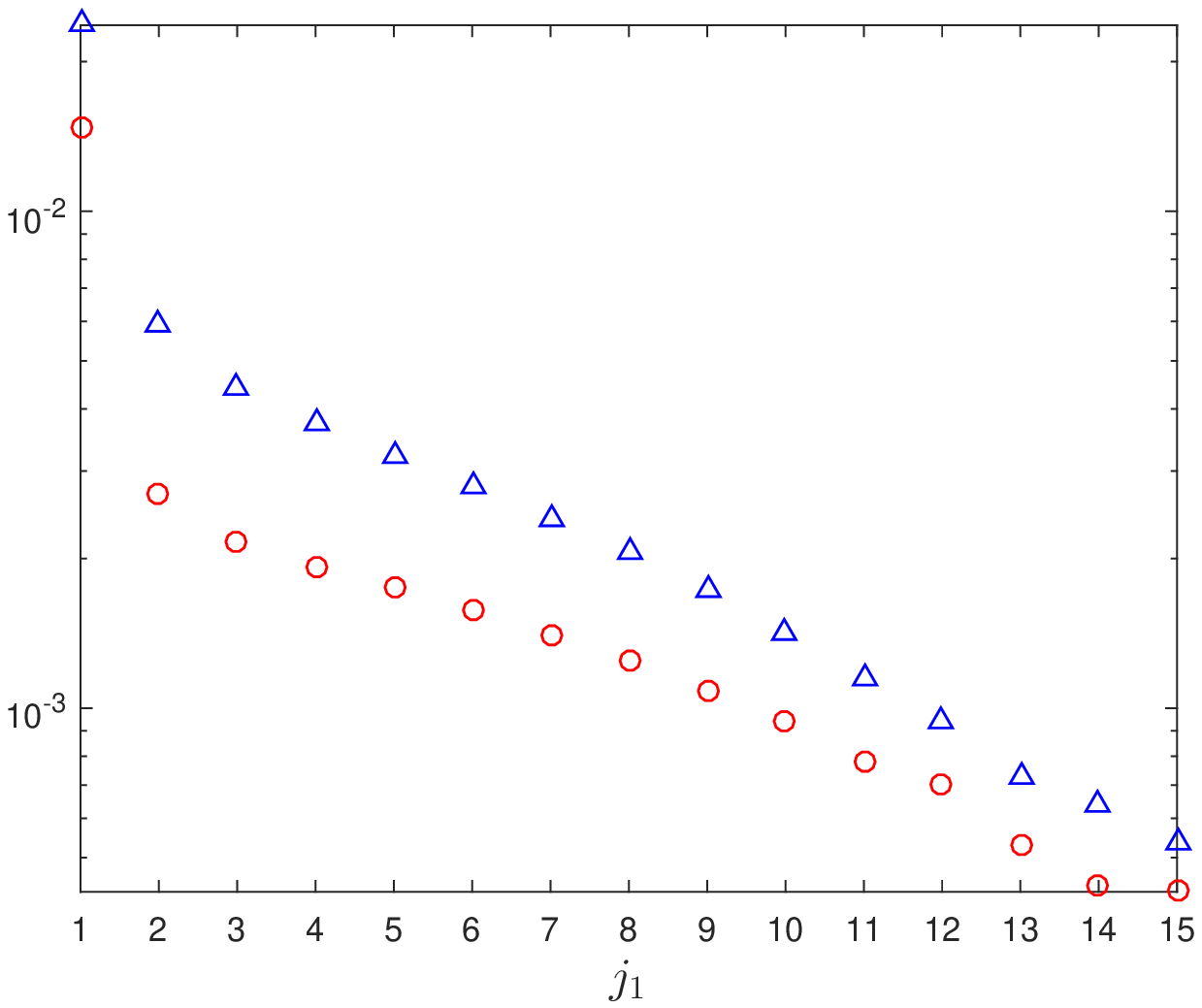}}
\caption{Numerical validation of Assumption \ref{Assumption:passlimit}.  The values of
$\mathbb{E}[|\mathcal{D}_{j_{1},j_{2}}(t)|^{2}]$ and $\mathbb{E}[|\widetilde{\mathcal{D}}_{j_{1},j_{2}}(t)|^{2}]$
are marked by {\color{red}$\bigcirc$} and {\color{blue}$\bigtriangleup$}, respectively.
Here,
$X=G+(G^{2}-1)+(G^{3}-3G)$, where $G$ is a long-range dependent Gaussian process with the spectral density function
$f_{G}(\lambda)= |\lambda|^{1-0.1}1_{(0,36)}(\lambda)$, and $\psi$ is the Daubechies wavelet (db8).}
\label{fig:assumption5}
\end{figure}

\begin{figure}[h]
\centering
\subfigure[][$j_{2}(j_{1})=j_{1}$]
{\includegraphics[scale=0.6]{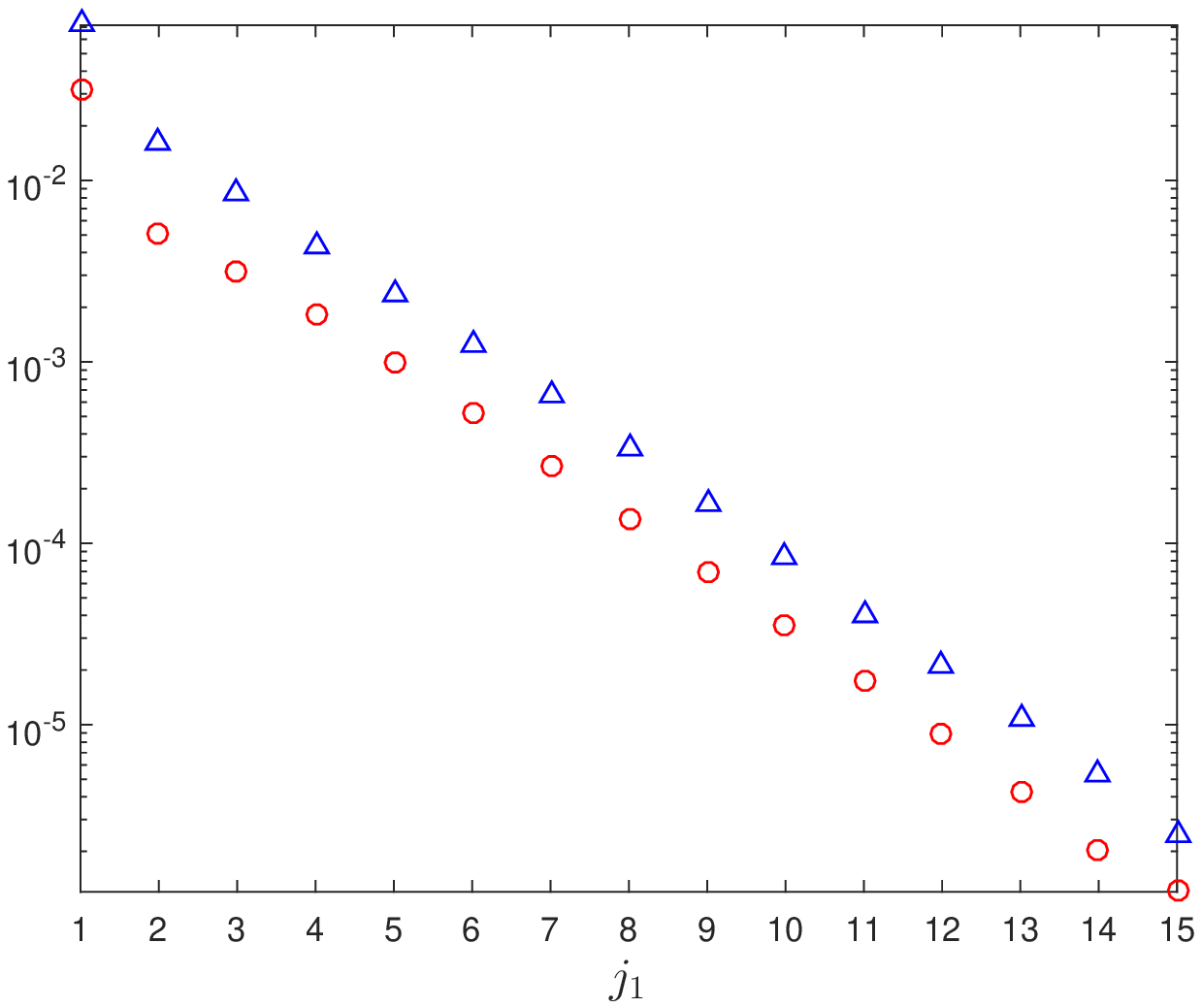}}
\hspace{0.1cm}
\subfigure[][$j_{2}(j_{1})=1.1j_{1}$]
{\includegraphics[scale=0.6]{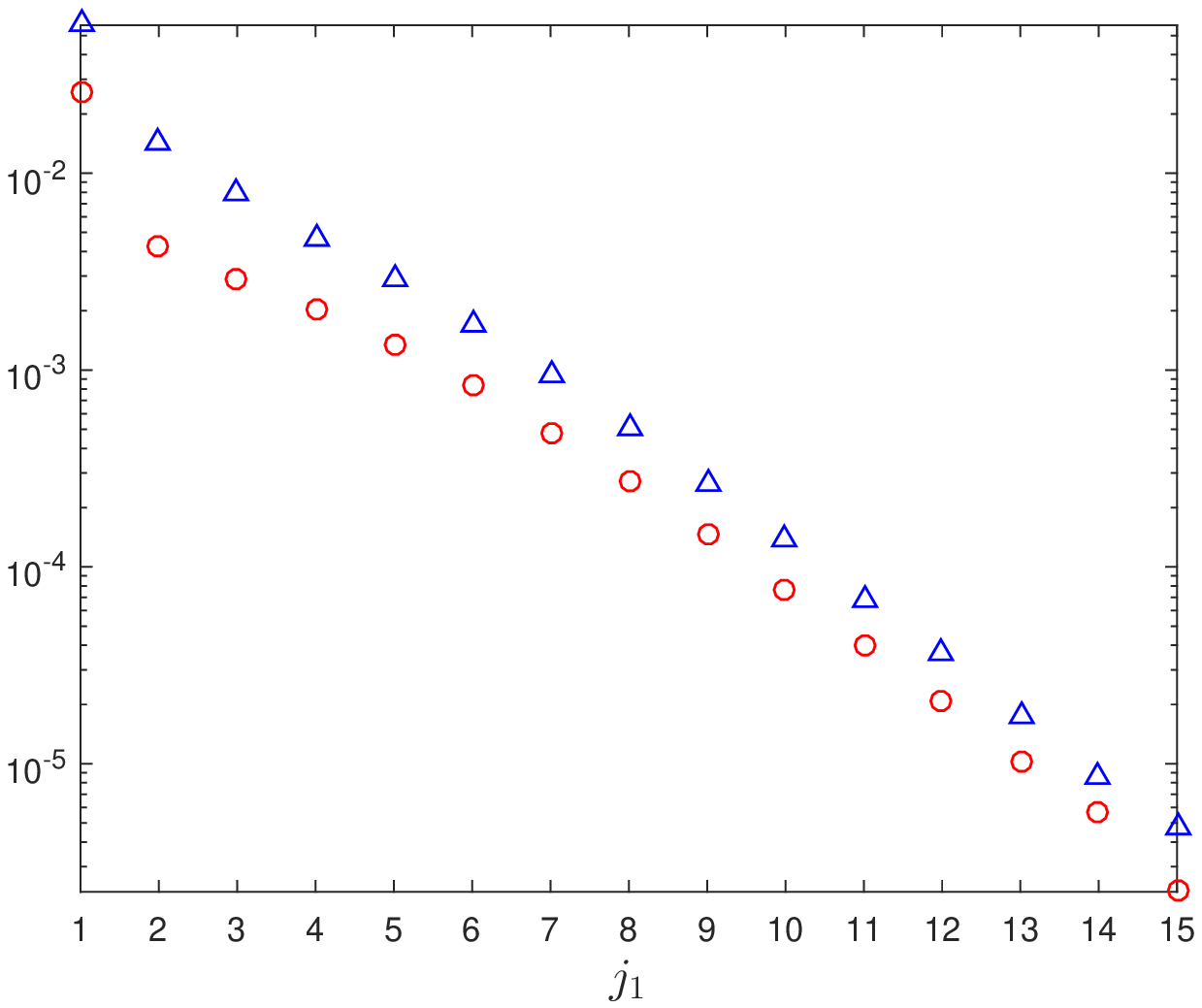}}
\caption{Numerical validation of Assumption \ref{Assumption:passlimit}.  The values of
$\mathbb{E}[|\mathcal{D}_{j_{1},j_{2}}(t)|^{2}]$ and $\mathbb{E}[|\widetilde{\mathcal{D}}_{j_{1},j_{2}}(t)|^{2}]$
are marked by {\color{red}$\bigcirc$} and {\color{blue}$\bigtriangleup$}, respectively.
Here,
$X=G+(G^{2}-1)+(G^{3}-3G)$, where $G$ is a short-range dependent Gaussian process with the spectral density function
$f_{G}(\lambda)= (1+|\lambda|^{2})^{-1}1_{(0,36)}(\lambda)$, and $\psi$ is the Daubechies wavelet (db8).}
\label{fig:assumption5_weak}
\end{figure}

\begin{Remark}
In Figure \ref{fig:assumption5}, we implement a long-range dependent Gaussian process $G$ with the Hurst index $\beta=0.1$.
In our main result, i.e., Theorem \ref{thm5a} below, we require that $\underset{j_{1}\rightarrow\infty}{\lim\sup}\ \frac{j_{2}(j_{1})}{j_{1}}<\frac{1}{1-\beta}$, so we consider the setup $j_{2}(j_{1}) = 1.1j_{1}$ in Figure \ref{fig:assumption5}(b).
On the other hand, from (\ref{|S+T|-|S|}) and the semi-log graphs in Figure \ref{fig:assumption5},
we expect that $M=1$, but how to prove it is still an open question.
In Figure \ref{fig:assumption5_weak}, we further implement a short-range dependent Gaussian process $G$
with the spectral density $f_{G}(\lambda)= (1+|\lambda|^{2})^{-1}1_{(0,36)}(\lambda)$. The numerical result also supports
Assumption \ref{Assumption:passlimit} though we only focus on the long-range dependent case $\beta\in(0,1)$.
\end{Remark}

\section{Main Results}\label{sec:mainresult}
For finding the double scaling limit of the second-order ST of the random process $X$,
we first estimate $\mathbb{E}[|\mathcal{D}_{j_{1},j_{2}}(t)|^{2}]$ through
$\mathbb{E}[|\widetilde{\mathcal{D}}_{j_{1},j_{2}}(t)|^{2}]$
and then explore the double scaling limit of $|G\star \psi_{j_{1}}|\star \psi_{j_{2}}$.
Finally, we apply Slutsky's argument (see, for example, the book of Leonenko \cite[p. 6]{leonenko1999limit}) and the continuous mapping theorem
to get the double scaling limit of $U[j_{1},j_{2}]X$.

\begin{Proposition}\label{lemma:estimateED}
Let $X=A(G)$, where $A$ is a non-random function satisfying Assumptions \ref{assumption:GH}-\ref{assumption:rank} and $G$ is a stationary Gaussian process satisfying Assumption \ref{assumption:spectral}, and let $\psi$ be a real-valued mother wavelet function satisfying Assumption \ref{Assumption:1:wavelet}.
%For any $t\in \mathbb{R}$, define
%\begin{align}\notag
%\mathcal{D}_{j_{1},j_{2}}(t)=\int_{\mathbb{R}}\big|X\star \psi_{j_{1}}(\tau)\big| \psi_{j_{2}}(t-\tau)d\tau-\int_{\mathbb{R}}\big|C_{1}G\star %\psi_{j_{1}}(\tau)\big|\psi_{j_{2}}(t-\tau)d\tau.
%\end{align}
%If Assumption \ref{Assumption:passlimit} also holds,
We have
\begin{equation}\label{estimate_D_combine}
\mathbb{E}[\widetilde{\mathcal{D}}_{j_{1},j_{2}}(0)^{2}] =
\left\{\begin{array}{ll} O\left( 2^{-j_{1}\beta}2^{-j_{2}\beta} +2^{j_{1}(1-2\beta)}2^{-j_{2}}\right)\ & \ \textup{for} \ \beta\in(0,\frac{1}{2}),\\
O\left(j_{1}2^{-j_{2}}+ 2^{-j_{1}/2}2^{-j_{2}/2} +2^{-j_{2}}\right)\ & \ \textup{for} \ \beta=\frac{1}{2},
\\
O\left(2^{-j_{2}}+ 2^{-j_{1}\beta}2^{-j_{2}\beta} +2^{j_{1}(1-2\beta)}2^{-j_{2}}\right)\ & \ \textup{for} \ \beta\in(\frac{1}{2},1),
\end{array}\right.
\end{equation}
when $j_{1},j_{2}\rightarrow \infty$ with $\underset{j_{1}\rightarrow\infty}{\lim\inf}\ \frac{j_{2}(j_{1})}{j_{1}}>1$.
\end{Proposition}

The proof of Proposition \ref{lemma:estimateED} is in Section \ref{sec:proof:lemma:estimateED}.
Next, we explore the double scaling limit of $|G\star \psi_{j_{1}}|\star \psi_{j_{2}}$.
The Wiener-It$\hat{\textup{o}}$ decomposition is applied to deal with the random process $|G\star \psi_{j_{1}}|$
as follows
\begin{align}
|G\star \psi_{j_{1}}| = \sigma_{G\star \psi_{j_{1}}}
\overset{\infty}{\underset{\ell=0}{\sum}}
\frac{C_{||,\ell}}{\sqrt{\ell!}}H_{\ell}\left(\frac{G\star \psi_{j_{1}}}{\sigma_{G\star \psi_{j_{1}}}}\right),
\end{align}
where $\sigma_{G\star \psi_{j_{1}}}$ is the standard deviation of $G\star \psi_{j_{1}}$
and
$\{C_{||,\ell}\}_{\ell=0,1,...}$ are the Hermite coefficients of the absolute value function.
Note that $C_{||,\ell}=0$ if $\ell$ is odd because $|\cdot|$ is an even function.
By the property (\ref{expectionhermite}),
\begin{align}\label{expect_H_Gpsi}
\mathbb{E}\left[
H_{\ell_{1}}\left(\frac{G\star \psi_{j_{1}}(t_{1})}{\sigma_{G\star \psi_{j_{1}}}}\right)
H_{\ell_{2}}\left(\frac{G\star \psi_{j_{1}}(t_{2})}{\sigma_{G\star \psi_{j_{1}}}}\right)\right]
= \ell!\sigma_{G\star \psi_{j_{1}}}^{-2\ell}R_{G\star \psi_{1}}^{\ell}(t_{1}-t_{2}),\ t_{1},t_{2}\in \mathbb{R},
\end{align}
if $\ell_{1}=\ell_{2}=\ell\geq0$, where $R_{G\star \psi_{1}}$ is the covariance function of the process $G\star\psi_{j_{1}}$.
Otherwise, i.e., $\ell_{1}\neq\ell_{2}$, the expectation in (\ref{expect_H_Gpsi}) is equal to zero.
Denote the spectral density of $R_{G\star \psi_{1}}$ by $f_{G\star\psi_{j_{1}}}$.
The following property about the spectral density of $R^{\ell}_{G\star \psi_{1}}$, i.e.,
the $\ell$-fold convolution of the spectral density of $G\star\psi_{j_{1}}$, is needed
in the proof of Proposition \ref{thm:third_order_Gaussian} below.

\begin{Lemma}\label{Lemma:fstarj1j2_limit}
Under Assumptions \ref{assumption:spectral} and \ref{Assumption:1:wavelet},
for any $\ell\in\{2,4,6,...\}$,
$
\|f^{\star \ell}_{G\star\psi_{j_{1}}}\|_{\infty}
= f^{\star \ell}_{G\star\psi_{j_{1}}}(0).
$
Furthermore, for any $\lambda\in \mathbb{R}$, $\ell\geq2$, and $j_{2}=j_{2}(j_{1})$ satisfying $\underset{j_{1}\rightarrow\infty}{\lim\inf}\ \frac{j_{2}(j_{1})}{j_{1}}>1$,
\begin{align}\notag
\gamma_{\ell}:=&\underset{\begin{subarray}{c}
j_{1}\rightarrow\infty\\
j_{2}=j_{2}(j_{1})\end{subarray}}{\lim}2^{-j_{1}}\sigma_{G\star j_{1}}^{-2\ell} f^{\star \ell}_{G\star\psi_{j_{1}}}(2^{-j_{2}}\lambda)
\\\notag=&\underset{j_{1}\rightarrow\infty}{\lim}2^{-j_{1}}\sigma_{G\star \psi_{j_{1}}}^{-2\ell} f^{\star \ell}_{G\star\psi_{j_{1}}}(0)
\\\label{lemma_f_cov_limit_equiv}
=&\frac{1}{2\pi} \int_{\mathbb{R}}\left[\int_{\mathbb{R}} \frac{|\hat{\psi}(\eta)|^2}{|\eta|^{1-\beta}}d\eta\right]^{-\ell}
\left[ \int_{\mathbb{R}}e^{it\zeta} \frac{|\hat{\psi}(\zeta)|^{2}}{|\zeta|^{1-\beta}}
d\zeta\right]^{\ell} dt
\end{align}
and $\gamma_{2}\geq \gamma_{4}\geq \gamma_{6} \geq \cdots\geq 0$.

%\begin{equation}\label{lemma_f_cov_limit_equiv}
%\underset{j_{2},j_{3}\rightarrow\infty: j_{3}=\delta j_{2}}{\lim}2^{-j_{2}}\sigma_{j_{1},j_{2}}^{-2N} f^{\star %N}_{X\star\psi_{j_{1}}\star\psi_{j_{2}}}(2^{-j_{3}}\eta)
%=\underset{j_{2}\rightarrow\infty}{\lim}2^{-j_{2}}\sigma_{j_{1},j_{2}}^{-2N} f^{\star N}_{X\star\psi_{j_{1}}\star\psi_{j_{2}}}(0),
%\end{equation}
%where $\sigma_{j_{1},j_{2}}$ is the standard deviation of $X\star\psi_{j_{1}}\star\psi_{j_{2}}$.
%Furthermore,
%\begin{equation}\label{lemma_f_cov_limit_explicit}
%\underset{j_{2}\rightarrow\infty}{\lim}2^{-j_{2}}\sigma_{j_{1},j_{2}}^{-2N} f^{\star N}_{X\star\psi_{j_{1}}\star\psi_{j_{2}}}(0)
%= \frac{1}{2\pi} \int_{\mathbb{R}}\left\{\left[\int_{\mathbb{R}} \frac{|\Psi(\eta)|^2}{|\eta|^{1-(2\alpha+\beta)}}d\eta\right]^{-1}
%\left[ \int_{\mathbb{R}}e^{it\zeta} \frac{|\Psi(\zeta)|^{2}}{|\zeta|^{1-(2\alpha+\beta)}}
%d\zeta\right]\right\}^{N} dt,
%\end{equation}
%where the integral converges under the conditions $\Psi\in L^{2}$ and $0<2\alpha+\beta<1$.
\end{Lemma}
The proof of Lemma \ref{Lemma:fstarj1j2_limit} is in \ref{sec:proof:Lemma:fstarj1j2_limit}.

\begin{Proposition}\label{thm:third_order_Gaussian}
Let $G$ be a stationary Gaussian process satisfying Assumptions \ref{assumption:spectral} and let $\psi$ be a real-valued mother wavelet function satisfying Assumption \ref{Assumption:1:wavelet}.
For any  $\beta\in(0,1)$
and $j_{2}=j_{2}(j_{1})$ satisfying $\underset{j_{1}\rightarrow\infty}{\lim\inf}\ \frac{j_{2}(j_{1})}{j_{1}}>1$,
the rescaled random process
\begin{equation}\label{normalized_UG}
2^{\frac{j_{1}(\beta-1)}{2}}2^{\frac{j_{2}}{2}}\big|G\star\psi_{j_{1}}\big|\star\psi_{j_{2}}(2^{j_{2}}t),\  t\in \mathbb{R},
\end{equation}
converges to a Gaussian process $V$ when $j_{1}\rightarrow \infty$ in the finite dimensional distribution sense.
Moreover, the limiting process
$V$ has the following representation
\begin{align}\label{thm:thirdorder:V3}
V(t) = \kappa\int_{\mathbb{R}}
e^{i\lambda t}
\hat{\psi}(\lambda)W(d\lambda)
\end{align}
where
\begin{align}\label{def:kappa}
\kappa=\left(\sigma^{2}\overset{\infty}{\underset{\ell=2}{\sum}}
\gamma_{\ell}
C^{2}_{||,\ell}\right)^{\frac{1}{2}}.
\end{align}
The constants $\sigma^{2}$ and $\{\gamma_{\ell}\}_{\ell=2,3,...}$ are defined in Lemma \ref{lemma:var_S_T} and Lemma \ref{Lemma:fstarj1j2_limit}.
%Here, $\{C_{\ell}\}_{\ell= 2,3,...}$ are the coefficients of the Hermite polynomials expansion of the absolute value function.
\end{Proposition}
The proof of this proposition can be found in Section \ref{sec:proof:thm:third_order_Gaussian}.
By the definition of $\mathcal{D}_{j_{1},j_{2}}$ in (\ref{D_true}),
\begin{align}
2^{\frac{j_{1}(\beta-1)}{2}}2^{\frac{j_{2}}{2}}|X\star \psi_{j_{1}}|\star \psi_{j_{2}} = 2^{\frac{j_{1}(\beta-1)}{2}}2^{\frac{j_{2}}{2}}\Big[|C_{A,1}G\star \psi_{j_{1}}|\star \psi_{j_{2}} + \mathcal{D}_{j_{1},j_{2}}\Big].
\end{align}
By Assumption \ref{Assumption:passlimit} and Proposition \ref{lemma:estimateED}, for all $\beta\in(0,1)$,
\begin{align}\notag
&\mathbb{E}\left[|2^{\frac{j_{1}(\beta-1)}{2}}2^{\frac{j_{2}}{2}}\mathcal{D}_{j_{1},j_{2}}(t)|^{2}\right]
\\\notag\leq &M 2^{j_{1}(\beta-1)}2^{j_{2}}\mathbb{E}\left[|\widetilde{\mathcal{D}}_{j_{1},j_{2}}(t)|^{2}\right]
\\\notag\leq&2^{j_{1}(\beta-1)}2^{j_{2}}\times O\left(j_{1}2^{-j_{2}}+ 2^{-j_{1}\beta}2^{-j_{2}\beta} +2^{j_{1}(1-2\beta)}2^{-j_{2}}\right)
\\\label{normalizedED}=&
O\left(j_{1}2^{j_{1}(\beta-1)}+2^{-j_{1}}2^{j_{2}(1-\beta)} +2^{-j_{1}\beta}\right)
\end{align}
when $j_{1},j_{2}\rightarrow\infty$ with $\underset{j_{1}\rightarrow\infty}{\lim\inf}\ \frac{j_{2}(j_{1})}{j_{1}}>1$.
If we additionally impose the condition that $\underset{j_{1}\rightarrow\infty}{\lim\sup}\ \frac{j_{2}(j_{1})}{j_{1}}<\frac{1}{1-\beta}$, then (\ref{normalizedED}) implies that
\begin{align}
\underset{\begin{subarray}{c}j_{1}\rightarrow\infty \\ j_{2}=j_{2}(j_{1})\end{subarray}}{\lim}\mathbb{E}\left[|2^{\frac{j_{1}(\beta-1)}{2}}2^{\frac{j_{2}}{2}}\mathcal{D}_{j_{1},j_{2}}(t)|^{2}\right]=0.
\end{align}
By Slutsky's argument \cite[p. 6]{leonenko1999limit} and the continuous mapping theorem,
Proposition \ref{lemma:estimateED} and Proposition \ref{thm:third_order_Gaussian}
lead to the final result of this work.
\begin{Theorem}\label{thm5a}
Let $X=A(G)$ be a subordinated Gaussian process satisfying Assumptions \ref{assumption:GH}-\ref{assumption:spectral} and let $\psi$ be a real-valued mother wavelet function satisfying Assumption \ref{Assumption:1:wavelet}.
Furthermore, suppose that Assumption \ref{Assumption:passlimit} holds.
For any  $\beta\in(0,1)$ and $j_{2}=j_{2}(j_{1})$ satisfying
$1<\underset{j_{1}\rightarrow\infty}{\lim\inf}\ \frac{j_{2}(j_{1})}{j_{1}}\leq \underset{j_{1}\rightarrow\infty}{\lim\sup}\ \frac{j_{2}(j_{1})}{j_{1}}<\frac{1}{1-\beta}$,
\begin{equation}
2^{\frac{j_{1}(\beta-1)}{2}}2^{\frac{j_{2}}{2}}U[j_{1},j_{2}]X(2^{j_{2}}\cdot) \overset{d}{\Rightarrow}|C_{A,1}||V(\cdot)|
\end{equation}
in the finite dimensional distribution sense when $j_{1}\rightarrow\infty$.
The process
$V$ is Gaussian and has the representation (\ref{thm:thirdorder:V3}).
%Here, $C_{A}= \int_{\mathbb{R}} A(z)H_{1}(z)\frac{1}{\sqrt{2\pi}}e^{-\frac{z^{2}}{2}}dz$ is the Hermite coefficient of $A$.
\end{Theorem}

\section{Proofs of Proposition \ref{lemma:estimateED} and Proposition \ref{thm:third_order_Gaussian}}\label{sec:proof}
%In the below, we provide the proof of Theorem 1 and Theorem 2.
%In the context henceforth, the notation $\overset{d}{\Rightarrow}$ denotes the
%convergence of random processes
%in the sense of finite-dimensional
%distributions.
In the below, we employ  the diagram
method (see, \cite{breuer1983central} or  \cite[p.72]{ivanov2012statistical}) to prove Proposition \ref{lemma:estimateED} and Proposition \ref{thm:third_order_Gaussian}. Given integers $\ell_{1},\ldots,\ell_{p}$, where $p\geq2$,
a graph $\Gamma$ with $\ell_{1}+\cdots+\ell_{p}$ vertices is
called a complete diagram
of order ($\ell_{1},\ldots,\ell_{p}$) if:\\
(i) the set of vertices $V$ of the graph $\Gamma$ is of the form $V=\overset{p}{\underset{i=1}{\cup}}W_{i}$,
where $W_{i}=\{(i,\ell):1\leq \ell\leq \ell_{i}\}$ is the $i$-th level of the graph $\Gamma$;
\\
(ii) each vertex is of degree 1, that is, each vertex is just an endpoint of an edge;
\\
(iii) if $((i,\ell),\ (i^{'},\ell^{'}))\in \Gamma$, then $i\neq i^{'}$, that is, the edges of the graph $\Gamma$
 connect only different levels.

\begin{figure}
\centering
\subfigure[][Regular diagram
($\# A_{1,2}=3$, $\# A_{1,3} =0$, $\# A_{1,4} =0$, $\# A_{2,3} =0$, $\# A_{2,4} =0$, $\# A_{3,4}=4$)]{\label{regular_diagram}\includegraphics[height=0.20\textheight,width=0.425\textwidth]{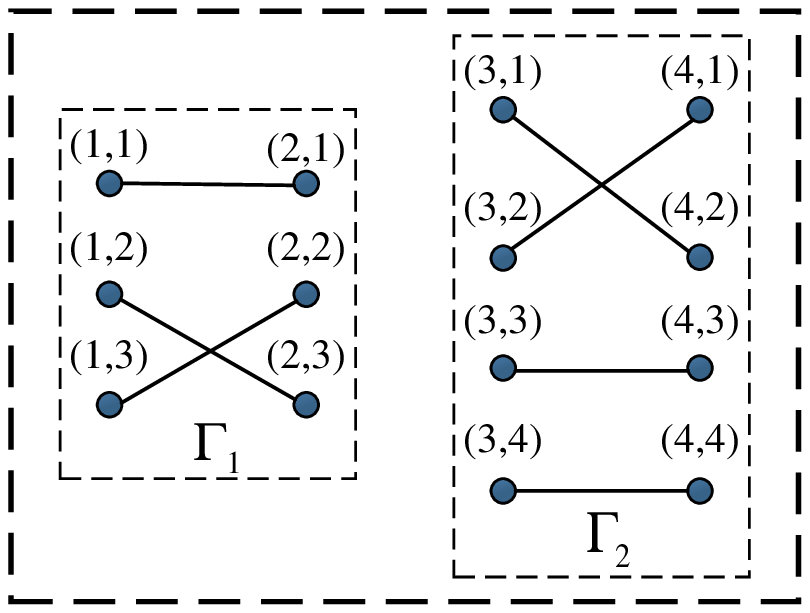}}
\hspace{0.2cm}
\subfigure[][Non-regular diagram
($\# A_{1,2}=2$, $\# A_{1,3} =0$, $\# A_{1,4} =1$, $\# A_{2,3} =1$, $\# A_{2,4} =0$, $\# A_{3,4}=3$)]{\label{nonregular_diagram}\includegraphics[height=0.20\textheight,width=0.425\textwidth]{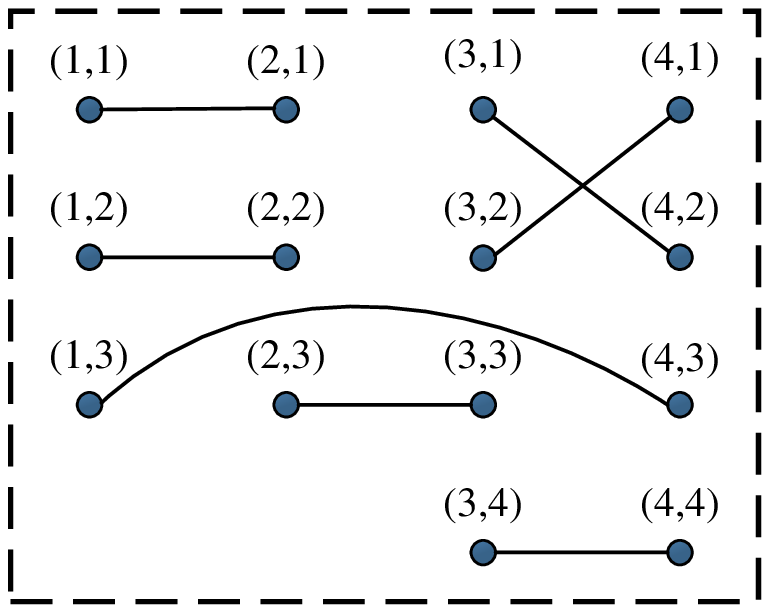}}
\caption{Illustration of regular and non-regular diagrams of order $(3,3,4,4)$. For integers $i<i'$, we denote the number of edges connecting the $i$-th and
$i'$-th levels by  $\# A_{i,i'}$.}
\label{diagram}
\end{figure}

Let $\mathrm{T}=\mathrm{T}(\ell_{1},\ldots,\ell_{p})$ be a set of
complete diagrams of order ($\ell_{1},\ldots,\ell_{p}$).
Denote by $E(\Gamma)$  the set of edges of the graph $\Gamma\in
\mathrm{T}$.
For the edge $e=((i,\ell),\
(i^{'},\ell^{'}))\in E(\Gamma)$ with $i<i^{'}$, $1\leq
\ell\leq \ell_{i}$ and $1\leq \ell^{'}\leq \ell_{i^{'}}$, we set
$d_{1}(e)=i$ and $d_{2}(e)=i^{'}$.
For any $i,i^{'}\in \{1,\ldots,p\}$ with $i<i^{'}$, let $A_{i,i^{'}}$ be the set of edges, which connect the $i$th and $i^{'}$th layers,
that is,
\begin{equation}\label{notationAB}
A_{i,i^{'}}=\big\{e\in E(\Gamma)\mid d_{1}(e)=i,\
d_{2}(e)=i^{'}\big\}.
\end{equation}
From the definition of $d_{1}(\cdot)$ and $d_{2}(\cdot)$, we know that $A_{i,i^{'}}= \phi$ for $i\geq i^{'}$.
Let $\# A_{i,i^{'}}$ be the cardinality of $A_{i,i^{'}}$, i.e., the number of edges in
$A_{i,i^{'}}$ (see Figure \ref{diagram}).
We call a diagram $\Gamma$ to be {\it
regular} if its levels can be split into pairs in such a manner
that no edge connects the levels belonging to different pairs (see Figure \ref{regular_diagram}).
Denote by $\mathrm{T}^{*}=\mathrm{T}^{*}(\ell_{1},\ldots,\ell_{p})$ the
set of all regular diagrams in $\mathrm{T}$.
If
$\Gamma\in \mathrm{T}^{*}$,
then $p$ is even and $\Gamma$ can be divided into $p/2$ sub-diagrams
(denoted by
$\Gamma_{1},\ldots,\Gamma_{p/2}$), which cannot be separated
again; in this case, we naturally define $d_{1}(\Gamma_{r})\equiv
d_{1}(e)$ and $d_{2}(\Gamma_{r})\equiv d_{2}(e)$ for any $e\in
E(\Gamma_{r}),\ r=1,\ldots,p/2$. We denote $\# E(\Gamma)$
(resp. $\# E(\Gamma_{r})$) the number of edges belonging to
the specific diagram $\Gamma$ (resp. the sub-diagram
$\Gamma_{r}$).

%%%%%%%%%%%%%%%%%%%%%%%%%%%%%%%%%%%%%%%%%%%%%%%%%%%%%%%%%%%%%%%%%%%%%%%%%%%%

\subsection{Proof of Proposition \ref{lemma:estimateED}}\label{sec:proof:lemma:estimateED}
\begin{figure}[h]
\centering
\subfigure[][$\# A_{1,2} = 3$, $\# A_{3,4} = 3$]
{\includegraphics[scale=0.5]{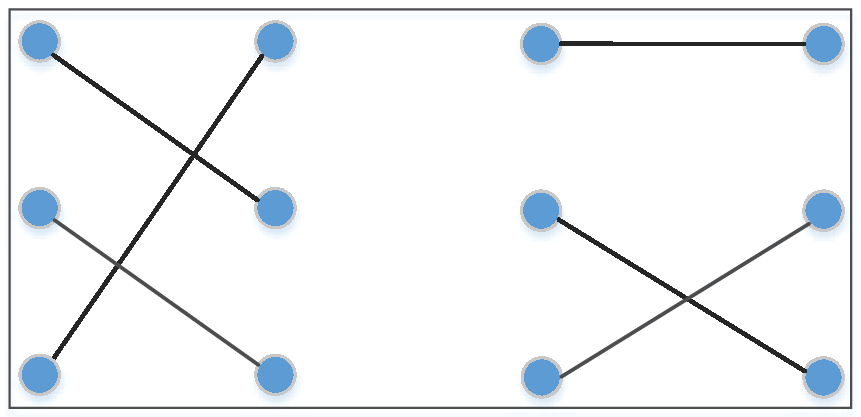}}
\subfigure[][$\# A_{1,3} = 3$, $\# A_{2,4} = 3$]
{\includegraphics[scale=0.5]{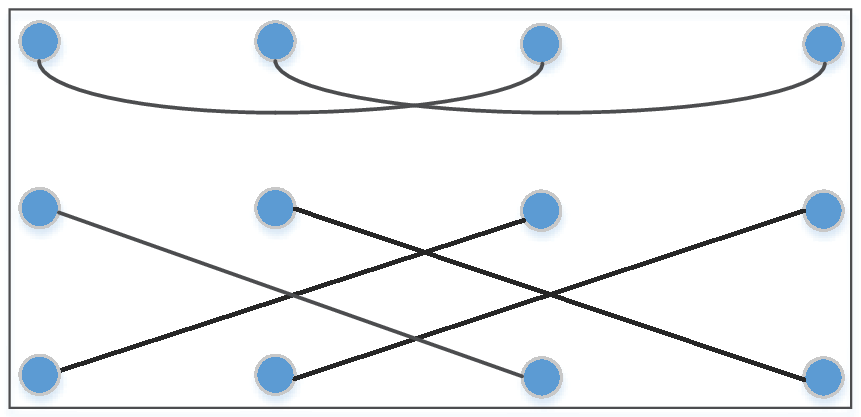}}
\subfigure[][$\# A_{1,4} = 3$, $\# A_{2,3} = 3$]
{\includegraphics[scale=0.5]{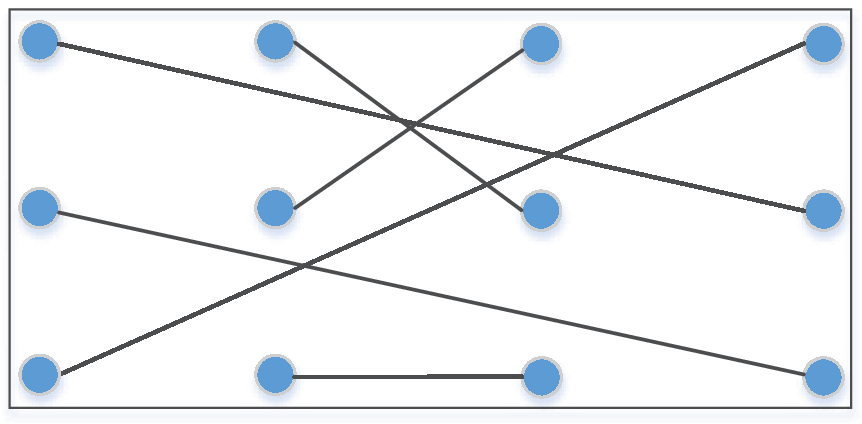}}
\caption{Illustration of three regular diagrams belonging to $\mathrm{T}^{*}_{1}$, $\mathrm{T}^{*}_{2}$, and $\mathrm{T}^{*}_{3}$, respectively.
}
\label{fig:3type_regular_diagram}
\end{figure}

%Because $\widetilde{\mathcal{D}}_{j_{1},j_{2}}$ is a stationary process, $\mathbb{E}[\widetilde{\mathcal{D}}_{j_{1},j_{2}}(t)^{2}] = %\mathbb{E}[\widetilde{\mathcal{D}}_{j_{1},j_{2}}(0)^{2}]$.
%Hence, it suffices to estimate $\mathbb{E}[\widetilde{\mathcal{D}}_{j_{1},j_{2}}(0)^{2}]$.
By denoting
%\begin{equation*}
$$Y_{j_{1}} = G\star\psi_{j_{1}}/\sigma_{G\star\psi_{j_{1}}},$$
%\end{equation*}
$\widetilde{\mathcal{D}}_{j_{1},j_{2}}(0)$ can be rewritten as
\begin{align}\label{def:D(t)}
\widetilde{\mathcal{D}}_{j_{1},j_{2}}(0)=\int_{\mathbb{R}}\textup{sign}\left(C_{A,1}\right)\textup{sign}\left(Y_{j_{1}}(\tau)\right)
\left[\overset{\infty}{\underset{\ell=2}{\sum}}
\frac{C_{A,\ell}}{\sqrt{\ell!}}H_{\ell}\left(G\right)\star \psi_{j_{1}}(\tau)\right]\psi_{j_{2}}(-\tau)d\tau.
\end{align}
Hence,
\begin{align}\notag
\mathbb{E}[\widetilde{\mathcal{D}}_{j_{1},j_{2}}(0)^2]=& \mathbb{E}\int_{\mathbb{R}}\int_{\mathbb{R}}\psi_{j_{2}}(-\tau)\psi_{j_{2}}(-\eta)
\textup{sign}\left(Y_{j_{1}}(\tau)\right)
\textup{sign}\left(Y_{j_{1}}(\eta)\right)
\\\notag&\times\left[\overset{\infty}{\underset{k=2}{\sum}}
\frac{C_{A,k}}{\sqrt{k!}}\int_{\mathbb{R}}H_{k}\left(G(s)\right) \psi_{j_{1}}(\tau-s)ds\right]
\left[\overset{\infty}{\underset{\ell=2}{\sum}}
\frac{C_{A,\ell}}{\sqrt{\ell!}}\int_{\mathbb{R}}H_{\ell}\left(G(t)\right) \psi_{j_{1}}(\eta-t)dt\right]
d\tau
d\eta
\\\notag=&
\overset{\infty}{\underset{k=2}{\sum}}\overset{\infty}{\underset{\ell=2}{\sum}}
\frac{C_{A,k}}{\sqrt{k!}}\frac{C_{A,\ell}}{\sqrt{\ell!}}
\int_{\mathbb{R}}\int_{\mathbb{R}}\int_{\mathbb{R}}\int_{\mathbb{R}}\psi_{j_{2}}(-\tau)\psi_{j_{2}}(-\eta)\psi_{j_{1}}(\tau-s)\psi_{j_{1}}(\eta-t)
\\\label{ED0_v1}&\times \mathbb{E}\left[H_{k}\left(G(s)\right)H_{\ell}\left(G(t)\right) \textup{sign}\left(Y_{j_{1}}(\tau)\right)
\textup{sign}\left(Y_{j_{1}}(\eta)\right)\right]
ds\
dt\
d\tau\
d\eta.
\end{align}
Because $\textup{sign}(\cdot)$ belongs to the Gaussian Hilbert space, it has the Hermite expansion
\begin{align}
\textup{sign}\left(y\right) = \overset{\infty}{\underset{\ell=1}{\sum}}\frac{C_{S,\ell}}{\sqrt{\ell!}}H_{\ell}(y),
\end{align}
where $\{C_{S,\ell}\}_{\ell=1}^{\infty}$ are the Hermite coefficients of $\textup{sign}(\cdot)$ with $C_{S,1}\neq0$.
Hence, (\ref{ED0_v1}) can be rewritten as
\begin{align}\label{ED}
\mathbb{E}[\widetilde{\mathcal{D}}_{j_{1},j_{2}}(0)^2]=&
\overset{\infty}{\underset{\ell_{1}=2}{\sum}}\ \overset{\infty}{\underset{\ell_{2}=2}{\sum}}\
\overset{\infty}{\underset{\ell_{3}=1}{\sum}}\ \overset{\infty}{\underset{\ell_{4}=1}{\sum}}
\frac{C_{A,\ell_{1}}}{\sqrt{\ell_{1}!}}
\frac{C_{A,\ell_{2}}}{\sqrt{\ell_{2}!}}
\frac{C_{S,\ell_{3}}}{\sqrt{\ell_{3}!}}
\frac{C_{S,\ell_{4}}}{\sqrt{\ell_{4}!}}
\\\notag
&\times\int_{\mathbb{R}}\int_{\mathbb{R}}\int_{\mathbb{R}}\int_{\mathbb{R}}\psi_{j_{2}}(-\tau)\psi_{j_{2}}(-\eta)\psi_{j_{1}}(\tau-s)\psi_{j_{1}}(\eta-t)
I_{\ell_{1},\ell_{2},\ell_{3},\ell_{4}}
\ ds
\ dt
\ d\tau
\ d\eta,
\end{align}
where
\begin{align*}%\label{expectation_I}
I_{\ell_{1},\ell_{2},\ell_{3},\ell_{4}} =
\mathbb{E}\left[H_{\ell_{1}}\left(G(s)\right) H_{\ell_{2}}\left(G(t)\right)H_{\ell_{3}}\left(Y_{j_{1}}(\tau)\right)
H_{\ell_{4}}\left(Y_{j_{1}}(\eta)\right)\right].
\end{align*}
The expectation $I_{\ell_{1},\ell_{2},\ell_{3},\ell_{4}}$ can be expanded by the Feynman diagram technique \cite[Theorem 5.3]{major1981lecture}.
By considering
a complete diagram $\Gamma$
of order ($\ell_{1},\ell_{2},\ell_{3},\ell_{4}$), where $\ell_{1}$, $\ell_{2}\geq2$ and $\ell_{3}$, $\ell_{4}\geq1$,
$I_{\ell_{1},\ell_{2},\ell_{3},\ell_{4}}$ can be expressed as
\begin{align}\label{expectation_I_v2}
I_{\ell_{1},\ell_{2},\ell_{3},\ell_{4}} =\underset{\Gamma\in \mathrm{T}}{\sum}I_{\ell_{1},\ell_{2},\ell_{3},\ell_{4}}(\Gamma),
\end{align}
where
\begin{align}\notag
I_{\ell_{1},\ell_{2},\ell_{3},\ell_{4}}(\Gamma)=&
R^{\#A_{1,2}}_{G}(s-t) R^{\#A_{1,3}}_{G,Y_{j_{1}}}(s-\tau)R^{\#A_{1,4}}_{G,Y_{j_{1}}}(s-\eta)
R^{\#A_{2,3}}_{G,Y_{j_{1}}}(t-\tau)
\\\label{def:I}&\times R^{\#A_{2,4}}_{G,Y_{j_{1}}}(t-\eta)
R^{\#A_{3,4}}_{Y_{j_{1}}}(\tau-\eta).
\end{align}
Here, $R_{Y_{j_{1}}}$ is the autocovariance function of the Gaussian process $Y_{j_{1}}$.
The cross covariance function between $G$ and $Y_{j_{1}}$ is denoted by
$R_{G,Y_{j_{1}}}(t-\tau)=\mathbb{E}\left[G(t)Y_{j_{1}}(\tau)\right]$, which has the following representation:
\begin{align}\notag
\mathbb{E}\left[G(t)Y_{j_{1}}(\tau)\right] =&\frac{1}{\sigma_{G\star\psi_{j_{1}}}}
\mathbb{E}\left[\int_{\mathbb{R}}e^{it\lambda}\sqrt{f_{G}(\lambda)}W(d\lambda)  \overline{\int_{\mathbb{R}}e^{i\tau\zeta}\sqrt{f_{G}(\zeta)}\hat{\psi}(2^{j_{1}}\zeta)W(d\zeta)}\right]
\\\notag=& \int_{\mathbb{R}}e^{i(t-\tau)\lambda}f_{c}(\lambda)d\lambda,
\end{align}
where
\begin{equation}\label{cross_spectral}
f_{c}(\lambda) = \frac{1}{\sigma_{G\star\psi_{j_{1}}}}f_{G}(\lambda)\overline{\hat{\psi}(2^{j_{1}}\lambda)}
\end{equation}
is the cross spectral density function.
For each fixed $(\ell_{1},\ell_{2},\ell_{3},\ell_{4})$, the set of complete diagrams $\mathrm{T} = \mathrm{T}(\ell_{1},\ell_{2},\ell_{3},\ell_{4})$
can be expressed as
\begin{equation}\label{decompose_diagram}
\mathrm{T} = \mathrm{T}^{*}\cup \left(\mathrm{T}\setminus \mathrm{T}^{*}\right) =\mathrm{T}^{*}_{1}\cup \mathrm{T}^{*}_{2}\cup \mathrm{T}^{*}_{3}\cup \left(\mathrm{T}\setminus \mathrm{T}^{*}\right),
\end{equation}
where $\mathrm{T}^{*}$ consists of three types of regular diagrams $\mathrm{T}^{*}_{1}$, $\mathrm{T}^{*}_{2}$, and $\mathrm{T}^{*}_{3}$.
If $\Gamma\in \mathrm{T}^{*}_{1}$, then $\# A_{1,2} = \ell_{1} = \ell_{2}$, $\# A_{3,4} = \ell_{3} = \ell_{4}$
and $\# A_{i,i^{'}} = 0$ for $(i,i^{'}) \notin \{(1,2),\ (3,4)\}$.
If $\Gamma\in \mathrm{T}^{*}_{2}$, then $\# A_{1,3} = \ell_{1} = \ell_{3}$, $\# A_{2,4} = \ell_{2} = \ell_{4}$
and $\# A_{i,i^{'}} = 0$ for $(i,i^{'}) \notin \{(1,3),\ (2,4)\}$.
If $\Gamma\in \mathrm{T}^{*}_{3}$, then $\# A_{1,4} = \ell_{1} = \ell_{4}$, $\# A_{2,3} = \ell_{2} = \ell_{3}$
and $\# A_{i,i^{'}} = 0$ for $(i,i^{'}) \notin \{(1,4),\ (2,3)\}$.
Some examples are illustrated in Figure. \ref{fig:3type_regular_diagram}.

By (\ref{expectation_I_v2}) and (\ref{decompose_diagram}),  (\ref{ED}) can be written as
\begin{align}\notag
\mathbb{E}[\widetilde{\mathcal{D}}_{j_{1},j_{2}}(0)^2]=&
\overset{\infty}{\underset{\ell_{1}=2}{\sum}}\ \overset{\infty}{\underset{\ell_{2}=2}{\sum}}\
\overset{\infty}{\underset{\ell_{3}=1}{\sum}}\ \overset{\infty}{\underset{\ell_{4}=1}{\sum}}
\frac{C_{A,\ell_{1}}}{\sqrt{\ell_{1}!}}
\frac{C_{A,\ell_{2}}}{\sqrt{\ell_{2}!}}
\frac{C_{S,\ell_{3}}}{\sqrt{\ell_{3}!}}
\frac{C_{S,\ell_{4}}}{\sqrt{\ell_{4}!}}\
\left(\underset{\Gamma\in \mathrm{T}^{*}_{1}}{\sum} + \underset{\Gamma\in \mathrm{T}^{*}_{2}}{\sum}
+ \underset{\Gamma\in \mathrm{T}^{*}_{3}}{\sum} + \underset{\Gamma\in \mathrm{T}\setminus\mathrm{T}^{*}}{\sum}\right)
\\\notag
&\times\int_{\mathbb{R}}\int_{\mathbb{R}}\int_{\mathbb{R}}\int_{\mathbb{R}}\psi_{j_{2}}(-\tau)\psi_{j_{2}}(-\eta)\psi_{j_{1}}(\tau-s)
\psi_{j_{1}}(\eta-t)
I_{\ell_{1},\ell_{2},\ell_{3},\ell_{4}}(\Gamma)
\ ds
\ dt
\ d\tau
\ d\eta
\\\label{ED_v2}=& E_{1}+E_{2}+E_{3}+E_{4}.
\end{align}
In the following, the asymptotic behavior of $E_{1}$, $E_{2}$, $E_{3}$, and $E_{4}$
when $j_{1}\rightarrow\infty$ and $\underset{j_{1}\rightarrow\infty}{\lim\inf}\ \frac{j_{2}(j_{1})}{j_{1}}>1$
is analyzed. \\

\noindent{\bf $\bullet$ Asymptotic behavior of $E_{1}$:}

By the definition of $\mathrm{T}_{1}^{*}$ and the spectral representation of the covariance functions,
\begin{align}\notag
E_{1} =&
\overset{\infty}{\underset{k=2}{\sum}}
C_{A,k}^{2}
\overset{\infty}{\underset{\ell=1}{\sum}}
C_{S,\ell}^{2}
\int_{\mathbb{R}^{4}}\psi_{j_{2}}(-\tau)\psi_{j_{2}}(-\eta)\psi_{j_{1}}(\tau-s)\psi_{j_{1}}(\eta-t)
R^{k}_{G}(s-t)
R^{\ell}_{Y_{j_{1}}}(\tau-\eta) \ ds
\ dt
\ d\tau
\ d\eta
\\\notag=&\overset{\infty}{\underset{k=2}{\sum}}
C_{A,k}^{2}
\overset{\infty}{\underset{\ell=1}{\sum}}
C_{S,\ell}^{2}
\int_{\mathbb{R}^{2}}|\hat{\psi}(2^{j_{2}}(\lambda_{1}+\lambda_{2}))|^{2}|\hat{\psi}(2^{j_{1}}\lambda_{1})|^{2}
f^{\star k}_{G}(\lambda_{1})
f^{\star \ell}_{Y_{j_{1}}}(\lambda_{2}) \ d\lambda_{1}
\ d\lambda_{2}.
\end{align}
By changing of variables $\omega_{1} = 2^{j_{1}}\lambda_{1}$ and $\omega_{2} = 2^{j_{1}}(\lambda_{1}+\lambda_{2})$,
\begin{align}\label{E1_v1}
E_{1}=& 2^{-2j_{1}}\overset{\infty}{\underset{k=2}{\sum}}
C_{A,k}^{2}
\overset{\infty}{\underset{\ell=1}{\sum}}
C_{S,\ell}^{2}\
B_{k,\ell}(j_{1},j_{2}),
\end{align}
where
\begin{align}
B_{k,\ell}(j_{1},j_{2})=\int_{\mathbb{R}^{2}}|\hat{\psi}(2^{j_{2}-j_{1}}\omega_{2})|^{2}|\hat{\psi}(\omega_{1})|^{2}
f^{\star k}_{G}(2^{-j_{1}}\omega_{1})
f^{\star \ell}_{Y_{j_{1}}}(2^{-j_{1}}(\omega_{2}-\omega_{1})) \ d\omega_{1}\ d\omega_{2}.
\end{align}
From Lemma \ref{lemma:convolution}, we know that the singularity strength of $f^{\star k}_{G}$
and $f^{\star\ell}_{Y_{j_{1}}}$ at the origin decreases along with $k$ and $\ell$, respectively.
Hence among $\{B_{k,\ell}\}_{k\geq 2,\ell\geq 1}$,
the integrand in $B_{2,1}$ has the strongest singularity at the origin.
Because $\overset{\infty}{\underset{k=2}{\sum}}
C_{A,k}^{2}<\infty$ and $\overset{\infty}{\underset{\ell=1}{\sum}}
C_{S,\ell}^{2}<\infty$, (\ref{E1_v1}) implies that
\begin{align}\label{E1_v2}
E_{1} = O\left(
2^{-2j_{1}}
B_{2,1}(j_{1},j_{2})\right).
\end{align}
Lemma \ref{lemma:convolution} shows that if $2\beta<1$ and $\omega_{1}\neq0$,
%\begin{equation*}
%f^{\star 2}_{G}(\omega) = O\left(|\omega|^{2\beta-1}\right)
%\end{equation*}
%when $|\omega|\rightarrow 0$, which implies that
\begin{equation}\label{E1_v2_component1}
f^{\star 2}_{G}(2^{-j_{1}}\omega_{1}) =O\left( 2^{j_{1}(1-2\beta)}|\omega_{1}|^{2\beta-1}\right)
\end{equation}
as $j_{1}\rightarrow \infty$.
On the other hand, by
\begin{equation*}
f_{Y_{j_{1}}}(\omega) = \frac{1}{\sigma_{G\star\psi_{j_{1}}}^{2}}f_{G}(\omega)|\hat{\psi}(2^{j_{1}}\omega)|^{2}
=\frac{1}{\sigma_{G\star\psi_{j_{1}}}^{2}}\frac{C_{G}(\omega)}{|\omega|^{1-\beta}}|\hat{\psi}(2^{j_{1}}\omega)|^{2},
\end{equation*}
and Lemma \ref{lemma:var_S_T}, that is, $\sigma_{G\star\psi_{j_{1}}}^{2}= O\left(2^{-j_{1}\beta}\right)$ when $j_{1}\rightarrow \infty$,
we have
\begin{equation*}
f_{Y_{j_{1}}}(2^{-j_{1}}\omega) = O\left(2^{ j_{1}}|\omega|^{\beta-1}|\hat{\psi}(\omega)|^{2}\right)
\end{equation*}
when $j_{1}\rightarrow\infty$.
It implies that
\begin{equation}\label{E1_v2_component2}
f_{Y_{j_{1}}}(2^{-j_{1}}(\omega_{2}-\omega_{1})) =O\left( 2^{ j_{1}}|\omega_{2}-\omega_{1}|^{(2\alpha+\beta)-1}|C_{\hat{\psi}}(\omega_{2}-\omega_{1})|^{2}\right)
\end{equation}
when $j_{1}$ is sufficiently large.
Under the situation $2\beta<1$,
by substituting (\ref{E1_v2_component1}) and (\ref{E1_v2_component2}) into (\ref{E1_v2}),
we obtain
\begin{align}\notag
 E_{1} =&O\left(
2^{-2j_{1}}
\int_{\mathbb{R}^{2}}|C_{\hat{\psi}}(2^{j_{2}-j_{1}}\omega_{2})|^{2}|2^{j_{2}-j_{1}}\omega_{2}|^{2\alpha} |C_{\hat{\psi}}(\omega_{1})|^{2}|\omega_{1}|^{2\alpha}\left[2^{j_{1}(1-2\beta)}|\omega_{1}|^{2\beta-1}\right]\right.
\\\notag&\left.\times\left[2^{ j_{1}}|\omega_{2}-\omega_{1}|^{(2\alpha+\beta)-1}|C_{\hat{\psi}}(\omega_{2}-\omega_{1})|^{2}\right]d\omega_{1}d\omega_{2}\right)
\\\label{E1_v2_new_v2}=&O\left(
2^{2j_{2}\alpha}2^{-2j_{1}(\alpha+\beta)}
\int_{\mathbb{R}}|C_{\hat{\psi}}(2^{j_{2}-j_{1}}\omega_{2})|^{2}|\omega_{2}|^{2\alpha}h_{1}(\omega_{2})d\omega_{2}\right),
%\right.
%\\&\times\left.\left[\int_{\mathbb{R}}|C_{\hat{\psi}}(\omega_{1})|^{2}|\omega_{1}|^{(2\alpha+2\beta)-1}
%|C_{\hat{\psi}}(\omega_{2}-\omega_{1})|^{2}
%|\omega_{2}-\omega_{1}|^{(2\alpha+\beta)-1}
%d\omega_{1}\right]d\omega_{2}\right).
\end{align}
where
\begin{equation*}
h_{1}(\omega_{2}):= \int_{\mathbb{R}}|C_{\hat{\psi}}(\omega_{1})|^{2}|\omega_{1}|^{(2\alpha+2\beta)-1}
|C_{\hat{\psi}}(\omega_{2}-\omega_{1})|^{2}
|\omega_{2}-\omega_{1}|^{(2\alpha+\beta)-1}
d\omega_{1}.
\end{equation*}
Under Assumption \ref{Assumption:1:wavelet}, Lemma \ref{lemma:convolution} implies that
the function $h_{1}$ in (\ref{E1_v2_new_v2})
is integrable and continuous.
By the change of variables,
\begin{align}
\underset{
\begin{subarray}{c}
j_{1}\rightarrow\infty
\\
j_{2}=j_{2}(j_{1})
\end{subarray}
}{\lim}2^{(j_{2}-j_{1})(2\alpha+1)}\int_{\mathbb{R}}|C_{\hat{\psi}}(2^{j_{2}-j_{1}}\omega_{2})|^{2}|\omega_{2}|^{2\alpha}
h_{1}(\omega_{2})d\omega_{2}=h_{1}(0)\|\hat{\psi}\|_{2}^{2},
\end{align}
where $j_{2}=j_{2}(j_{1})$ satisfies $\underset{j_{1}\rightarrow\infty}{\lim\inf}\ \frac{j_{2}(j_{1})}{j_{1}}>1$.
Therefore, for $\beta\in(0,\frac{1}{2})$,
\begin{align}\label{final_E1_b<05}
E_{1} =O\left(2^{j_{1}(1-2\beta)}2^{-j_{2}}\right)
\end{align}
when $j_{1}\rightarrow\infty$ and
$\underset{j_{1}\rightarrow\infty}{\lim\inf}\ \frac{j_{2}(j_{1})}{j_{1}}>1$.
Similarly,
\begin{align}\label{final_E1_b>05}
E_{1} = \left\{\begin{array}{ll}
O\left(j_{1}2^{-j_{2}}\right) & \textup{for}\ \beta=\frac{1}{2},\\
O\left(2^{-j_{2}}\right) & \textup{for}\ \beta\in(\frac{1}{2},1)
\end{array}\right.
\end{align}
when $j_{1}\rightarrow\infty$ and
$\underset{j_{1}\rightarrow\infty}{\lim\inf}\ \frac{j_{2}(j_{1})}{j_{1}}>1$.

\noindent{\bf $\bullet$ Asymptotic behavior of $E_{2}$:}
From the definition of $\mathrm{T}_{2}^{*}$,
\begin{align}\notag
E_{2} =&
\overset{\infty}{\underset{k=2}{\sum}}
C_{A,k}C_{S,k}
\overset{\infty}{\underset{\ell=2}{\sum}}
C_{A,\ell}C_{S,\ell}
\int_{\mathbb{R}^{4}}\psi_{j_{2}}(-\tau)\psi_{j_{2}}(-\eta)\psi_{j_{1}}(\tau-s)\psi_{j_{1}}(\eta-t)
\\\notag &\times R^{k}_{G,{Y_{j_{1}}}}(s-\tau)
R^{\ell}_{G,Y_{j_{1}}}(t-\eta) \ ds
\ dt
\ d\tau
\ d\eta.
\end{align}
By changing of variables $(\tau,\eta,s,t) \rightarrow (\tau',\eta',\tau^{'}-s^{'},\eta^{'}-t^{'})$,
we know that the integral above vanishes because $\int_{\mathbb{R}}\psi_{j_{2}}(-\tau^{'})d\tau^{'} = 0.$
Hence, $E_{2} = 0.$\\

\noindent{\bf $\bullet$ Asymptotic behavior of $E_{3}$:}
From the definition of $\mathrm{T}_{3}^{*}$,
\begin{align}\notag
E_{3} =&
\overset{\infty}{\underset{k=2}{\sum}}
C_{A,k}C_{S,k}
\overset{\infty}{\underset{\ell=2}{\sum}}
C_{A,\ell}C_{S,\ell}
\int_{\mathbb{R}^{4}}\psi_{j_{2}}(-\tau)\psi_{j_{2}}(-\eta)\psi_{j_{1}}(\tau-s)\psi_{j_{1}}(\eta-t)
\\\notag &\times R^{k}_{G,{Y_{j_{1}}}}(s-\eta)
R^{\ell}_{G,Y_{j_{1}}}(t-\tau) \ ds
\ dt
\ d\tau
\ d\eta
\\=&\overset{\infty}{\underset{k=2}{\sum}}
C_{A,k}C_{S,k}
\overset{\infty}{\underset{\ell=2}{\sum}}
C_{A,\ell}C_{S,\ell}
\ \tilde{B}_{k,\ell},
\end{align}
where
\begin{align}\notag
\tilde{B}_{k,\ell} =& \int_{\mathbb{R}^{2}}|\hat{\psi}(2^{j_{2}}(\omega_{1}-\omega_{2}))|^{2}
\hat{\psi}(2^{j_{1}}\omega_{1})\hat{\psi}(2^{j_{1}}\omega_{2})f^{\star k}_{c}(\omega_{1})
f^{\star \ell}_{c}(\omega_{2})\ d\omega_{2}\ d\omega_{1}
\\\label{E3_origin_B}=& \int_{\mathbb{R}^{2}}|\hat{\psi}(2^{j_{2}}\lambda_{2})|^{2}
\hat{\psi}(2^{j_{1}}\lambda_{1})\hat{\psi}(2^{j_{1}}(\lambda_{1}-\lambda_{2}))f^{\star k}_{c}(\lambda_{1})
f^{\star \ell}_{c}(\lambda_{1}-\lambda_{2})\ d\lambda_{2}\ d\lambda_{1}
\end{align}
and
\begin{align}\label{recall_fc}
f_{c}(\lambda) =\frac{2^{\alpha j_{1}}}{\sigma_{G\star\psi_{j_{1}}}}
C_{G}(\lambda)\overline{C_{\hat{\psi}}(2^{j_{1}}\lambda)}|\lambda|^{(\alpha+\beta)-1},\ \lambda\in \mathbb{R},
\end{align}
is the cross-spectral density function of $G$ and $Y_{j_{1}}$
initially defined in (\ref{cross_spectral}).
By the change of variables $\lambda_{1}\rightarrow 2^{-j_{1}}\lambda_{1}$
and $\lambda_{2}\rightarrow 2^{-j_{1}}\lambda_{2}$,
(\ref{E3_origin_B}) can be rewritten as
\begin{align*}%\label{E3_case1}
\tilde{B}_{k,\ell} = 2^{-2j_{1}}\int_{\mathbb{R}^{2}}|\hat{\psi}(2^{j_{2}-j_{1}}\lambda_{2})|^{2}
\hat{\psi}(\lambda_{1})\hat{\psi}(\lambda_{1}-\lambda_{2})f^{\star k}_{c}(2^{-j_{1}}\lambda_{1})
f^{\star \ell}_{c}(2^{-j_{1}}(\lambda_{1}-\lambda_{2}))\ d\lambda_{2}\ d\lambda_{1}.
\end{align*}
Because the iterative convolution reduces the singularity strength of $\{f^{\star \ell}_{c}\}_{\ell =2,3,...}$ at the origin,
the asymptotic behavior of $E_{3}$ is dominated by  $\tilde{B}_{2,2}$, that is,
\begin{align}\label{E3_case1}
E_{3} =&O\left(
2^{-2j_{1}}\int_{\mathbb{R}^{2}}|\hat{\psi}(2^{j_{2}-j_{1}}\lambda_{2})|^{2}
\hat{\psi}(\lambda_{1})\hat{\psi}(\lambda_{1}-\lambda_{2})f^{\star 2}_{c}(2^{-j_{1}}\lambda_{1})
f^{\star 2}_{c}(2^{-j_{1}}(\lambda_{1}-\lambda_{2}))\ d\lambda_{2}\ d\lambda_{1}
\right)\end{align}
when $j_{1}\rightarrow\infty$ and $\underset{j_{1}\rightarrow\infty}{\lim\inf}\ \frac{j_{2}(j_{1})}{j_{1}}>1$.
By (\ref{recall_fc}), we first have
\begin{align}\notag
f^{\star 2}_{c}(2^{-j_{1}}\lambda) =&\left[\frac{2^{\alpha j_{1}}}{\sigma_{G\star\psi_{j_{1}}}}
\right]^{2}\int_{\mathbb{R}}\frac{C_{G}(2^{-j_{1}}\lambda-\eta)\overline{C_{\hat{\psi}}(2^{j_{1}}(2^{-j_{1}}\lambda-\eta))}}
{|2^{-j_{1}}\lambda-\eta|^{1-(\alpha+\beta)}}
\frac{C_{G}(\eta)\overline{C_{\hat{\psi}}(2^{j_{1}}\eta)}}{|\eta|^{1-(\alpha+\beta)}}d\eta
\\\label{E3_fc_cov}=
&\left[\frac{2^{\alpha j_{1}}}{\sigma_{G\star\psi_{j_{1}}}}
\right]^{2}2^{j_{1}(1-2\alpha-2\beta)}\int_{\mathbb{R}}\frac{C_{G}(2^{-j_{1}}(\lambda-\eta))
\overline{C_{\hat{\psi}}(\lambda-\eta)}}{|\lambda-\eta|^{1-(\alpha+\beta)}}
\frac{C_{G}(2^{-j_{1}}\eta)\overline{C_{\hat{\psi}}(\eta)}}{|\eta|^{1-(\alpha+\beta)}}d\eta.
\end{align}
Because $\sigma_{G\star\psi_{j_{1}}}^{-2}=O( 2^{j_{1}\beta})$ when $j_{1}\rightarrow \infty$ and $\beta<1$, from (\ref{E3_fc_cov}),
we have
\begin{align}\label{E3_fc_cov_v2}
f^{\star 2}_{c}(2^{-j_{1}}\lambda) = &O\left(2^{j_{1}(1-\beta)}h_{2}(\lambda)\right),
\end{align}
where
\begin{equation}\label{E3_def_K}
h_{2}(\lambda) = \int_{\mathbb{R}}
\overline{C_{\hat{\psi}}(\lambda-\eta)}|\lambda-\eta|^{(\alpha+\beta)-1}
\overline{C_{\hat{\psi}}(\eta)}|\eta|^{(\alpha+\beta)-1}d\eta,\ \lambda\in \mathbb{R},
\end{equation}
is a continuous and integrable function because $\alpha\geq1$ (Assumption \ref{Assumption:1:wavelet}).
By applying (\ref{E3_fc_cov_v2}) to the functions $f^{\star 2}_{c}(2^{-j_{1}}\lambda_{1})$ and
$f^{\star 2}_{c}(2^{-j_{1}}(\lambda_{1}-\lambda_{2}))$ in (\ref{E3_case1}),
we get
\begin{align}\notag
E_{3} =&O\Big(
2^{-2j_{1}}\int_{\mathbb{R}^{2}}|C_{\hat{\psi}}(2^{j_{2}-j_{1}}\lambda_{2})|^{2}|2^{j_{2}-j_{1}}\lambda_{2}|^{2\alpha}
C_{\hat{\psi}}(\lambda_{1})|\lambda_{1}|^{\alpha}C_{\hat{\psi}}(\lambda_{1}-\lambda_{2})|\lambda_{1}-\lambda_{2}|^{\alpha}
\\\notag&\times\left[2^{j_{1}(1-\beta)}h_{2}(\lambda_{1})\right]
\left[2^{j_{1}(1-\beta)}h_{2}(\lambda_{1}-\lambda_{2})\right]d\lambda_{2}d\lambda_{1}\Big)
\\\label{E3_case1_v2}=&
O\Big(2^{-j_{1}(2\alpha+2\beta)}2^{2j_{2}\alpha}
\int_{\mathbb{R}}
|C_{\hat{\psi}}(2^{j_{2}-j_{1}}\lambda_{2})|^{2}|\lambda_{2}|^{2\alpha}h_{3}(\lambda_{2})d\lambda_{2}\Big),
\end{align}
where
\begin{align*}
h_{3}(\lambda_{2}):=
\int_{\mathbb{R}}h_{2}(\lambda_{1})C_{\hat{\psi}}(\lambda_{1})|\lambda_{1}|^{\alpha}h_{2}(\lambda_{1}-\lambda_{2})C_{\hat{\psi}}(\lambda_{1}-\lambda_{2})
|\lambda_{1}-\lambda_{2}|^{\alpha}d\lambda_{1}
\end{align*}
is an integrable and continuous function.
By the change of variables,  (\ref{E3_case1_v2}) implies that
\begin{align}\label{E3_C31_final}
E_{3} =&O\left(2^{j_{1}(1-2\beta)}2^{-j_{2}}\right)
\end{align}
when $\beta<1$, $j_{1}\rightarrow\infty$, and $\underset{j_{1}\rightarrow\infty}{\lim\inf}\ \frac{j_{2}(j_{1})}{j_{1}}>1$.

\begin{figure}[h]
\centering
\subfigure[][]
{\includegraphics[scale=0.55]{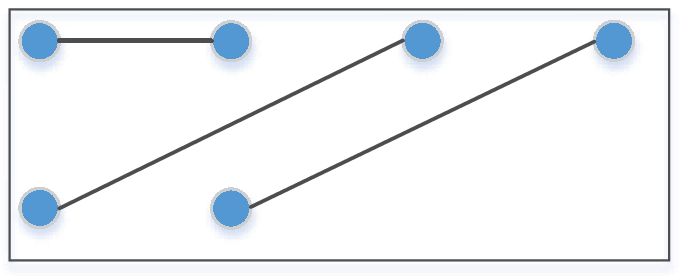}}
\subfigure[][]
{\includegraphics[scale=0.55]{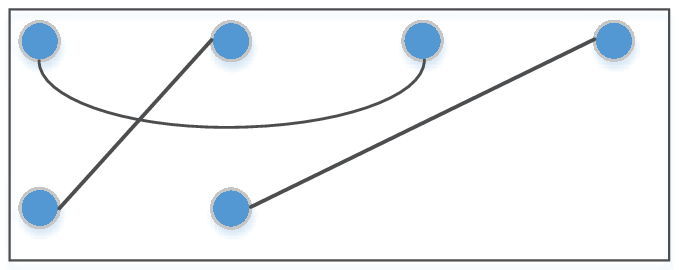}}
\subfigure[][]
{\includegraphics[scale=0.55]{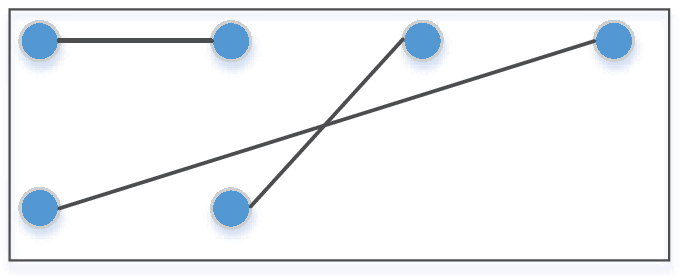}}
\subfigure[][]
{\includegraphics[scale=0.55]{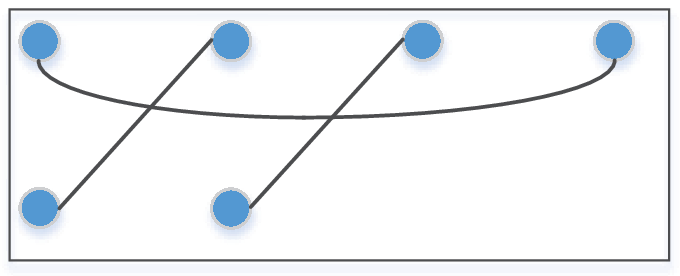}}
\caption{Illustration of non-regular diagrams
}
\label{fig:nonregular_diagram}
\end{figure}

\noindent{\bf $\bullet$ Asymptotic behavior of $E_{4}$:}
From the definition of $E_{4}$ in (\ref{ED_v2})
\begin{equation}\label{E4:sum}
E_{4} = \frac{C^{2}_{A,2}}{2}C^{2}_{S,1}\Big[E_{4,a,b}+E_{4,c,d}\Big]+F,
\end{equation}
where
\begin{equation*}
E_{4,a,b} =
2\int_{\mathbb{R}^{4}}\psi_{j_{2}}(-\tau)\psi_{j_{2}}(-\eta)\psi_{j_{1}}(\tau-s)\psi_{j_{1}}(\eta-t)
R_{G}(s-t) R_{G,Y_{j_{1}}}(s-\tau)
R_{G,Y_{j_{1}}}(t-\eta)
ds
dt
d\tau
d\eta
\end{equation*}
is contributed by the non-regular diagrams (a) and (b)
in Figure \ref{fig:nonregular_diagram}
and
\begin{equation*}
E_{4,c,d} =
2\int_{\mathbb{R}^{4}}\psi_{j_{2}}(-\tau)\psi_{j_{2}}(-\eta)\psi_{j_{1}}(\tau-s)\psi_{j_{1}}(\eta-t)
R_{G}(s-t) R_{G,Y_{j_{1}}}(s-\eta)
R_{G,Y_{j_{1}}}(t-\tau)
ds
dt
d\tau
d\eta
\end{equation*}
is contributed by the non-regular diagrams (c) and (d)
in Figure \ref{fig:nonregular_diagram}.
The third component $F$ in (\ref{E4:sum}) is contributed by non-regular diagrams,
each of which contains more than three edges.
From the definition of $I_{\ell_{1},\ell_{2},\ell_{3},\ell_{4}}(\Gamma)$ in (\ref{def:I}), we know that the more edges $\Gamma$ has, the faster decay
$R^{\# A_{1,2}}_{Y_{j_{1}}}$, $R^{\# A_{1,3}}_{Y_{j_{1}},Y_{j_{1},j_{2}}}$, $R^{\# A_{1,4}}_{Y_{j_{1}},Y_{j_{1},j_{2}}}$,
$R^{\# A_{2,3}}_{Y_{j_{1}},Y_{j_{1},j_{2}}}$, or $R^{\# A_{2,4}}_{Y_{j_{1}},Y_{j_{1},j_{2}}}$ have.
Hence, for understanding the asymptotic behavior of $E_{4}$ when $j_{1},j_{2}\rightarrow \infty$,
it suffices to analyze $E_{4,a,b}$ and $E_{4,c,d}$.
By the spectral representation for the covariance functions $R_{G}$ and $R_{G,Y_{j_{1}}}$
(see also (\ref{cross_spectral})),
\begin{align}\notag
E_{4,a,b} =&
2\int_{\mathbb{R}^{3}}|\hat{\psi}(2^{j_{2}}\lambda_{1})|^{2}\hat{\psi}(2^{j_{1}}(\lambda_{1}+\lambda_{2}))\hat{\psi}(2^{j_{1}}(\lambda_{3}-\lambda_{1}))
f_{G}(\lambda_{1}) f_{c}(\lambda_{2})
f_{c}(\lambda_{3})
d\lambda_{1}
d\lambda_{2}
d\lambda_{3}
\\\notag
=&2^{-j_{2}-2j_{1}}\int_{\mathbb{R}^{3}}|\hat{\psi}(\lambda_{1})|^{2}\hat{\psi}(\frac{\lambda_{1}}{2^{j_{2}-j_{1}}}+\lambda_{2})
\hat{\psi}(\lambda_{3}-\frac{\lambda_{1}}{2^{j_{2}-j_{1}}})
 f_{G}(\frac{\lambda_{1}}{2^{j_{2}}}) f_{c}(\frac{\lambda_{2}}{2^{j_{1}}})
f_{c}(\frac{\lambda_{3}}{2^{j_{1}}})
d\lambda_{1}
d\lambda_{2}
d\lambda_{3}.
\end{align}
It implies that when $j_{1}\rightarrow\infty$ and $\underset{j_{1}\rightarrow\infty}{\lim\inf}\ \frac{j_{2}(j_{1})}{j_{1}}>1$,
\begin{align}\label{E41}
E_{4,a,b} =O\left(2^{-j_{2}-2j_{1}}\Big[\int_{\mathbb{R}}|\hat{\psi}(\lambda_{1})|^{2}
f_{G}(\frac{\lambda_{1}}{2^{j_{2}}})d\lambda_{1}\Big]
\Big[\int_{\mathbb{R}}
\hat{\psi}(\lambda_{2})f_{c}(\frac{\lambda_{2}}{2^{j_{1}}})d\lambda_{2}\Big]^{2}\right).
\end{align}
For the first integral in (\ref{E41}), we have
\begin{align}\label{E41a}
\int_{\mathbb{R}}|\hat{\psi}(\lambda)|^{2}
f_{G}(2^{-j_{2}}\lambda)d\lambda
=
\int_{\mathbb{R}}|\hat{\psi}(\lambda)|^{2}
\frac{C_{G}(2^{-j_{2}}\lambda)}{|2^{-j_{2}}\lambda|^{1-\beta}}d\lambda_{1}
=
O\left(2^{j_{2}(1-\beta)}\right)
\end{align}
when $j_{1}\rightarrow\infty$.
For the second integral in (\ref{E41}), by (\ref{recall_fc}),
\begin{align}\label{fc_asymp_pre}
f_{c}(2^{-j_{1}}\lambda) =\frac{2^{j_{1}(1-\beta)}}{\sigma_{G\star\psi_{j_{1}}}}
C_{G}(2^{-j_{1}}\lambda)\overline{C_{\hat{\psi}}(\lambda)}|\lambda|^{(\alpha+\beta)-1}.
\end{align}
Because $\sigma_{G\star \psi_{j_{1}}}=O( 2^{-j_{1}\beta/2})$ when $j_{1}\rightarrow\infty$ and $\beta<1$,
(\ref{fc_asymp_pre}) implies that
\begin{align}\label{fc_asymp}
f_{c}(2^{-j_{1}}\lambda)
=O\left(2^{j_{1}(1-\frac{\beta}{2})}\overline{C_{\hat{\psi}}(\lambda)}|\lambda|^{(\alpha+\beta)-1}\right).
\end{align}
Hence,
\begin{align}\label{E41b}
\Big[\int_{\mathbb{R}}
\hat{\psi}(\lambda)f_{c}(2^{-j_{1}}\lambda)d\lambda\Big]^{2}
=O\left(
2^{j_{1}(2-\beta) }
\Big[\int_{\mathbb{R}}
\frac{|\hat{\psi}(\lambda)|^{2}}{|\lambda|^{1-\beta}}d\lambda\Big]^{2}\right).
\end{align}
By substituting (\ref{E41a}) and (\ref{E41b}) into (\ref{E41}) yields
\begin{align}\label{limit_E41}
 E_{4,a,b}= O\left(2^{-j_{1}\beta} 2^{-j_{2}\beta} \right)
\end{align}
when $\beta<1$, $j_{1}\rightarrow\infty$ and $\underset{j_{1}\rightarrow\infty}{\lim\inf}\ \frac{j_{2}(j_{1})}{j_{1}}>1$.

Next, we analyze $E_{4,c,d}$.
By the spectral representation for the covariance functions $R_{G}$ and $R_{G,Y_{j_{1}}}$ again,
we have
\begin{align}\notag
E_{4,c,d} =&
2\int_{\mathbb{R}^{3}}|\hat{\psi}(2^{j_{2}}(\lambda_{1}+\lambda_{2}-\lambda_{3}))|^{2}\hat{\psi}(2^{j_{1}}(\lambda_{1}+\lambda_{2}))
\hat{\psi}(2^{j_{1}}(\lambda_{3}-\lambda_{1}))
f_{G}(\lambda_{1}) f_{c}(\lambda_{2})
f_{c}(\lambda_{3})
d\lambda_{1}
d\lambda_{2}
d\lambda_{3}
\\\notag
=&
2\int_{\mathbb{R}^{3}}|\hat{\psi}(2^{j_{2}}\omega_{1})|^{2}\hat{\psi}(2^{j_{1}}\omega_{2})\hat{\psi}(2^{j_{1}}\omega_{3})
f_{G}(\omega_{2}-\omega_{3}-\omega_{1}) f_{c}(\omega_{3}+\omega_{1})
f_{c}(\omega_{2}-\omega_{1})
d\omega_{1}
d\omega_{2}
d\omega_{3}
\\\label{E42}
=&2^{1-3j_{1}}\int_{\mathbb{R}^{3}}|\hat{\psi}(2^{-j_{1}+j_{2}}\omega_{1})|^{2}\hat{\psi}(\omega_{2})\hat{\psi}(\omega_{3})
f_{G}(\frac{\omega_{2}-\omega_{3}-\omega_{1}}{2^{j_{1}}}) f_{c}(\frac{\omega_{3}+\omega_{1}}{2^{j_{1}}})
f_{c}(\frac{\omega_{2}-\omega_{1}}{2^{j_{1}}})
d\omega_{1}
d\omega_{2}
d\omega_{3}.
\end{align}
%It implies that when $j_{3} = \delta j_{2}$ for some $\delta>1$ and $j_{2}\rightarrow\infty$,
%\begin{align}\label{E42}
%E_{4,2} \asymp 2^{-j_{3}-2j_{2}}\int_{\mathbb{R}^{2}}\Psi(\omega_{2})\Psi(\omega_{3})
%f_{Y_{j_{1}}}(\frac{\omega_{2}-\omega_{3}}{2^{j_{2}}})
%f_{c}(\frac{\omega_{3}}{2^{j_{2}}})
%f_{c}(\frac{\omega_{2}}{2^{j_{2}}})
%d\omega_{2}
%d\omega_{3}.
%\end{align}
By Assumption \ref{assumption:spectral}, for $\omega_{2}-\omega_{3}-\omega_{1}\neq0$,
\begin{align}\label{asy_f_Y1}
f_{G}(\frac{\omega_{2}-\omega_{3}-\omega_{1}}{2^{j_{1}}}) =O\left(2^{j_{1}(1-\beta)}|\omega_{2}-\omega_{3}-\omega_{1}|^{\beta-1}\right)
\end{align}
when $j_{1}$ is large enough.
By (\ref{fc_asymp}),
\begin{align}\label{asy_f_Y12}
f_{c}(\frac{\omega_{3}+\omega_{1}}{2^{j_{1}}})
f_{c}(\frac{\omega_{2}-\omega_{1}}{2^{j_{1}}})=O\left(
2^{j_{1}(2-\beta) }\frac{\overline{\hat{\psi}(\omega_{3}+\omega_{1})}\overline{\hat{\psi}(\omega_{2}-\omega_{1})}}{|\omega_{3}+\omega_{1}|^{1-\beta}|\omega_{2}-\omega_{1}|^{1-\beta}}
\right).\end{align}
By substituting (\ref{asy_f_Y1}) and (\ref{asy_f_Y12}) into (\ref{E42}), we get
\begin{align}\notag
E_{4,c,d}=& O\Big(
2^{-2j_{1}\beta}\int_{\mathbb{R}}
|\hat{\psi}(2^{-j_{1}+j_{2}}\omega_{1})|^{2}
h_{4}(\omega_{1})
d\omega_{1}\Big),
\end{align}
where
\begin{align}
h_{4}(\omega_{1}) = \int_{\mathbb{R}^{2}} \frac{\hat{\psi}(\omega_{2})
\hat{\psi}(\omega_{3})\overline{\hat{\psi}(\omega_{3}+\omega_{1})}\overline{\hat{\psi}(\omega_{2}-\omega_{1})}}{|\omega_{2}-\omega_{3}-\omega_{1}|^{1-\beta}|\omega_{3}+\omega_{1}|^{1-\beta}|\omega_{2}-\omega_{1}|^{1-\beta}}
d\omega_{2}
d\omega_{3}.
\end{align}
By the change of variables,
\begin{align}\notag
E_{4,c,d}=& O\Big(
2^{-2j_{1}\beta}2^{j_{1}-j_{2}}\int_{\mathbb{R}}
|\hat{\psi}(\omega_{1})|^{2}
h_{4}(2^{j_{1}-j_{2}}\omega_{1})
d\omega_{1}\Big).
\end{align}
Note that
\begin{align}\notag
\underset{\begin{subarray}{c}j_{1}\rightarrow\infty\\j_{2}=j_{2}(j_{1})\end{subarray}}
{\lim}h_{4}(2^{j_{1}-j_{2}}\omega_{1}) =& \int_{\mathbb{R}} \int_{\mathbb{R}}\frac{|\hat{\psi}(\omega_{2})|^{2}
|\hat{\psi}(\omega_{3})|^{2}}{|\omega_{2}-\omega_{3}|^{1-\beta}|\omega_{3}|^{1-\beta}|\omega_{2}|^{1-\beta}}
d\omega_{2}
d\omega_{3}
\\\label{app:Riesz}=&\int_{\mathbb{R}}\frac{|\hat{\psi}(\omega_{2})|^{2}}{|\omega_{2}|^{1-\beta}}\left[ \int_{\mathbb{R}}\frac{1
}{|\omega_{2}-\omega_{3}|^{1-\beta}}\frac{|\hat{\psi}(\omega_{3})|^{2}}{|\omega_{3}|^{1-\beta}}
d\omega_{3}
\right]d\omega_{2}.
\end{align}
Under Assumption \ref{Assumption:1:wavelet}, $\frac{|\hat{\psi}(\cdot)|^{2}}{|\cdot|^{1-\beta}}$ is continuous and integrable,
so $\frac{|\hat{\psi}(\cdot)|^{2}}{|\cdot|^{1-\beta}}\in L^{\frac{1}{\beta+\varepsilon}}(\mathbb{R})$ for all $\varepsilon\in(0,1-\beta)$.
By \cite[p. 119]{stein1970singular},
\begin{align}
h_{5}(\omega_{2}):=\int_{\mathbb{R}}\frac{1
}{|\omega_{2}-\omega_{3}|^{1-\beta}}\frac{|\hat{\psi}(\omega_{3})|^{2}}{|\omega_{3}|^{1-\beta}}
d\omega_{3}\in L^{\frac{1}{\varepsilon}}(\mathbb{R},d\omega_{2})
\end{align}
and $\|h_{5}\|_{L^{\frac{1}{\varepsilon}}(\mathbb{R})}\leq A_{\frac{1}{\beta+\varepsilon},\frac{1}{\varepsilon}} \|\frac{|\hat{\psi}(\cdot)|^{2}}{|\cdot|^{1-\beta}}\|_{L^{\frac{1}{\beta+\varepsilon}}(\mathbb{R})}$
where $A_{\frac{1}{\beta+\varepsilon},\frac{1}{\varepsilon}}$ is an universal constant.
Hence, by the H$\ddot{\textup{o}}$lder inequality, the limit in (\ref{app:Riesz}) is finite.
Therefore, we get
\begin{align}\label{E42s}
E_{4,c,d} =O\left(2^{j_{1}(1-2\beta)}2^{-j_{2}}\right).
\end{align}
Combining (\ref{limit_E41}) and (\ref{E42s}) yields
\begin{align}\label{limit_E4}
E_{4} =O\left( 2^{-j_{1}\beta}2^{-j_{2}\beta} +2^{j_{1}(1-2\beta)}2^{-j_{2}}\right).
\end{align}
when $\beta<1$, $j_{1}\rightarrow\infty$ and $\underset{j_{1}\rightarrow\infty}{\lim\inf}\ \frac{j_{2}(j_{1})}{j_{1}}>1$.

Finally, by combining (\ref{final_E1_b<05}), (\ref{final_E1_b>05}), (\ref{E3_C31_final}), and (\ref{limit_E4}),
we get
\begin{equation*}
E_{1}+E_{2}+E_{3}+E_{4} = \left\{\begin{array}{ll} O\left( 2^{-j_{1}\beta}2^{-j_{2}\beta} +2^{j_{1}(1-2\beta)}2^{-j_{2}}\right)\ & \ \textup{for} \ \beta\in(0,\frac{1}{2}),\\
O\left(j_{1}2^{-j_{2}}+ 2^{-j_{1}/2}2^{-j_{2}/2} +2^{-j_{2}}\right)\ & \ \textup{for} \ \beta=\frac{1}{2},
\\
O\left(2^{-j_{2}}+ 2^{-j_{1}\beta}2^{-j_{2}\beta} +2^{j_{1}(1-2\beta)}2^{-j_{2}}\right)\ & \ \textup{for} \ \beta\in(\frac{1}{2},1),
\end{array}\right.
\end{equation*}
when $j_{1}\rightarrow\infty$ and $\underset{j_{1}\rightarrow\infty}{\lim\inf}\ \frac{j_{2}(j_{1})}{j_{1}}>1$.
It leads to (\ref{estimate_D_combine}).

%%%%%%%%%%%%%%%%%%%%%%%%%%%%%%%%%%%%%%%%%%%%%%%%%%%%%%%%%%%%%%%%%%%%%%%%%%%

\subsection{Proof of Proposition \ref{thm:third_order_Gaussian}}\label{sec:proof:thm:third_order_Gaussian}
%The method of the following proof can
%be traced back to \cite{BreuerMajor}.

For any $M\in \mathbb{N}$ and any set of real numbers
$\{a_{1},a_{2},\ldots,a_{M}\}$, denote
\begin{align}
\xi_{j_{1},j_{2}}=
\overset{M}{\underset{k=1}{\sum}}a_{k}2^{\frac{j_{1}(\beta-1)}{2}}2^{\frac{j_{2}}{2}}|G\star\psi_{j_{1}}|\star\psi_{j_{2}}(2^{j_{2}}t_{k}),
\end{align}
where $t_{1},\ldots,t_{M}\in \mathbb{R}$
are arbitrary.
We divide $\xi_{j_{1},j_{2}}$ into two parts as follows
\begin{equation}\label{proof_3rd_ST:1}
\xi_{j_{1},j_{2}} = \xi_{j_{1},j_{2},\leq N}+\xi_{j_{1},j_{2},>N},
\end{equation}
where $N\in\{2,4,6,...\}$,
\begin{equation}\label{proof_3rd_ST:main}
\xi_{j_{1},j_{2},\leq N} = \sigma_{G\star \psi_{j_{1}}}2^{\frac{j_{1}(\beta-1)}{2}}2^{\frac{j_{2}}{2}}
\overset{M}{\underset{k=1}{\sum}}a_{k}
\sum_{\ell=2}^{N}
\frac{C_{||,\ell}}{\sqrt{\ell!}}H_{\ell}(\frac{G\star\psi_{j_{1}}}{\sigma_{G\star \psi_{j_{1}}}})\star\psi_{j_{2}}(2^{j_{2}}t_{k}),
\end{equation}
and
\begin{equation}\label{proof_3rd_ST:res}
\xi_{j_{2},j_{3},> N} = \sigma_{G\star \psi_{j_{1}}}2^{\frac{j_{1}(\beta-1)}{2}}2^{\frac{j_{2}}{2}}
\overset{M}{\underset{k=1}{\sum}}a_{k}
\sum_{\ell=N+1}^{\infty}
\frac{C_{||,\ell}}{\sqrt{\ell!}}H_{\ell}(\frac{G\star\psi_{j_{1}}}{\sigma_{G\star \psi_{j_{1}}}})\star\psi_{j_{2}}(2^{j_{2}}t_{k}).
\end{equation}
Note that $\{C_{||,\ell}\}_{\ell=0}^{\infty}$ are the Hermite coefficients of $|\cdot|$.
Because $|\cdot|$ is an even function, $C_{||,\ell}=0$ if $\ell$ is odd.
%Thus, the Gaussian process $X\star\psi_{j_{1}}\star\psi_{j_{2}}/\sigma_{j_{1},j_{2}}$ has variance one.
%{\color{blue} For any fixed $\delta>1$, we will analyze the double scaling limits
%of the statistical moments of $\xi_{j_{2},j_{3},\leq N}$ and $\xi_{j_{2},j_{3},> N}$, where $j_{3} = \delta j_{2}$ and $j_{2}\rightarrow\infty$.
%Although $j_{2}$ and $j_{3}$ are restricted to tend to infinity simultaneously,
%the constraint $\delta>1$ leads to that $2^{j_{3}}$ tends to infinity faster than $2^{j_{2}}$. Hence, we will first take the limit %$j_{3}\rightarrow \infty$ in the following computation when $j_{2}$ and $j_{3}$ appear
%in the equations at the same time.
%}

For any even integer $N\geq2$,
we define a truncated version of the limiting process as follows
\begin{align}
V_{\leq N}(t)=
\kappa_{N}\int_{\mathbb{R}}
e^{i\lambda t}
\hat{\psi}(\lambda)W(d\lambda),
\end{align}
where
\begin{align}\label{df:kappa_truncated_o3}
\kappa_{N}=\left[\sigma^{2}\overset{N}{\underset{\ell=2}{\sum}}
\gamma_{\ell}
C^{2}_{||,\ell}\right]^{\frac{1}{2}}.
\end{align}
The constant $\sigma^{2}=\underset{j_{1}\rightarrow\infty}{\lim}2^{j_{1}\beta}\sigma^{2}_{G\star \psi_{j_{1}}}$ is defined in Lemma \ref{lemma:var_S_T}. The constants $\{\gamma_{\ell}\}_{\ell=2}^{\infty}$ are defined in Lemma \ref{Lemma:fstarj1j2_limit} and
satisfy $\gamma_{2}\geq \gamma_{4} \geq \gamma_{6}  \geq \cdots\geq 0.$
Also note that
$\overset{\infty}{\underset{\ell=0}{\sum}}
C^{2}_{||,\ell}= \mathbb{E}[|G|^{2}] = 1$. Hence, $\kappa_{N}$ converges to $\kappa$ when $N\rightarrow\infty$, where $\kappa$ is defined in (\ref{def:kappa}).
For $j_{2}=j_{2}(j_{1})$ satisfying $\underset{j_{1}\rightarrow\infty}{\lim\inf}\ \frac{j_{2}(j_{1})}{j_{1}}>1$, we will prove that
\begin{equation}\label{proof:3rd_ST:3claim}
\begin{array}{ll}
\textup{(a)}\ \underset{N\rightarrow\infty}{\lim}\ \underset{\begin{subarray}{c}j_{1}\rightarrow\infty\\
j_{2}=j_{2}(j_{1})\end{subarray}}{\lim} \mathbb{E}[\xi_{j_{1},j_{2},>N}^2]=0;\\
%\textup{(b)}\ \textup{The limit of}\ \kappa_{3,N}\ \textup{exists when}\ N\rightarrow\infty;\\
\textup{(b)}\ \xi_{j_{1},j_{2},\leq N}\ \textup{converges in distribution to} \overset{M}{\underset{k=1}{\sum}}a_{k}V_{\leq N}(t_{k})\ \textup{as}\ j_{1}\rightarrow\infty.
\end{array}
\end{equation}
By Slutsky's argument, (a), (b) and the convergence of $\kappa_{N}$ when $N\rightarrow\infty$ imply that
$\xi_{j_{1},j_{2}}$ converges in distribution to $\overset{M}{\underset{k=1}{\sum}}a_{k}V(t_{k})$ when $j_{1}\rightarrow\infty$
and $\underset{j_{1}\rightarrow\infty}{\lim\inf}\ \frac{j_{2}(j_{1})}{j_{1}}>1$.

\noindent{\it Proof of (a):}
In the following, we prove that $\underset{N\rightarrow\infty}{\lim}\ \underset{\begin{subarray}{c}j_{1}\rightarrow\infty\\
j_{2}=j_{2}(j_{1})\end{subarray}}{\lim}\mathbb{E}[\xi_{j_{1},j_{2},>N}^2]=0$.
By (\ref{expect_H_Gpsi}),
\begin{align}\notag
\mathbb{E}[\xi_{j_{1},j_{2},>N}^2]=&
 2^{j_{1}(\beta-1)}2^{j_{2}}\overset{M}{\underset{k=1}{\sum}}\overset{M}{\underset{h=1}{\sum}}a_{k}a_{h}
\sum_{\ell=N+1}^{\infty}C^{2}_{||,\ell}\ \sigma_{G\star \psi_{j_{1}}}^{-2(\ell-1)}
\\\notag&\times\int_{\mathbb{R}}\int_{\mathbb{R}} R^{\ell}_{G\star\psi_{j_{1}}}(s-s^{'})\psi_{j_{2}}(2^{j_{2}}t_{k}-s)
\psi_{j_{2}}(2^{j_{2}}t_{h}-s^{'})dsds^{'}
\\\notag=& 2^{j_{1}(\beta-1)}2^{j_{2}}\overset{M}{\underset{k=1}{\sum}}\overset{M}{\underset{h=1}{\sum}}a_{k}a_{h}
\sum_{\ell=N+1}^{\infty}C^{2}_{||,\ell}\int_{\mathbb{R}} e^{i2^{j_{2}}\lambda(t_{k}-t_{h})}\sigma_{G\star \psi_{j_{1}}}^{-2(\ell-1)}f^{\star\ell}_{G\star\psi_{j_{1}}}(\lambda)
|\hat{\psi}(2^{j_{2}}\lambda)|^2
d\lambda
\\\label{third_order_Gaussian_proof_a}=&
2^{j_{1}(\beta-1)}\int_{\mathbb{R}}
|\overset{M}{\underset{k=1}{\sum}}a_{k}e^{i\lambda t_{k}}|^{2}\left[\sum_{\ell=N+1}^{\infty}C^{2}_{||,\ell}\ \sigma_{G\star \psi_{j_{1}}}^{-2(\ell-1)}f^{\star\ell}_{G\star\psi_{j_{1}}}(\frac{\lambda}{2^{j_{2}}})
\right]|\hat{\psi}(\lambda)|^2
d\lambda,
\end{align}
where the last equality follows from applying the monotone convergence theorem to change the order of $\overset{\infty}{\underset{\ell=N+1}{\sum}}$ and integration.
From the nonnegativity of the spectral densities,
\begin{align}\notag
f_{G\star\psi_{j_{1}}}^{\star\ell}(\lambda)
= \int_{\mathbb{R}}f_{G\star\psi_{j_{1}}}^{\star(\ell-1)}(\lambda-\zeta)f_{G\star\psi_{j_{1}}}(\zeta)d\zeta
\leq \|f_{G\star\psi_{j_{1}}}^{\star(\ell-1)}\|_{\infty} \sigma_{G\star \psi_{j_{1}}}^{2}.
\end{align}
It implies that
\begin{align}\label{third_order_conv_decrease}
\sigma_{G\star \psi_{j_{1}}}^{-2\ell}\|f_{G\star\psi_{j_{1}}}^{\star\ell}\|_{\infty}
\leq \sigma_{G\star \psi_{j_{1}}}^{-2N}\|f_{G\star\psi_{j_{1}}}^{\star N}\|_{\infty}
\end{align}
for all $\ell\geq N+1$.
Hence, the second summation in (\ref{third_order_Gaussian_proof_a}) can be estimated as follows
\begin{align}\label{proof:3rd_ST:(a)}
\sum_{\ell=N+1}^{\infty}C^{2}_{||,\ell}\ \sigma_{G\star \psi_{j_{1}}}^{-2(\ell-1)}f_{G\star\psi_{j_{1}}}^{\star\ell}(2^{-j_{2}}\lambda)
\leq \sigma_{G\star \psi_{j_{1}}}^{2-2N} \|f_{G\star\psi_{j_{1}}}^{\star N}\|_{\infty}\sum_{\ell=N+1}^{\infty}C^{2}_{||,\ell}.
\end{align}
where $\|f_{G\star\psi_{j_{1}}}^{\star N}\|_{\infty} = f_{G\star\psi_{j_{1}}}^{\star N}(0)<\infty$ by Lemma \ref{Lemma:fstarj1j2_limit}.
Thanks to (\ref{proof:3rd_ST:(a)}), we get an upper bound for (\ref{third_order_Gaussian_proof_a}):
\begin{align}\label{proof_third_order_a}
\mathbb{E}[\xi_{j_{1},j_{2},>N}^2]
\leq
\left[2^{j_{1}\beta}\sigma_{G\star \psi_{j_{1}}}^{2}\right]
\left[2^{-j_{1}}\sigma_{G\star \psi_{j_{1}}}^{-2N}f^{\star N}_{G\star\psi_{j_{1}}}(0)\right]
\left[\sum_{\ell=N+1}^{\infty}C^{2}_{||,\ell}
\right]
\int_{\mathbb{R}}
|\overset{M}{\underset{k=1}{\sum}}a_{k}e^{i\lambda t_{k}}|^{2}|\hat{\psi}(\lambda)|^2
d\lambda.
\end{align}
By (\ref{lemma:varS}), we know that
\begin{equation}\label{proof_third_order_a1}
\underset{j_{1}\rightarrow\infty}{\lim} 2^{j_{1}\beta}\sigma_{G\star \psi_{j_{1}}}^{2} = C_{G}(0) \int_{\mathbb{R}} \frac{|\hat{\psi}(\eta)|^2}{|\eta|^{1-\beta}}d\eta.
\end{equation}
By Lemma \ref{Lemma:fstarj1j2_limit},
\begin{align}\label{proof_third_order_a2}
\gamma_{N} = \underset{j_{1}\rightarrow\infty}{\lim}2^{-j_{1}}\sigma_{G\star \psi_{j_{1}}}^{-2N} f^{\star N}_{G\star\psi_{j_{1}}}(0)
=\frac{1}{2\pi} \int_{\mathbb{R}}\left[\int_{\mathbb{R}} \frac{|\hat{\psi}(\eta)|^2}{|\eta|^{1-\beta}}d\eta\right]^{-N}
\left[ \int_{\mathbb{R}}e^{it\zeta} \frac{|\hat{\psi}(\zeta)|^{2}}{|\zeta|^{1-\beta}}
d\zeta\right]^{N} dt
\end{align}
and $\gamma_{N}\geq\gamma_{N+2}\geq0 $ for all $N\in\{2,4,...\}$.
By (\ref{proof_third_order_a}), (\ref{proof_third_order_a1}), (\ref{proof_third_order_a2}), and the fact $\sum_{\ell=0}^{\infty}C^{2}_{||,\ell}=1,$ we get
\begin{equation}\notag
\underset{N\rightarrow\infty}{\lim}\ \underset{\begin{subarray}{c}j_{1}\rightarrow\infty\\
j_{2}=j_{2}(j_{1})\end{subarray}}{\lim}\mathbb{E}[\xi_{j_{1},j_{2},>N}^2] = 0.
\end{equation}

%%%%%%%%%%%%%%%%%%%
\noindent{\it Proof of (b):}
We apply the method of moments \cite{breuer1983central} (see also \cite[Theorem 6.5]{dasgupta2008asymptotic}) to prove
the claim (b). Because the linear combination of $V_{\leq N}(t)$ has the Gaussian distribution, which can be uniquely determined by its moments, it suffices to prove that
\begin{align}\label{Markovmethod2_o3}
\underset{\begin{subarray}{c}j_{1}\rightarrow\infty\\
j_{2}=j_{2}(j_{1})\end{subarray}}{\lim}\ \mathbb{E}\xi_{j_{1},j_{2},\leq N}^{p}=
\left
\{\begin{array}{lr}0 &\textup{if}\ p\ \textup{is odd},
\\ (p-1)!!
\left\{\mathbb{E}\Big[\big(\overset{M}{\underset{k=1}{\sum}}a_{k}V_{\leq N}(t_{k})\big)^{2}\Big]\right\}^{p/2}
 &\textup{if}\ p\ \textup{is even}.
\end{array}
\right.
\end{align}
% as follows
%\begin{equation}
%\underset{N_{0}\rightarrow\infty}{\lim}\kappa_{\textup{trunc}} = \frac{1}{2\pi}\int_{\mathbb{R}}R_{|X\star \psi_{j_{1}}|}(\tau)d\tau.
%\end{equation}
From the definition of
$\xi_{j_{1},j_{2},\leq
N}$ in (\ref{proof_3rd_ST:main}),
\begin{align}\notag
\mathbb{E}\xi_{j_{1},j_{2},\leq N}^{p}=&2^{\frac{j_{1}(\beta-1)}{2}p}\ 2^{\frac{j_{2}}{2}p}\ \sigma_{G\star \psi_{j_{1}}}^{p}
\overset{M}{\underset{{k_{1},\cdots,k_{p}=1}}{\sum}}\
\overset{N}{\underset{\ell_{1},\cdots,\ell_{p}=2}{\sum}}
\Big[\overset{p}{\underset{i=1}{\prod}}a_{k_{i}} \frac{C_{||,\ell_{i}}}{\sqrt{\ell_{i}!}}
\Big]
\\\label{proofthm:weaksmall6_o3a}
\times&
\int_{\mathbb{R}^{p}}
\Big[
\overset{p}{\underset{i=1}{\prod}}
2^{-j_{2}}\psi(\frac{2^{j_{2}}t_{k_{i}}-s_{i}}{2^{j_{2}}})\Big]
\Big[\mathbb{E}\overset{p}{\underset{i=1}{\prod}} H_{\ell_{i}}\Big(\frac{G\star \psi_{j_{1}}(s_{i})}{\sigma_{G\star \psi_{j_{1}}}}\Big)\Big]
ds_{1}\cdots ds_{p}.
\end{align}
By considering a complete diagram of order ($\ell_{1},\ldots,\ell_{p}$), where $p\geq2$ and $2\leq \ell_{i}\leq N$ for $1\leq i\leq p$,
%The definition of a complete diagram can be found at the beginning of this section.
(\ref{proofthm:weaksmall6_o3a}) can be rewritten as
\begin{align}\label{proofthm:weaksmall6_o3}
\mathbb{E}\xi_{j_{1},j_{2},\leq N}^{p}&=\underset{(K,L)\in D_{p}}{\sum}\hspace{-0.35cm}C(K,L)\underset{\Gamma\in \mathrm{T}^{*}}{\sum}\hspace{-0.05cm}F_{\Gamma}(K,L,j_{1},j_{2})
+\hspace{-0.05cm}\underset{(K,L)\in D_{p}}{\sum}\hspace{-0.35cm}C(K,L)\underset{\Gamma\in \mathrm{T}\backslash \mathrm{T}^{*}}{\sum}\hspace{-0.15cm}F_{\Gamma}(K,L,j_{1},j_{2}),
\end{align}
where
\begin{align*}\notag
C(K,L)=&\overset{p}{\underset{i=1}{\prod}}a_{k_{i}}\ \frac{C_{||,\ell_{i}}}{\sqrt{\ell_{i}!}},
\end{align*}
\begin{align}\label{proofthm:weakK(J,L)_o3}
F_{\Gamma}(K,L,j_{1},j_{2})=&2^{\frac{j_{1}(\beta-1)}{2}p}\ 2^{\frac{j_{2}}{2}p}\ \sigma_{G\star \psi_{j_{1}}}^{p}
\int_{\mathbb{R}^{p}}
\Big[
\overset{p}{\underset{i=1}{\prod}}
2^{-j_{2}}\psi(\frac{2^{j_{2}}t_{k_{i}}-s_{i}}{2^{j_{2}}})\Big]
\\\notag&\times\Big[
\underset{e\in E(\Gamma)}{\prod}R_{Y_{j_{1}}}(s_{d_{1}(e)}-s_{d_{2}(e)})
\Big]ds_{1}\cdots ds_{p},
\end{align}
$$Y_{j_{1}}(s) = \frac{G\star \psi_{j_{1}}(s)}{\sigma_{G\star \psi_{j_{1}}}},\ s \in \mathbb{R},$$
and
\begin{equation*}
D_{p}=\{(K,L):K=(k_{1},\ldots,k_{p}),1\leq k_{i}\leq M,\
L=(\ell_{1},\ldots, \ell_{p}),\ 2\leq \ell_{i}\leq N,\ i=1,\ldots,p\}.
\end{equation*}
To prove (\ref{Markovmethod2_o3}),  by (\ref{proofthm:weaksmall6_o3}), it suffices to verify the following two claims:
\begin{align*}
\left\{
\begin{array}{l}
\textup{(b1)}\ \underset{\begin{subarray}{c}j_{1}\rightarrow\infty\\
j_{2}=j_{2}(j_{1})\end{subarray}}{\lim}
\underset{(K,L)\in D_{p}}{\sum}C(K,L)\underset{\Gamma\in \mathrm{T}^{*}}{\sum}F_{\Gamma}(K,L,j_{1},j_{2})
=(p-1)!!
\Big\{\mathbb{E}\Big[\big(\overset{M}{\underset{k=1}{\sum}}a_{k}V_{N}(t_{k})\big)^{2}\Big]\Big\}^{p/2},
\\
\textup{(b2)}\ \underset{\begin{subarray}{c}j_{1}\rightarrow\infty\\
j_{2}=j_{2}(j_{1})\end{subarray}}{\lim}\underset{(K,L)\in D_{p}}{\sum}C(K,L)\underset{\Gamma\in \mathrm{T}\backslash \mathrm{T}^{*}}{\sum}F_{\Gamma}(K,L,j_{1},j_{2})
=0.
\end{array}\right.
\end{align*}

\noindent{\it Proof of (b1):} If $\Gamma$ is a regular diagram in $\mathrm{T}^*(\ell_{1},\ldots,\ell_{p})$,
then $\Gamma$ has a unique
decomposition
$\Gamma=(\Gamma_{1},\ldots,\Gamma_{p/2})$, where
$\Gamma_{1}$,$\ldots$,$\Gamma_{p/2}$
cannot be further decomposed.
Accordingly, $F_{\Gamma}(K,L,j_{1},j_{2})$
can be rewritten as the following $p/2$ products
\begin{align}\notag%\label{proofthm:weaksmall7_o3}
F_{\Gamma}(K,L,j_{1},j_{2})=&
2^{\frac{j_{1}(\beta-1)}{2}p}\ 2^{\frac{j_{2}}{2}p}\ \sigma_{G\star \psi_{j_{1}}}^{p}\ \overset{p/2}{\underset{r=1}{\prod}}\int_{\mathbb{R}^{2}}
2^{-j_{2}}
\psi(\frac{2^{j_{2}}t_{k_{d_{1}(\Gamma_{r})}}-s}{2^{j_{2}}})
\\\notag&\times
2^{-2j_{2}} \psi(\frac{2^{j_{2}}t_{k_{d_{2}(\Gamma_{r})}}-s^{'}}{2^{j_{2}}})
R_{Y_{j_{1}}}^{\# E(\Gamma_{r})}(s-s^{'})
\ ds\ ds^{'}.
\end{align}
By
\begin{align*}
R_{Y_{j_{1}}}^{\# E(\Gamma_{r})}(s-s^{'})
=\int_{\mathbb{R}}e^{i\lambda(s-s^{'})}f_{Y_{j_{1}}}^{\star\# E(\Gamma_{r})}(\lambda)d\lambda
\ \textup{for}\ r=1,\ldots,p/2,
\end{align*}
the limit of $F_{\Gamma}(K,L,j_{1},j_{2})$  can be expressed as
\begin{align}\notag
&\underset{\begin{subarray}{c}j_{1}\rightarrow\infty\\
j_{2}=j_{2}(j_{1})\end{subarray}}{\lim}F_{\Gamma}(K,L,j_{1},j_{2})
\\\notag=&
\underset{\begin{subarray}{c}j_{1}\rightarrow\infty\\
j_{2}=j_{2}(j_{1})\end{subarray}}{\lim}2^{\frac{j_{1}(\beta-1)}{2}p}\ 2^{\frac{j_{2}}{2}p}\ \sigma_{G\star \psi_{j_{1}}}^{p}
\overset{p/2}{\underset{r=1}{\prod}}
\Big[\int_{\mathbb{R}}
e^{i2^{j_{2}}\lambda(t_{k_{d_{1}(\Gamma_{r})}}-t_{k_{d_{2}(\Gamma_{r})}})}
|\hat{\psi}(2^{j_{2}}\lambda)|^{2}
 f_{Y_{j_{1}}}^{\star\# E(\Gamma_{r})}(\lambda)d\lambda\Big]
\\\label{proofthm:weaksmall7s_o3}=&\left[\underset{j_{1}\rightarrow\infty}{\textup{lim}}
2^{\frac{j_{1}\beta}{2}p}\sigma_{G\star \psi_{j_{1}}}^{p}\right]\overset{p/2}{\underset{r=1}{\prod}}\Big[
\int_{\mathbb{R}}
e^{i\lambda(t_{k_{d_{1}(\Gamma_{r})}}-t_{k_{d_{2}(\Gamma_{r})}})}
|\hat{\psi}(\lambda)|^{2}\underset{\begin{subarray}{c}j_{1}\rightarrow\infty\\
j_{2}=j_{2}(j_{1})\end{subarray}}{\lim}2^{-j_{1}}f_{Y_{j_{1}}}^{\star\# E(\Gamma_{r})}(2^{-j_{2}}\lambda)
 d\lambda\Big].
\end{align}
Note that $2\alpha+\beta>1/2$ under Assumption \ref{Assumption:1:wavelet} and $\# E(\Gamma_{r})\geq2$,
(\ref{third_order_conv_decrease}) implies that
\begin{align}
\|f_{Y_{j_{1}}}^{\star\# E(\Gamma_{r})}\|_{\infty}= \|\sigma_{G\star\psi_{j_{1}}}^{-2\# E(\Gamma_{r})}f_{G\star \psi_{j_{1}}}^{\star\# E(\Gamma_{r})}\|_{\infty} \leq \|\sigma_{G\star \psi_{j_{1}}}^{-4}f_{G\star \psi_{j_{1}}}^{\star2}\|_{\infty}.
\end{align}
By Lemma \ref{Lemma:fstarj1j2_limit},
$
\underset{j_{1}\rightarrow\infty}{\textup{lim}}2^{-j_{1}}\|\sigma_{G\star \psi_{j_{1}}}^{-4}f_{G\star \psi_{j_{1}}}^{\star2}\|_{\infty}
$
converges. Hence, we can take the limit inside the integral in
(\ref{proofthm:weaksmall7s_o3}).
By (\ref{lemma:varS}) and Lemma \ref{Lemma:fstarj1j2_limit},
\begin{align}\label{proofthm:weaksmall7s_o3_v2}
\underset{\begin{subarray}{c}j_{1}\rightarrow\infty\\
j_{2}=j_{2}(j_{1})\end{subarray}}{\lim}F_{\Gamma}(K,L,j_{1},j_{2})
=\sigma^{p}
\Big[\overset{p/2}{\underset{r=1}{\prod}}
\gamma_{\# E(\Gamma_{r})}\int_{\mathbb{R}}
e^{i\lambda(t_{k_{d_{1}(\Gamma_{r})}}-t_{k_{d_{2}(\Gamma_{r})}})}
|\hat{\psi}(\lambda)|^{2}
 d\lambda\Big],
\end{align}
where $\sigma^{2}$ is defined in (\ref{lemma:varS})
and
$\gamma_{\# E(\Gamma_{r})}$ is defined in (\ref{proof_third_order_a2}), i.e.,
\begin{align*}
\gamma_{\# E(\Gamma_{r})} = \frac{1}{2\pi} \int_{\mathbb{R}}\left[\int_{\mathbb{R}} \frac{|\hat{\psi}(\eta)|^2}{|\eta|^{1-(2\alpha+\beta)}}d\eta\right]^{-\# E(\Gamma_{r})}
\left[ \int_{\mathbb{R}}e^{it\zeta} \frac{|\hat{\psi}(\zeta)|^{2}}{|\zeta|^{1-(2\alpha+\beta)}}
d\zeta\right]^{\# E(\Gamma_{r})} dt,\ r= 1,2,...,\frac{p}{2}.
\end{align*}
Meanwhile, because $\Gamma$ is a  regular
diagram in $\mathrm{T}(L)$, $C(K,L)$ can be rewritten as follows:
\begin{equation}\label{proofthm:regularcoeff_o3}
C(K,L)=\overset{p/2}{\underset{r=1}{\prod}}
a_{k_{d_{1}(\Gamma_{r})}} a_{k_{d_{2}(\Gamma_{r})}}
\frac{C^{2}_{||,\# E(\Gamma_{r})}}{\# E(\Gamma_{r})!}.
\end{equation}
Combining (\ref{proofthm:weaksmall7s_o3_v2}) and (\ref{proofthm:regularcoeff_o3}) yields
\begin{align}\notag%\label{proofthm:weaksmall8_o3}
&\underset{\begin{subarray}{c}j_{1}\rightarrow\infty\\
j_{2}=j_{2}(j_{1})\end{subarray}}{\lim}
\underset{(K,L)\in D_{p}}{\sum}C(K,L)\underset{\Gamma\in \mathrm{T}^{*}}{\sum}F_{\Gamma}(K,L,j_{1},j_{2})
\\\notag=&\sigma^{p}\underset{(K,L)\in D_{p}}{\sum}
\ \underset{\Gamma\in \mathrm{T}^{*}}{\sum}\Big[
\overset{p/2}{\underset{r=1}{\prod}}
a_{k_{d_{1}(\Gamma_{r})}} a_{k_{d_{2}(\Gamma_{r})}}
\int_{\mathbb{R}}
e^{i\lambda(t_{k_{d_{1}(\Gamma_{r})}}-t_{k_{d_{2}(\Gamma_{r})}})}
|\hat{\psi}(\lambda)|^{2}
d\lambda
\Big]
\Big[
\overset{p/2}{\underset{r=1}{\prod}}
\gamma_{\# E(\Gamma_{r})}
\frac{C^{2}_{||,\# E(\Gamma_{r})}}{\# E(\Gamma_{r})!}
\Big]
\\\label{proofthm:regualr2_o3}=&(p-1)!!
\left[\left(\sigma^{2}\overset{N}{\underset{\ell=2}{\sum}}
\gamma_{\ell}
C^{2}_{||,\ell}\right)
\overset{M}{\underset{k,h=1}{\sum}}a_{k} a_{h}
\int_{\mathbb{R}}
e^{i\lambda(t_{k}-t_{h})}
|\hat{\psi}(\lambda)|^{2}
d\lambda
\right]^{p/2},
\end{align}
where $(p-1)!! = 1\cdot 3\cdot 5 \cdots (p-1).$
By the orthogonal property of the Gaussian random measure $W$ (see (\ref{sample path represent})) and
$\kappa_{N}^{2}=\sigma^{2}\overset{N}{\underset{\ell=2}{\sum}}
\gamma_{\ell}
C^{2}_{||,\ell}$,
the right hand side of (\ref{proofthm:regualr2_o3}) is equal to
\begin{align}\notag
(p-1)!!
\Big[
\mathbb{E}
\Big|
\overset{M}{\underset{k=1}{\sum}}a_{k}
\kappa_{N}\int_{\mathbb{R}}
e^{i\lambda t_{k}}
\hat{\psi}(\lambda)W(d\lambda)
\Big|^{2}
\Big]^{p/2}.
\end{align}
The proof of (b1) is complete.\\

%%%%%%%%%%%%% proof part 2
\noindent{\it Proof of (b2)}:
%$
%\underset{\substack{j_{3}\rightarrow \infty\\j_{2}\rightarrow \infty}}{\textup{lim}}\ \underset{(K,L)\in D_{p}}{\sum}C(K,L)
%\underset{\Gamma\in \mathrm{T}\backslash \mathrm{T}^{*}}{\sum}F_{\Gamma}(K,L,j_{2},j_{3})
%=0.$\\
Because the number of elements in the summation
$$
\underset{(K,L)\in D_{p}}{\sum}C(K,L)
\underset{\Gamma\in \mathrm{T}\backslash \mathrm{T}^{*}}{\sum}F_{\Gamma}(K,L,j_{1},j_{2})
$$
is finite and $C(K,L)$ is independent of $j_{1}$ and $j_{2}$, it  suffices to
show that
$$\underset{\begin{subarray}{c}j_{1}\rightarrow\infty\\
j_{2}=j_{2}(j_{1})\end{subarray}}{\lim}\
F_{\Gamma}(K,L,j_{1},j_{2})=0
$$
for each $\Gamma\in \mathrm{T}(\ell_{1},\ldots,\ell_{p})\backslash
\mathrm{T}^{*}$ and $(K,L)\in D_{p}$.
To prove this,
we rearrange the $p$ levels of $\Gamma$ such that $\ell_{1}\leq \ell_{2}\leq \ell_{3}\leq \cdots \leq \ell_{p}.$
By (\ref{proofthm:weakK(J,L)_o3}) and the Fourier transform relationship between $R_{Y_{j_{1}}}^{\# A_{i,i^{'}}}$
and $f_{Y_{j_{1}}}^{\star\# A_{i,i^{'}}}$,
\begin{align}\notag
&F_{\Gamma}(K,L,j_{1},j_{2})
\\\notag
=&2^{j_{1}\frac{\beta-1}{2}p}\ 2^{j_{2}\frac{p}{2}}\
\sigma_{G\star \psi_{j_{1}}}^{p}\
\int_{\mathbb{R}^{p}}
\Big[\overset{p}{\underset{i=1}{\prod}}2^{-j_{2}}\psi(\frac{2^{j_{2}}t_{k_{i}}-s_{i}}{2^{j_{2}}})\Big]
\Big[
\underset{(i,i^{'}): A_{i,i^{'}}\neq \phi}{\prod}
R_{Y_{j_{1}}}^{\# A_{i,i^{'}}}(s_{i}-s_{i^{'}})
\Big]ds_{1}\cdots ds_{p}
\\\notag
=&
\left[2^{j_{1}\frac{\beta}{2}}
\sigma_{G\star \psi_{j_{1}}}\right]^{p}2^{(j_{2}-j_{1})\frac{p}{2}}
\\\notag
&
\int_{\mathbb{R}^{p}}
\Big[\overset{p}{\underset{i=1}{\prod}}2^{-j_{2}}\psi(\frac{2^{j_{2}}t_{k_{i}}-s_{i}}{2^{j_{2}}})\Big]
\Big[
\underset{(i,i^{'}): A_{i,i^{'}}\neq \phi}{\prod}\hspace{-0.1cm}
\int_{\mathbb{R}}e^{i\lambda_{i,i^{'}}(s_{i}-s_{i^{'}})}f_{Y_{j_{1}}}^{\star\# A_{i,i^{'}}}(\lambda_{i,i^{'}})d\lambda_{i,i^{'}}
\Big]ds_{1}\cdots ds_{p}
\\\notag
=&\left[2^{j_{1}\frac{\beta}{2}}
\sigma_{G\star \psi_{j_{1}}}\right]^{p}2^{(j_{2}-j_{1})\frac{p}{2}}2^{-j_{2}\#\{(i,i^{'}): A_{i,i^{'}}\neq \phi\}}\\\notag
&\times\int_{\mathbb{R}^{p}}
\Big[\overset{p}{\underset{i=1}{\prod}}\psi(t_{k_{i}}-s_{i})\Big]
\Big[
\underset{(i,i^{'}): A_{i,i^{'}}\neq \phi}{\prod}
\int_{\mathbb{R}}e^{i\lambda_{i,i^{'}}(s_{i}-s_{i^{'}})}f_{Y_{j_{1},j_{2}}}^{\star\# A_{i,i^{'}}}(2^{-j_{2}}\lambda_{i,i^{'}})d\lambda_{i,i^{'}}
\Big]ds_{1}\cdots ds_{p}
\\\notag=&
\left[2^{j_{1}\frac{\beta}{2}}
\sigma_{G\star \psi_{j_{1}}}\right]^{p}
\int_{\mathbb{R}^{p}}
\Big[\overset{p}{\underset{i=1}{\prod}}\psi(t_{k_{i}}-s_{i})\Big]
\int_{\mathbb{R}}\cdots \int_{\mathbb{R}} \left[\overset{p}{\underset{i=1}{\prod}}
e^{i s_{i}\Theta_{i}}\right]H(j_{1},j_{2})
\left[\underset{(i,i^{'}): A_{i,i^{'}}\neq \phi}{\prod}d\lambda_{i,i^{'}}\right]
ds_{1}\cdots ds_{p}
\\\label{third_order_F}=&\left[2^{j_{1}\frac{\beta}{2}}
\sigma_{G\star \psi_{j_{1}}}\right]^{p}
\int_{\mathbb{R}}\cdots \int_{\mathbb{R}}\
\left[\overset{p}{\underset{i=1}{\prod}}e^{i t_{k_{i}} \Theta_{i}}\hat{\psi}\left(\Theta_{i}\right)\right]
H(j_{1},j_{2})\left[\underset{(i,i^{'}): A_{i,i^{'}}\neq \phi}{\prod}d\lambda_{i,i^{'}}\right].
%\\\label{proofthm:nonregular1}\rightarrow &
%\int_{\mathbb{R}^{p}}
%\Big[\overset{p}{\underset{i=1}{\prod}}\psi(t_{k_{i}}-s_{i})\Big]
%\Big[
%\underset{(i,i^{'}): A_{i,i^{'}}\neq \phi}{\prod}\delta(s_{i}-s_{i^{'}})
%\Big]ds_{1}\cdots ds_{p}
%\left[\underset{(i,i^{'}): A_{i,i^{'}}\neq \phi}{\prod}2\pi f_{Y_{j_{1},j_{2}}}^{\star\# A_{i,i^{'}}}(0)\right],
\end{align}
where
\begin{equation*}
\Theta_{i}=\Big(\underset{i^{'}: A_{i,i^{'}}\neq \phi}{\sum}\lambda_{i,i^{'}}\Big)-\Big(\underset{i^{'}: A_{i^{'},i}\neq \phi}{\sum}\lambda_{i',i}\Big)
\end{equation*}
and
\begin{equation}\notag
H(j_{1},j_{2}) := 2^{(j_{2}-j_{1})\frac{p}{2}}2^{-j_{2}\#\{(i,i^{'}): A_{i,i^{'}}\neq \phi\}}\left[\underset{(i,i^{'}): A_{i,i^{'}}\neq \phi}{\prod}f_{Y_{j_{1}}}^{\star\# A_{i,i^{'}}}(2^{-j_{2}}\lambda_{i,i^{'}})\right].
\end{equation}
Because (\ref{lemma:varS}) implies that
\begin{equation}\label{H0_o3}
\underset{j_{1}\rightarrow\infty}{\lim}\left[2^{j_{1}\frac{\beta}{2}}
\sigma_{G\star \psi_{j_{1}}}\right]^{p} = \sigma^{p}
\end{equation}
we pay attention to the limiting behavior of
$H(j_{1},j_{2})$.
Note that $H(j_{1},j_{2})$ can be rewritten as
\begin{align}\notag
H(j_{1},j_{2}) =& 2^{(j_{2}-j_{1})\frac{p}{2}}
2^{(j_{1}-j_{2})\#\{(i,i^{'}): \#A_{i,i^{'}}\geq 2\}}
2^{(j_{1}-j_{2})(2\alpha+\beta)\#\{(i,i^{'}): \#A_{i,i^{'}}=1\}}
\\\label{H_o3}\times&\left[\underset{(i,i^{'}): \#A_{i,i^{'}}=1}{\prod}2^{j_{2}(2\alpha+\beta-1)-j_{1}(2\alpha+\beta)}f_{Y_{j_{1}}}(2^{-j_{2}}\lambda_{i,i^{'}})\right]
\left[\underset{(i,i^{'}): \#A_{i,i^{'}}\geq2}{\prod}\hspace{-0.5cm}2^{-j_{1}}f_{Y_{j_{1},j_{2}}}^{\star\# A_{i,i^{'}}}(2^{-j_{2}}\lambda_{i,i^{'}})\right].
\end{align}
Because
\begin{align*}
f_{Y_{j_{1}}}(2^{-j_{2}}\lambda) =& \sigma_{G\star \psi_{j_{1}}}^{-2}
\frac{C_{G}(2^{-j_{2}}\lambda)}{|2^{-j_{2}}\lambda|^{1-\beta}}
|C_{\hat{\psi}}(2^{j_{1}-j_{2}}\lambda)|^{2} |2^{j_{1}-j_{2}}\lambda|^{2\alpha}
\\\notag=&
2^{2j_{1}\alpha}2^{j_{2}
(1-2\alpha-\beta)}\sigma_{G\star \psi_{j_{1}}}^{-2}C_{G}(2^{-j_{2}}\lambda)|C_{\hat{\psi}}(2^{j_{1}-j_{2}}\lambda)|^{2}|\lambda|^{2\alpha+\beta-1}
\end{align*}
and $\underset{j_{1}\rightarrow\infty}{\lim} 2^{-j_{1}\beta }\sigma_{G\star \psi_{j_{1}}}^{-2}=\sigma^{-2}$,
\begin{align}\label{H1_o3}
\underset{\begin{subarray}{c}j_{1}\rightarrow\infty\\
j_{2}=j_{2}(j_{1})\end{subarray}}{\lim}
2^{j_{2}(2\alpha+\beta-1)-j_{1}(2\alpha+\beta)}f_{Y_{j_{1}}}(2^{-j_{2}}\lambda) = \sigma^{-2}C_{G}(0)|C_{\hat{\psi}}(0)|^{2}
|\lambda|^{2\alpha+\beta-1}.
\end{align}
By Lemma \ref{Lemma:fstarj1j2_limit}, if $\# A_{i,i^{'}}\geq2$,
\begin{equation}\label{H2_o3}
\underset{\begin{subarray}{c}j_{1}\rightarrow\infty\\
j_{2}=j_{2}(j_{1})\end{subarray}}{\lim}
2^{-j_{1}}f_{Y_{j_{1}}}^{\star\# A_{i,i^{'}}}(2^{-j_{2}}) = \gamma_{\# A_{i,i^{'}}}.
\end{equation}
By (\ref{third_order_F}), (\ref{H0_o3}), (\ref{H_o3}), (\ref{H1_o3}), and (\ref{H2_o3}), we get
\begin{align}\label{asym_F}
F_{\Gamma}(K,L,j_{1},j_{2}) =O\left(2^{(j_{2}-j_{1})\frac{p}{2}} 2^{-(j_{2}-j_{1})\#\{(i,i^{'}): \#A_{i,i^{'}}\geq2\}}2^{-(j_{2}-j_{1})(2\alpha+\beta)\#\{(i,i^{'}): \#A_{i,i^{'}}=1\}}\right)
\end{align}

Note that $j_{2}(j_{1})$ satisfies $\underset{j_{1}\rightarrow\infty}{\lim\inf}\ \frac{j_{2}(j_{1})}{j_{1}}>1$. To prove $F_{\Gamma}(K,L,j_{1},j_{2})\rightarrow 0$
when $j_{1}\rightarrow\infty$,
(\ref{asym_F}) implies that it suffices to show that
\begin{align}\label{inequality_p/2_order3}
(2\alpha+\beta)\#\{(i,i^{'}): \# A_{i,i^{'}}=1\}+\#\{(i,i^{'}): \# A_{i,i^{'}}\geq 2\}>\frac{p}{2}
\end{align}
for all non-regular diagrams. The proof of (\ref{inequality_p/2_order3}) is completely the same as the proof of (102) in \cite{liu2020central},
so we omit it.
The proof of (b2) is complete.
We thus finish the proof of Proposition \ref{thm:third_order_Gaussian}.

\vskip 20 pt
%%%%%%%%%%%%%%%%%%%%%%%%%%%  Proof %%%%%%%%%%%%%%%%%%%%%%%%%%%%%%%%%%%%%%

%
%
%----------- --Proof of Weakly dependence Large scaling--------------------------

\bigskip

%\section*{Reference}
\bibliographystyle{elsarticle-num}
\bibliography{reference}

\bigskip

%\newpage

%\section*{Acknowledgement}
%The authors are indebted to the inspiring lectures of Professor N. N. Leonenko
%at National Taiwan University for the perspective on {\it Random Fields: modeling, inferences, and applications} in 2006.
%The authors also thank the anonymous reviewers for their valuable comments and suggestions to make the paper more precise and readable.
%The preparation of this manuscript was supported in part by NCTS/TPE and the Taiwan
%Ministry of Science and Technology under grant numbers MOST 104-2115-M-006-016-MY2.

\begin{appendix}

%%%%%%%%%%%%%%%%%%%%%%%%%%%%%%%%%%%%%%%%%%%%%%%%%%%%%%%%%%%%%%%%%%%%%%%%
\section{Proof of Lemma \ref{lemma:var_S_T}}\label{sec:proof:var_S_T}
By the definition $S_{j_{1}} = C_{A,1}G\star \psi_{j_{1}}$,
\begin{align*}\notag
2^{j_{1}\beta}\sigma_{S_{j_{1}}}^{2} =& 2^{j_{1}\beta}\mathbb{E}\left[|C_{A,1}G\star \psi_{j_{1}}|^{2}\right]
\\=&2^{j_{1}\beta}C_{A,1}^{2}\int_{\mathbb{R}}f_{G\star \psi_{j_{1}}}(\lambda)d\lambda
\\=&2^{j_{1}\beta}C_{A,1}^{2}\int_{\mathbb{R}}\frac{C_{G}(\lambda)}{|\lambda|^{1-\beta}} |\hat{\psi}(2^{j_{1}}\lambda)|^{2}d\lambda
\\=&C_{A,1}^{2}\int_{\mathbb{R}}\frac{C_{G}(2^{-j_{1}}\lambda)}{|\lambda|^{1-\beta}} |\hat{\psi}(\lambda)|^{2}d\lambda
\\\rightarrow& C_{A,1}^{2}C_{G}(0)\int_{\mathbb{R}}\frac{|\hat{\psi}(\lambda)|^{2}}{|\lambda|^{1-\beta}} d\lambda
\end{align*}
when $j_{1}\rightarrow\infty$. On the other hand, by (\ref{expectionhermite}),
\begin{align*}\notag
\sigma_{T_{j_{1}}}^{2} =& \overset{\infty}{\underset{\ell_{1} =2}{\sum}}\ \overset{\infty}{\underset{\ell_{2} =2}{\sum}}
\frac{C_{A,\ell_{1}}}{\sqrt{\ell_{1}!}}\frac{C_{A,\ell_{2}}}{\sqrt{\ell_{2}!}}
\int_{\mathbb{R}}\int_{\mathbb{R}}\mathbb{E}\left[H_{\ell_{1}}\left(G(\eta_{1})\right)H_{\ell_{2}}\left(G(\eta_{2})\right)
\right]\psi_{j_{1}}(-\eta_{1})\psi_{j_{1}}(-\eta_{2})d\eta_{1}d\eta_{2}
\\=&
\overset{\infty}{\underset{\ell=2}{\sum}}
C^{2}_{A,\ell}
\int_{\mathbb{R}}\int_{\mathbb{R}}R^{\ell}_{G}\left(\eta_{1}-\eta_{2}\right)
\psi_{j_{1}}(-\eta_{1})\psi_{j_{1}}(-\eta_{2})d\eta_{1}d\eta_{2}
\\=&\overset{\infty}{\underset{\ell=2}{\sum}}
C^{2}_{A,\ell}
\int_{\mathbb{R}}f^{\star\ell}_{G}\left(\lambda\right)
|\hat{\psi}_{j_{1}}(\lambda)|^{2}d\lambda.
\end{align*}
If $2\beta<1$, by the property of multi-fold convolutions of spectral densities in Lemma \ref{lemma:convolution},
there exists a bounded and continuous function $B_{2}$ such that
\begin{align}\notag
&\underset{j_{1}\rightarrow\infty}{\lim}2^{2j_{1}\beta}\int_{\mathbb{R}}f^{\star2}_{G}\left(\lambda\right)
|\hat{\psi}_{j_{1}}(\lambda)|^{2}d\lambda
\\\notag=&\underset{j_{1}\rightarrow\infty}{\lim} 2^{2j_{1}\beta}2^{-j_{1}}
\int_{\mathbb{R}}f^{\star2}_{G}\left(2^{-j_{1}}\lambda\right)
|\hat{\psi}(\lambda)|^{2}d\lambda
\\\notag=&\underset{j_{1}\rightarrow\infty}{\lim} 2^{2j_{1}\beta}2^{-j_{1}}
\int_{\mathbb{R}}B_{2}(2^{-j_{1}}\lambda)|2^{-j_{1}}\lambda|^{2\beta-1}
|\hat{\psi}(\lambda)|^{2}d\lambda
\\\label{2b<1_a}=& B_{2}(0)\int_{\mathbb{R}}|\lambda|^{2\beta-1}
|\hat{\psi}(\lambda)|^{2}d\lambda.
\end{align}
For the integers $\ell\geq3$ satisfying $2<\ell \beta\leq1$,
by Lemma \ref{lemma:convolution},
\begin{align}\label{2b<1_b}
\underset{j_{1}\rightarrow\infty}{\lim}2^{2j_{1}\beta}\int_{\mathbb{R}}f^{\star\ell}_{G}\left(\lambda\right)
|\hat{\psi}_{j_{1}}(\lambda)|^{2}d\lambda
=0.
\end{align}
Denote $\ell^{*} = \min\{\ell\mid \ell\beta>1\}$. For all $\ell\geq\ell^{*}$,
by Lemma \ref{lemma:convolution}, there exists a bounded and continuous function $B_{\ell^{*}}$ such that
\begin{align}\label{2b<1_c}
\underset{j_{1}\rightarrow\infty}{\lim}2^{j_{1}}\int_{\mathbb{R}}f^{\star\ell}_{G}\left(\lambda\right)
|\hat{\psi}_{j_{1}}(\lambda)|^{2}d\lambda
\leq \underset{j_{1}\rightarrow\infty}{\lim}
\int_{\mathbb{R}}B_{\ell^{*}}\left(2^{-j_{1}}\lambda\right)
|\hat{\psi}(\lambda)|^{2}d\lambda
=B_{\ell^{*}}\left(0\right)
\|\hat{\psi}\|_{2}^{2}.
\end{align}
By (\ref{2b<1_a})-(\ref{2b<1_c}) and $\overset{\infty}{\underset{\ell=2}{\sum}}
C^{2}_{A,\ell}<\infty$,
we have
\begin{equation*}
\underset{j_{1}\rightarrow\infty}{\lim} 2^{2j_{1}\beta}\sigma_{T_{j_{1}}}^{2}=
C^{2}_{A,2}
B_{2}(0)\int_{\mathbb{R}}|\lambda|^{2\beta-1}
|\hat{\psi}(\lambda)|^{2}d\lambda
\end{equation*}
under the situation $2\beta<1$.

If $2\beta=1$, by Lemma \ref{lemma:convolution},
\begin{align}\notag
&\int_{\mathbb{R}}f^{\star2}_{G}\left(\lambda\right)
|\hat{\psi}_{j_{1}}(\lambda)|^{2}d\lambda
=
2^{-j_{1}}\int_{\mathbb{R}}f^{\star2}_{G}\left(2^{-j_{1}}\lambda\right)
|\hat{\psi}(\lambda)|^{2}d\lambda
\\\notag=& 2^{-j_{1}}
\int_{\mathbb{R}}B_{2}(2^{-j_{1}}\lambda)\ln\left(2+|2^{-j_{1}}\lambda|^{-1}\right)
|\hat{\psi}(\lambda)|^{2}d\lambda.
\end{align}
Because the logarithmic function is a slowly varying function and
\begin{align*}
\underset{j_{1}\rightarrow\infty}{\lim}
\frac{\ln\left(2+|2^{-j_{1}}\lambda|^{-1}\right)}{\ln 2^{j_{1}}}=1,
\end{align*}
\begin{align}\notag
&\underset{j_{1}\rightarrow\infty}{\lim} 2^{j_{1}}\left(\ln2^{j_{1}}\right)^{-1}
\int_{\mathbb{R}}f^{\star2}_{G}\left(\lambda\right)
|\hat{\psi}_{j_{1}}(\lambda)|^{2}d\lambda
=B_{2}(0)
\|\hat{\psi}\|^{2}.
\end{align}
For the integers $\ell\geq3$, because $2\beta=1$ implies that $\ell\beta>1$,
we have
\begin{align*}
\underset{j_{1}\rightarrow\infty}{\lim}2^{j_{1}}\left(\ln2^{j_{1}}\right)^{-1}\int_{\mathbb{R}}f^{\star\ell}_{G}\left(\lambda\right)
|\hat{\psi}_{j_{1}}(\lambda)|^{2}d\lambda
=0.
\end{align*}
Hence, for the case $2\beta=1$,
\begin{equation*}
\underset{j_{1}\rightarrow\infty}{\lim} 2^{j_{1}}\left(\ln2^{j_{1}}\right)^{-1}\sigma_{T_{j_{1}}}^{2}=
B_{2}(0)
\|\hat{\psi}\|^{2}.
\end{equation*}

If $2\beta>1$, by (\ref{2b<1_c}) directly, we have
\begin{equation*}
\underset{j_{1}\rightarrow\infty}{\lim} 2^{j_{1}}\sigma_{T_{j_{1}}}^{2}=
\overset{\infty}{\underset{\ell=2}{\sum}}
C^{2}_{A,\ell}
B_{\ell}(0).
\end{equation*}

%\begin{align}\notag
%\gamma_{N}:=&\underset{j_{2},j_{3}\rightarrow\infty: j_{3}=\delta j_{2}}{\lim}2^{-j_{2}}\sigma_{j_{1},j_{2}}^{-2N} f^{\star %N}_{X\star\psi_{j_{1}}\star\psi_{j_{2}}}(2^{-j_{3}}\eta)
%\\\notag=&\hspace{0.7cm}\underset{j_{2}\rightarrow\infty}{\lim}2^{-j_{2}}\sigma_{j_{1},j_{2}}^{-2N} f^{\star %N}_{X\star\psi_{j_{1}}\star\psi_{j_{2}}}(0)
%\\\label{lemma_f_cov_limit_equiv}
%=&\hspace{0.7cm}\frac{1}{2\pi} \int_{\mathbb{R}}\left[\int_{\mathbb{R}} \frac{|\Psi(\eta)|^2}{|\eta|^{1-(2\alpha+\beta)}}d\eta\right]^{-N}
%\left[ \int_{\mathbb{R}}e^{it\zeta} \frac{|\Psi(\zeta)|^{2}}{|\zeta|^{1-(2\alpha+\beta)}}
%d\zeta\right]^{N} dt,
%\end{align}
%where $\sigma_{j_{1},j_{2}}$ is the standard deviation of $X\star\psi_{j_{1}}\star\psi_{j_{2}}$
%and

\section{Proof of Lemma \ref{lemma:path_dominate}}\label{sec:proof:lemma:io_path}

First of all, for all $\tau\in \mathbb{R}$ and $\nu\in \mathbb{R}$,
we have
\begin{align}\notag
P\left(|S_{j_{1}}(\tau)|<|T_{j_{1}}(\tau)|\right) =& P\left(|S_{j_{1}}(\tau)|<|T_{j_{1}}(\tau)|, |S_{j_{1}}(\tau)|\geq\sigma_{S_{j_{1}}}^{\nu} \right)
\\\notag&+ P\left(|S_{j_{1}}(\tau)|<|T_{j_{1}}(\tau)|,|S_{j_{1}}(\tau)|<\sigma_{S_{j_{1}}}^{\nu} \right)
\\\notag\leq & P\left(|T_{j_{1}}(\tau)|\geq\sigma_{S_{j_{1}}}^{\nu} \right) + P\left( |S_{j_{1}}(\tau)|<\sigma_{S_{j_{1}}}^{\nu} \right)
\\\label{est:ST}\leq & \frac{\mathbb{E}|T_{j_{1}}(\tau)|^{2}}{\sigma_{S_{j_{1}}}^{2\nu}}  + 2\int_{0}^{\sigma_{S_{j_{1}}}^{\nu}}\frac{1}{\sqrt{2\pi \sigma_{S_{j_{1}}}^{2}}}e^{-\frac{r^{2}}{2\sigma_{S_{j_{1}}}^{2}}}dr.
\end{align}
By Lemma \ref{lemma:var_S_T},  there exist a constant $C>0$ and a threshold $J$
such that
\begin{align}\label{est:T}
\frac{\mathbb{E}|T_{j_{1}}(\tau)|^{2}}{\sigma_{S_{j_{1}}}^{2\nu}}
\leq
C2^{\nu\beta j_{1}}\times\left\{\begin{array}{ll}
 2^{-2\beta j_{1}} & \textup{for}\ \beta\in(0,\frac{1}{2}),
\\
 j_{1}2^{-j_{1}} & \textup{for}\ \beta=\frac{1}{2},
\\
 2^{-j_{1}} & \textup{for}\ \beta\in(\frac{1}{2},1),
\end{array}\right.
\end{align}
and
\begin{align}\label{est:S}
2\int_{0}^{\sigma_{S_{j_{1}}}^{\nu}}\frac{1}{\sqrt{2\pi \sigma_{S_{j_{1}}}^{2}}}e^{-\frac{r^{2}}{2\sigma_{S_{j_{1}}}^{2}}}dr
= 2 \int_{0}^{\sigma_{S_{j_{1}}}^{\nu-1}}\frac{1}{\sqrt{2\pi}}e^{-\frac{r^{2}}{2}}dr\leq C 2^{-\beta\frac{\nu-1}{2}j_{1}}
\end{align}
when $j_{1}\geq J$. By choosing $\nu=\frac{5}{3}$ for $\beta\in(0,\frac{1}{2}]$ and $\nu=\frac{1}{3}+\frac{2}{3\beta}$ for $\beta\in(\frac{1}{2},1)$, (\ref{est:ST})-(\ref{est:S}) imply that when $j_{1}\geq J$
\begin{align}\label{est:S+est:T}
P\left(|S_{j_{1}}(\tau)|<|T_{j_{1}}(\tau)|\right) \leq 2C\times\left\{\begin{array}{ll}
2^{-\frac{\beta}{3}j_{1}} &\ \textup{for}\ \beta\in(0,\frac{1}{2}),\\
j_{1}2^{-\frac{1}{6}j_{1}} &\ \textup{for}\ \beta=\frac{1}{2},
\\
2^{-\frac{1-\beta}{3}j_{1}} &\ \textup{for}\ \beta\in(\frac{1}{2},1).
\end{array}\right.
\end{align}
Define $n_{j_{1}}=[2^{\frac{\beta}{6}j_{1}}]$ for $\beta\in(0,\frac{1}{2}]$
and $n_{j_{1}}=[2^{j_{1}\frac{1-\beta}{6}}]$ for $\beta\in(\frac{1}{2},1)$.
For each $j_{1}\geq J$, we select $n_{j_{1}}$ time points $\tau_{j_{1},1},\tau_{j_{1},2},...,\tau_{j_{1},n_{j_{1}}}$
from $\mathbb{R}$ arbitrarily.
%Below we show that
%\begin{align*}
%P\left(|S_{j_{1}}(\tau_{j_{1},k})|<|T_{j_{1}}(\tau_{j_{1},k})|\ \textup{for some}\ k\in\{1,2,...,n\}\ \textup{i.o.}\right)=0.
%\end{align*}
By (\ref{est:S+est:T}) and Boole's inequality,
\begin{align}\notag
&\overset{\infty}{\underset{j_{1}=J}{\sum}}P\left(|S_{j_{1}}(\tau_{j_{1},k})|<|T_{j_{1}}(\tau_{j_{1},k})|\ \textup{for some}\ k\in\{1,...,n_{j_{1}}\}\right)
\\\label{BC_lemma_condition}
\leq& 2 C \times\left\{\begin{array}{ll}
\overset{\infty}{\underset{j_{1}=J}{\sum}}2^{\frac{\beta}{6}j_{1}} 2^{-\frac{\beta}{3}j_{1}}<\infty &\ \textup{for}\ \beta\in(0,\frac{1}{2}),
\\
\overset{\infty}{\underset{j_{1}=J}{\sum}}2^{\frac{1}{12}j_{1}} j_{1}2^{-\frac{1}{6}j_{1}}<\infty &\ \textup{for}\ \beta=\frac{1}{2},
\\
\overset{\infty}{\underset{j_{1}=J}{\sum}}2^{\frac{1-\beta}{6}j_{1}} 2^{-\frac{1-\beta}{3}j_{1}}<\infty &\ \textup{for}\ \beta\in(\frac{1}{2},1).
\end{array}\right.
\end{align}
By (\ref{BC_lemma_condition}) and the Borel-Cantelli lemma, we obtain (\ref{statement_io}).

%%%%%%%%%%%%%%%%%%%%%%%%%%%%%%%%%%%%%%%%%%%%%%%%%%%%%%%%%%%%%%%%%%%%%%%%5

\section{Proof of Lemma \ref{Lemma:fstarj1j2_limit}}\label{sec:proof:Lemma:fstarj1j2_limit}

Because $f_{G\star\psi_{j_{1}}}\in L^{1}$, by \cite[Theorem 8.22]{folland1999real},
\begin{equation}\label{R_N_spec}
R^{\ell}_{G\star\psi_{j_{1}}}(t) = \int_{\mathbb{R}} e^{i\lambda t} f^{\star \ell}_{G\star\psi_{j_{1}}}(\lambda)d\lambda\ \ \textup{for all}\ \ell\in\{2,3,4,...\}.
\end{equation}
Under Assumption \ref{Assumption:1:wavelet}, which implies that $2\alpha+\beta>1/2$ for all $\beta>0$, $f_{G\star\psi_{j_{1}}}$
is square integrable. Hence,
$R_{G\star\psi_{j_{1}}}\in L^{2}$ by Parseval's identity, i.e., $R^{2}_{G\star\psi_{j_{1}}}\in L^{1}$.
Note that $|R_{G\star\psi_{j_{1}}}|\leq \sigma^{2}_{G\star\psi_{j_{1}}}$.
Hence,
\begin{equation}\label{R_N_integrable}
R^{\ell}_{G\star\psi_{j_{1}}}\in L^{1}\ \textup{for all}\ \ell\in\{2,3,4,...\}.
\end{equation}
By (\ref{R_N_spec}), (\ref{R_N_integrable}) and the inverse Fourier transform,
\begin{equation}\label{R_N_inverse}
f_{G\star\psi_{j_{1}}}^{\star \ell}(\eta)
=\frac{1}{2\pi}
\int_{\mathbb{R}}e^{-i\eta t}\left[R_{G\star\psi_{j_{1}}}(t)\right]^{\ell}dt
\end{equation}
for all $\eta\in \mathbb{R}$.
If $\ell$ is even, (\ref{R_N_inverse}) implies that
\begin{equation*}%\label{R_N_inverse_even}
\|f_{G\star\psi_{j_{1}}}^{\star \ell}\|_{\infty}
\leq\frac{1}{2\pi}
\int_{\mathbb{R}}\left[R_{G\star\psi_{j_{1}}}(t)\right]^{\ell}dt
= f_{G\star\psi_{j_{1}}}^{\star \ell}(0).
\end{equation*}
Next, we prove that (\ref{lemma_f_cov_limit_equiv}) holds. By (\ref{R_N_inverse}),
\begin{align}\notag
2^{-j_{1}}\sigma_{G\star\psi_{j_{1}}}^{-2\ell} f^{\star \ell}_{G\star\psi_{j_{1}}}(2^{-j_{2}}\lambda)
=&
2^{-j_{1}}\sigma_{G\star\psi_{j_{1}}}^{-2\ell}\left[\frac{1}{2\pi}
\int_{\mathbb{R}}e^{-i2^{-j_{2}}\lambda t}\left[R_{G\star\psi_{j_{1}}}(t)\right]^{\ell}dt\right]
\\\label{int_R_power_N}=&\frac{1}{2\pi}
\int_{\mathbb{R}}e^{-i2^{-(j_{2}-j_{1})}\lambda t}\left[\sigma_{G\star\psi_{j_{1}}}^{-2}R_{G\star\psi_{j_{1}}}(2^{j_{1}}t)\right]^{\ell}dt.
\end{align}
By the generalized dominated convergence theorem, whose details are presented in \ref{sec:appendix:DCT},
\begin{align}\label{int_R_power_N_limit}
\underset{\begin{subarray}{c} j_{1}\rightarrow\infty\\ j_{2}=j_{2}(j_{1})\end{subarray}}
{\lim}2^{-j_{1}}\sigma_{G\star\psi_{j_{1}}}^{-2\ell} f^{\star \ell}_{G\star\psi_{j_{1}}}(2^{-j_{2}}\lambda)
=\frac{1}{2\pi}
\int_{\mathbb{R}}\underset{j_{1}\rightarrow\infty}{\lim}\left[\sigma_{G\star\psi_{j_{1}}}^{-2}R_{G\star\psi_{j_{1}}}(2^{j_{1}}t)\right]^{\ell}dt.
\end{align}
For the integrand of the integral in the right hand side of (\ref{int_R_power_N_limit}),
\begin{align}\notag
\underset{j_{1}\rightarrow\infty}{\lim}
\sigma_{G\star j_{1}}^{-2}R_{G\star\psi_{j_{1}}}(2^{j_{1}}t)
=&\underset{j_{1}\rightarrow\infty}{\lim} \sigma_{G\star \psi_{j_{1}}}^{-2} \int_{\mathbb{R}}e^{i2^{j_{1}}t\lambda} f_{G}(\lambda) |\hat{\psi}(2^{j_{1}}\lambda)|^{2}d\lambda
\\\notag=&\underset{j_{1}\rightarrow\infty}{\lim} \sigma_{G\star \psi_{j_{1}}}^{-2} 2^{-j_{1}}\int_{\mathbb{R}}e^{it\zeta} f_{G}(2^{-j_{1}}\zeta)
 |\hat{\psi}(\zeta)|^{2}d\zeta
\\\notag=&\underset{j_{1}\rightarrow\infty}{\lim}\sigma_{G\star \psi_{j_{1}}}^{-2} 2^{-j_{1}}\int_{\mathbb{R}}e^{it\zeta} \frac{C_{G}(2^{-j_{1}}\zeta)}{|2^{-j_{1}}\zeta|^{1-\beta}}
|\hat{\psi}(\zeta)|^{2}d\zeta
\\\notag=&\left[\underset{j_{1}\rightarrow\infty}{\lim}\sigma_{G\star \psi{j_{1}}}^{-2} 2^{-j_{1}\beta}\right]
\left[\underset{j_{1}\rightarrow\infty}{\lim} \int_{\mathbb{R}}e^{it\zeta} C_{G}(2^{-j_{1}}\zeta)
\frac{|\hat{\psi}(\zeta)|^{2}}{|\zeta|^{1-\beta}}d\zeta\right]
\\\label{L2_limit_cov_R}=&\left[\int_{\mathbb{R}} \frac{|\hat{\psi}(\eta)|^2}{|\eta|^{1-\beta}}d\eta\right]^{-1}
\left[ \int_{\mathbb{R}}e^{it\zeta} \frac{|\hat{\psi}(\zeta)|^{2}}{|\zeta|^{1-\beta}}
d\zeta\right],
\end{align}
where the last equality follows from the proof of (\ref{lemma:varS}).
Hence,
\begin{align}\label{int_R_power_N_limit_done}
\underset{\begin{subarray}{c} j_{1}\rightarrow\infty\\ j_{2}=j_{2}(j_{1})\end{subarray}}
{\lim}2^{-j_{1}}\sigma_{G\star\psi_{j_{1}}}^{-2\ell} f^{\star \ell}_{G\star\psi_{j_{1}}}(2^{-j_{2}}\lambda)
=\frac{1}{2\pi}\int_{\mathbb{R}}
\left[\int_{\mathbb{R}} \frac{|\hat{\psi}(\eta)|^2}{|\eta|^{1-\beta}}d\eta\right]^{-\ell}
\left[ \int_{\mathbb{R}}e^{it\zeta} \frac{|\hat{\psi}(\zeta)|^{2}}{|\zeta|^{1-\beta}}
d\zeta\right]^{\ell}dt.
\end{align}
Because the right hand side of (\ref{int_R_power_N_limit_done}) does not depend on $\lambda$,
we have
\begin{align}\notag
\gamma_{\ell}=\underset{\begin{subarray}{c} j_{1}\rightarrow\infty\\ j_{2}=j_{2}(j_{1})\end{subarray}}
{\lim}2^{-j_{1}}\sigma_{G\star \psi_{j_{1}}}^{-2\ell} f^{\star \ell}_{G\star\psi_{j_{1}}}(2^{-j_{2}}\lambda)
=
\underset{ j_{1}\rightarrow\infty}{\lim}2^{-j_{1}}\sigma_{G\star \psi_{j_{1}}}^{-2\ell} f^{\star \ell}_{G\star\psi_{j_{1}}}(0)
\end{align}
for all integer $\ell\geq2$.

From (\ref{int_R_power_N}), we know that $\gamma_{\ell}$ can also be expressed as
\begin{align}\notag
\gamma_{\ell}=
\frac{1}{2\pi}
\int_{\mathbb{R}}\underset{ j_{1}\rightarrow\infty}{\lim}\left[\sigma_{G\star\psi_{j_{1}}}^{-2}R_{G\star\psi_{j_{1}}}(2^{j_{1}}t)\right]^{\ell}dt.
\end{align}
The fact that $0\leq\left[\sigma_{G\star \psi_{j_{1}}}^{-2}R_{G\star\psi_{j_{1}}}(2^{j_{1}}t)\right]^{\ell}\leq 1$ for $\ell\in\{2,4,...\}$ implies
$\gamma_{2}\geq \gamma_{4} \geq \cdots \geq 0$.

%%%%%%%%%%%%%%%%%%%%%%%%%%%%%%%%%%%%%%%%%%%%%%%%%%%%%%%%%%%%%%%%%%%

\section{Dominated convergence theorem used in the proof of Lemma \ref{Lemma:fstarj1j2_limit}}\label{sec:appendix:DCT}
Before taking the limit $j_{1}\rightarrow\infty$ inside the integral in (\ref{int_R_power_N}),
we need to check twos conditions for the generalized dominated convergence theorem:
\begin{equation*}
\left\{\begin{array}{ll}
\textup{(i)}\ g_{j_{1}}(t)\ \textup{converges to}\ g(t)\ \textup{as}\ j_{1}\rightarrow\infty\ \textup{almost everywhere};\\
\textup{(ii)}\ \int_{\mathbb{R}} g_{j_{1}}(t) dt \rightarrow \int_{\mathbb{R}} g(t) dt  \ \textup{as}\ j_{1}\rightarrow\infty,
\end{array}\right.
\end{equation*}
where
\begin{equation}\label{def:g_j2}
g_{j_{1}}(t) = \left[\sigma_{G\star j_{1}}^{-2}R_{G\star\psi_{j_{1}}}(2^{j_{1}}t)\right]^{2}
\end{equation}
is the dominant function for the integral in (\ref{int_R_power_N}) and
\begin{equation}\notag
g(t) = \left[\int_{\mathbb{R}} \frac{|\hat{\psi}(\eta)|^2}{|\eta|^{1-\beta}}d\eta\right]^{-2}
\left[ \int_{\mathbb{R}}e^{it\zeta} \frac{|\hat{\psi}(\zeta)|^{2}}{|\zeta|^{1-\beta}}
d\zeta\right]^{2}.
\end{equation}

The proof of (i) is the same as (\ref{L2_limit_cov_R}). It remains to be proved that
(ii) holds.
By the definition of $g_{j_{1}}$ in (\ref{def:g_j2}) and by (\ref{int_R_power_N}) with $\eta=0$ and $\ell=2$, we have
\begin{align}\notag
\int_{\mathbb{R}}g_{j_{1}}(t)dt  =& \int_{\mathbb{R}}\left[\sigma_{G\star \psi_{j_{1}}}^{-2}R_{G\star\psi_{j_{1}}}(2^{j_{1}}t)\right]^{2} dt
\\\notag=& 2\pi \left[2^{-j_{1}}\sigma_{G\star\psi_{j_{1}}}^{-4} f^{\star 2}_{G\star\psi_{j_{1}}}(0)\right]
\\\label{proof_ii_generalized_L}=& 2\pi \left[2^{-j_{1}2\beta}\sigma_{G\star\psi{j_{1}}}^{-4} \right]\left[2^{j_{1}(2\beta-1)}f^{\star 2}_{G\star\psi_{j_{1}}}(0)\right].
\end{align}
By the proof of (\ref{lemma:varS}), we have
\begin{equation}\label{sigma4_limit}
\underset{j_{1}\rightarrow\infty}{\lim} 2^{-j_{1}2\beta}\sigma_{G\star\psi_{j_{1}}}^{-4} = \left[C_{G}(0)
\int_{\mathbb{R}} \frac{|\hat{\psi}(\eta)|^2}{|\eta|^{1-\beta}}d\eta\right]^{-2}.
\end{equation}
On the other hand,
\begin{align}\notag
2^{j_{1}(2\beta-1)}f^{\star 2}_{G\star\psi_{j_{1}}}(0)
=& 2^{j_{1}(2\beta-1)}\int_{\mathbb{R}}f_{G}(\lambda)^2|\hat{\psi}(2^{j_{1}}\lambda)|^4d\lambda
\\\notag=&2^{j_{1}(2\beta-2)}
\int_{\mathbb{R}}f_{G}(2^{-j_{1}}\eta)^2|\hat{\psi}(\eta)|^4d\eta
\\\notag=&
\int_{\mathbb{R}}C_{G}(2^{-j_{1}}\eta)^2  \frac{|\hat{\psi}(\eta)|^4}{|\eta|^{2-2\beta}}d\eta.
\end{align}
Under Assumption \ref{Assumption:1:wavelet},
\begin{align}\label{thm3:proof_f_bound_end}
\underset{ j_{1}\rightarrow\infty}{\lim}2^{j_{1}(2\beta-1)}f^{\star2}_{G\star\psi_{j_{1}}}(0)
=
C_{G}(0)^2 \int_{\mathbb{R}} \frac{|\hat{\psi}(\eta)|^4}{|\eta|^{2-2\beta}}d\eta<\infty.
\end{align}
By substituting (\ref{sigma4_limit}) and (\ref{thm3:proof_f_bound_end}) into (\ref{proof_ii_generalized_L}),
we get
\begin{align}\notag
&\underset{ j_{1}\rightarrow\infty}{\lim}\int_{\mathbb{R}}g_{j_{1}}(t)dt
=
2\pi
\left[\int_{\mathbb{R}} \frac{|\hat{\psi}(\eta)|^2}{|\eta|^{1-\beta}}d\eta\right]^{-2}\left[ \int_{\mathbb{R}} \frac{|\hat{\psi}(\eta)|^4}{|\eta|^{2-2\beta}}d\eta
\right]
=\int_{\mathbb{R}} g(t) dt,
\end{align}
where the last equality follows from Parseval's identity. Hence, the claim (ii) is proved.

%%%%%%%%%%%%%%%%%%%%%%%%%%%%%%%%%%%%%%%%%%%%5

\end{appendix}

\end{document}